\documentclass[11pt,twoside]{article}
\usepackage{setspace}
\usepackage{amsmath,amsthm,amssymb}
\usepackage{color}

\newtheorem{theorem}{Theorem}[section]
\newtheorem{definition}[theorem]{Definition}

\newtheorem{lemma}[theorem]{Lemma}
\newtheorem{remark}[theorem]{Remark}
\newtheorem{proposition}[theorem]{Proposition}

\newcommand{\Hh}{{\mathcal H}}

\newcommand{\Mm}{{\mathcal M}}
\newcommand{\Aa}{{\mathcal A}}

\date{\rule{0pt}{15pt}\\ \today}

\newcommand{\bg}{\begin{equation}}
\newcommand{\ed}{\end{equation}}
\newcommand{\bga}{\begin{eqnarray}}
\newcommand{\eda}{\end{eqnarray}}

\title{ {\bf\Large Large time decay properties of solutions to a viscous Boussinesq system in a half space}}

\author{
\vspace{5pt} Pigong Han$^{1}$ and Maria E. Schonbek$^{2}$ \\
\vspace{5pt}\small{$^{1}$ Academy of Mathematics and Systems Science}\\[-8pt]
\small{Chinese Academy of Sciences, Beijing 100190 P.R. China}\\
\vspace{5pt}\small{$^{2}$ Department of Mathematics}\\[-8pt]
\small{UC Santa Cruz, Santa Cruz, CA 95064, USA}\\
}
\date {}
\begin{document}

\baselineskip=18pt

\maketitle

Abstract\quad\small  We consider the long time  behavior  of weak
and strong solutions of  the $n$-dimensional viscous Boussinesq system  in the half space, with $n\geq3$ . The $L^r(\mathbb{R}^n_+)$-asymptotics  of  strong solutions and  their first three derivatives, with $1\leq r\leq\infty$,  are derived combining   $L^q-L^r$ estimates and  properties of the fractional powers of the Stokes operator. For the $L^\infty-$asymptotics of the second order derivatives the unboundedness of the projection operator $P:\,L^\infty(\mathbb{R}^n_+)\rightarrow L^\infty_\sigma(\mathbb{R}^n_+)$ is dealt by an appropriate  decomposition of the nonlinear term.

\vskip 2mm


\let\thefootnote\relax\footnotetext{Mathematics Subject Classification 35Q35, 76D05}

\section{Introduction and main results}\label{1}

\setcounter{equation}{0}

We consider questions regarding heat transfer for viscous incompressible flows in the half space. The model used is given by a  Boussinesq approximation. As customary the variations of the density in the continuity equation and the local heat source due to
the viscous dissipation are  neglected. The variations of the temperature are dealt by the
addition of a vertical buoyancy force  in the the fluid motion equation.
This yields the following Dirichlet problem:
\begin{equation} \label{1.1}
\left\{
\begin{array}{lll}
\partial_t\theta-k\Delta\theta+(u\cdot\nabla)\theta=0& \mbox{in} & \mathbb{R}^n_+\times (0, \infty),\\
\partial_tu-\nu\Delta u+(u\cdot\nabla)u+\nabla
p=\beta_1\theta e_n& \mbox{in} & \mathbb{R}^n_+\times (0, \infty),\\
\nabla\cdot u=0 & \mbox{in} & \mathbb{R}^n_+\times (0, \infty),\\
\theta(x, t)=u(x, t)=0 & \mbox{on} & \partial\mathbb{R}^n_+\times (0, \infty),\\
u(x, 0)=a,\,\,\,\,\theta(x, 0)=b & \mbox{in} & \mathbb{R}^n_+,\\
\end{array}
\right.
\end{equation}
where $n\geq3$, and $\mathbb{R}^n_+=\{x=(x^\prime,
x_n)\in\mathbb{R}^n\,|\,\,x_n>0 \}$ is the upper-half space of
$\mathbb{R}^n$; $(a, b)$ is given initial data; the velocity $u=u(x,t)=(u_1(x, t), u_2(x, t),\cdots, u_n(x, t) )$  is an $n$-component divergence free  vector field, the scalar function $\theta=\theta(x, t)$ denotes the density or the temperature and $p=p(x, t)$ is the pressure of the fluid. The Reynolds number $k$  takes into account the strength of heat conductivity. Here, $e_n=(0, 0,\cdots, 0, 1)$, and $\beta_1\in\mathbb{R}^1$ is a physical constant. The constant $\nu>0$ and $k>0$ are the viscous and the thermal diffusion coefficient. By rescaling, without loss of generality, we let
$\nu=k=\beta_1=1$.

The Boussinesq system is commonly used to model ocean and atmospheric dynamics (see \cite{Pe}). It arises from the
density-dependent incompressible Navier-Stokes equations though the Boussinesq approximation, a model where
the density dependence is neglected  in all the terms except the one involving  gravity. Recently, due to its connection
to three-dimensional incompressible flows, this system has received considerable attention in the math dynamical milieu.

When the initial density $b$ is identically zero (or constant), the incompressible Navier-Stokes equations are recovered. The existence of global weak solutions in the energy space for Navier-Stokes equations goes back to J. Leray \cite{Ler}, and the uniqueness of these solutions is only known in two  space dimensions. It is also well known that smooth solutions are global in dimension two  while  for higher dimensions the existence of regular solutions  is only known if  the data are small in some appropriate  spaces. See  \cite{Lem} for more detailed discussions.

Decay results for Navier-Stokes flows have been  widely studied,  readers are referred to \cite{Br1, Br2, BrV, Ka, Sc0, Sc1, Sc2, Wi} and the references therein. The results for solutions to the Navier-Stokes equations  are good guide of what to expect for solutions to the Boussinesq system. The global existence of weak solutions with large data, and  strong solutions with  small data has been studied by several authors. See, e.g. \cite{AH1, Ch, DP1, DP2}. Conditional regularity results for weak solutions (of Serrin type) can be found in \cite{CaD}. The smoothness of solutions arising from large axisymmetric data is addressed in \cite{AH2, HR}. Further regularity issues on the solutions have been discussed. \cite{CFD, FaZ}.

The goal of this paper is to study in which way the variations of the temperature affect the asymptotic behavior of the velocity field. Different models are known in the literature under the name of viscous (or dissipative) Boussinesq system.
The asymptotic behavior of viscous Boussinesq systems of different nature have been addressed recently, e.g., in \cite{AbA, ChG}. The results therein are not comparable  with ours.

The large time behavior of solutions to (\ref{1.1}) has  many open questions.
Self-similarity issues and stability results for solutions in critical spaces (with respect to the scaling) are dealt  for instance in \cite{FeV, KP}. By using Fourier transform and a straightforward adaptation of Caffarelli, Kohn and Nirenberg's method in \cite{CKN}, Brandolese and Schonbek \cite{BrS} recently considered the decay properties of weak and strong solutions of system (\ref{1.1}) in the three-dimensional whole space. The methods employed in  \cite{BrS} (Fourier splitting method for example) seem not applicable to the present case. 

$L^2$-decay for weak solutions has been established by Schonbek \cite{Sc0, Sc1}, Wiegner \cite{Wi} for the Navier-Stokes equations in the whole space; by Borchers and Miyakawa \cite{BM} in the half space; by Brandolese and Schonbek \cite{BrS} for the Boussinesq system in the whole space. The following result generalizes the $L^2$-decay in \cite{BrS} from the whole space to half space.
\begin{theorem}\label{th:1} Let $a\in L^2_\sigma(\mathbb{R}^n_+)$ and $b\in L^1(\mathbb{R}^n_+)\bigcap L^2(\mathbb{R}^n_+)$, $n\geq3$. Then (\ref{1.1}) admits a weak solution $(u, \theta)$ satisfying for $t>0$
\begin{equation} \label{1.2}
\|\theta(t)\|_{L^2(\mathbb{R}^n_+)}\leq
C(1+t)^{-\frac{n}{4}},
\end{equation}
and
\begin{equation} \label{1.3}
\|u(t)\|_{L^2(\mathbb{R}^n_+)}\leq
\left\{
\begin{array}{lll}
C(1+t)^\frac{1}{4}&\mbox{if}& n=3,\vspace{2mm}\\
C\log_e(1+t)&\mbox{if}& n=4,\vspace{2mm}\\
C&\mbox{if}& n\geq5.
\end{array}
\right.
\end{equation}
If there exists $\delta>0$ depending on $\|a\|_{L^2(\mathbb{R}^n_+)}$ such that conditions $[a.]$ and $[b.]$  below hold for $n=3$, and for $n\geq 4$  $[b.]$ holds,
\begin{equation} \label{1.4}
[a.]\;\;\|b\|_{L^1(\mathbb{R}^n_+)}\leq\delta;\quad \mbox{and}\; \;[b.] \quad\|x_nb\|_{L^1(\mathbb{R}^n_+)}<\infty,
\end{equation}
 then the weak solution $(u, \theta)$ satisfies
\begin{equation} \label{1.5}
\|\theta(t)\|_{L^2(\mathbb{R}^n_+)}\leq
C(1+t)^{-\frac{n+2}{4}},\,\,\,\forall\,t>0;\quad and \quad \|u(t)\|_{L^2(\mathbb{R}^n_+)}\longrightarrow0\quad as \quad t\longrightarrow+\infty.
\end{equation}
Further if $a\in L^\frac{n}{n-1}(\mathbb{R}^n_+)$, then for any small $\epsilon>0$ and $t>0$
\begin{equation} \label{1.6}
\|u(t)\|_{L^2(\mathbb{R}^n_+)}\leq C_\epsilon(1+t)^{-\frac{n-2}{4}+\epsilon}.
\end{equation}
\end{theorem}
\begin{remark} \label{r:1}The above estimate (\ref{1.2}) for the temperature looks optimal, since the decay agrees with that
of the heat kernel. On the other hand, the optimality of the estimate (\ref{1.3}) for the velocity field
is not so clear.  Given additional assumptions on the data  we expect that  estimates  (\ref{1.2}), (\ref{1.3})  can be improved to (\ref{1.5}) and (\ref{1.6}).
\end{remark}
\begin{remark} \label{r:2}It is not clear whether (\ref{1.6}) holds for $\epsilon=0$,  the main problem arises from the boundary $\partial\mathbb{R}^n_+$. However, if the initial data satisfy further suitable assumptions (small condition in Theorem \ref{th:4} for example), (\ref{1.6}) holds true with $\epsilon=0$ for strong solutions of (\ref{1.1}) (see Theorem \ref{th:5} below).
\end{remark}
\begin{remark} \label{r:d}$L^q$-decay behavior has been obtained for the Navier-Stokes equations in the whole space,  see \cite{Ka}, and  in the half space, see \cite{FM, Han, Han1} . As far as we know, few $L^q$-decay results are known in the half space  for the Boussinesq system. We recall that  for the whole space,  the $L^2$ decay solutions to the Boussinesq system  required a zero mass condition,  \cite{BrS}. For the half space, in this paper  the assumption that  $b\in L^1$, plus some adequate smallness conditions  insure the decay of the solution, without the necessity    of zero initial mass as for the whole space. \end{remark}
\begin{theorem} \label{th:4}\it Let $a\in L^n_\sigma(\mathbb{R}^n_+)$ and $b\in L^\frac{n}{3}(\mathbb{R}^n_+)$, $n\geq3$.
There exists a number $\epsilon_0>0$ sufficiently small  so that if  $\|a\|_{L^n(\mathbb{R}^n_+)}+\|b\|_{L^\frac{n}{3}(\mathbb{R}^n_+)}\leq\epsilon_0$, then  (1.1) admits a strong solution $(u, \theta)$, which satisfies for $t>0$
\begin{equation}\label{}\notag
\|\nabla^k u(t)\|_{L^q(\mathbb{R}^n_+)}\leq C\epsilon_0t^{-\frac{k}{2}-\frac{n}{2}(\frac{1}{n}-\frac{1}{q})} \,\quad for\,\,\,\,\, k=0, 1\,\,\,\,and\,\,\,\,n\leq q\leq\infty;
\end{equation}
\begin{equation}\label{}\notag
\|\nabla^k\theta(t)\|_{L^q(\mathbb{R}^n_+)}\leq C\epsilon_0t^{-\frac{k}{2}-1-\frac{n}{2}(\frac{1}{n}-\frac{1}{q})} \,\quad for\,\,\,\,\, k=0, 1\,\,\,\,and\,\,\,\,n\leq q\leq\infty.
\end{equation}
\end{theorem}
\begin{theorem} \label{th:5} Let $a\in L^n_\sigma(\mathbb{R}^n_+)\bigcap L^\frac{n}{n-1}(\mathbb{R}^n_+)$ and $b\in L^1(\mathbb{R}^n_+)$ $\bigcap L^2(\mathbb{R}^n_+)$ $\bigcap L^\frac{n}{3}(\mathbb{R}^n_+)$, $n\geq3$.  Assume  (\ref{1.4}) holds. Let $(u, \theta)$ be the strong solution of (\ref{1.1}) obtained in Theorem \ref{th:4}. Then for $t>0$

\begin{equation}\label{}\notag
\|\nabla^k u(t)\|_{L^q(\mathbb{R}^n_+)}\leq
Ct^{-\frac{k-1}{2}-\frac{n}{2}(1-\frac{1}{q})}\,\,\,\,for\,\,\,\,\, k=0, 1\,\,\,\,and\,\,\,\,\frac{n}{n-1}\leq  q \leq+\infty;
\end{equation}

\begin{equation}\label{}\notag
\|\nabla^k\theta(t)\|_{L^q(\mathbb{R}^n_+)}\leq Ct^{-\frac{k+1}{2}-\frac{n}{2}(1-\frac{1}{q})}\,\,\,\,\, for\,\,\,\,\, k=0, 1\,\,\,\,and\,\,\,\,1\leq q\leq+\infty.
\end{equation}

\end{theorem}

For the Navier-Stokes equations, Schonbek and Wiegner \cite{SW} established the decay of higher-order norms of the solutions of Navier-Stokes equations. Here due to existence of the boundary $\partial\mathbb{R}^n_+$, we only can show the decay for the second-order derivatives of the velocity to the Boussinesq system. Further, using the properties of the operators of $E(t), F(t)$ (see (2.2) for the definition of $E(t)$, and (3.18) for $F(t)$), we find the time behavior of the temperature for the Boussinesq system.

\begin{theorem}  \label{th:6}Let $a\in L^n_\sigma(\mathbb{R}^n_+)\bigcap L^\frac{n}{n-1}(\mathbb{R}^n_+)$ and $b\in L^1(\mathbb{R}^n_+)$ $\bigcap L^2(\mathbb{R}^n_+)$ $\bigcap L^\frac{n}{3}(\mathbb{R}^n_+)$, $n\geq3$.  Assume  (\ref{1.4}) holds. Let $(u, \theta)$ be the strong solution of (\ref{1.1}) obtained in Theorem \ref{1.4}. Then for $t>0$
\begin{equation}\label{}\notag
\|\nabla^2u(t)\|_{L^r (\mathbb{R}^n_+)}+\|Au(t)\|_{L^r(\mathbb{R}^n_+)}+\|\partial_tu(t)\|_{L^r(\mathbb{R}^n_+)}+\|\nabla p(t)\|_{L^r(\mathbb{R}^n_+)}\leq Ct^{-\frac{1}{2}-\frac{n}{2}(1-\frac{1}{r})}(1+t^{-n+2})
\end{equation}
provided that $\frac{n}{n-1}\leq r<\infty$; where $A$ is the Stokes operator given in the beginning of section 2.
\begin{equation}\label{}\notag
\|\nabla^2\theta(t)\|_{L^r (\mathbb{R}^n_+)}\leq Ct^{-\frac{3}{2}-\frac{n}{2}(1-\frac{1}{r})}(1+t^{-\frac{n-2}{2}})\quad for\quad 1\leq r\leq\infty.
\end{equation}

Further if $y_nb, y_n^2b\in L^1(\mathbb{R}^n_+)$ and $y_nb\in L^\infty(\mathbb{R}^n_+)$, then for any sufficiently small number $\epsilon>0$ and $t\geq1$
\begin{equation}\label{}\notag
\|\nabla^2u(t)\|_{L^\infty (\mathbb{R}^n_+)}\leq C_\epsilon t^{-\frac{n}{2}+\epsilon}.
\end{equation}

\begin{equation}\label{}\notag
\|\nabla^3\theta(t)\|_{L^r (\mathbb{R}^n_+)}\leq Ct^{-2-\frac{n}{2}(1-\frac{1}{r})}\quad for\quad 1\leq r\leq\infty.
\end{equation}
\end{theorem}

\subsection{Notation and definitions}
The following notations will be used.\\
$C^\infty_{0, \sigma}(\mathbb{R}^n_+)= \{  \phi=(\phi_1, \phi_2, \cdots ,\phi_n) \in C^\infty, \nabla\cdot\phi=0, \; \mbox{real vector-valued functions  with compact support in}\;  \mathbb{R}^n_+\}.$\\
$ L^q_\sigma(\mathbb{R}^n_+),1<q<\infty= \;\mbox{closure of}\; C^\infty_{0, \sigma}(\mathbb{R}^n_+) \;\mbox{with respect to}\; \|\cdot\|_{L^q(\mathbb{R}^n_+)}\; (\mbox{usual Lebesgue}\;\; L^q$ norm.)\\
The $L^\infty(\mathbb{R}^n_+)$ norm is given by  $\|u\|_{L^\infty(\mathbb{R}^n_+)}=ess\,\sup\limits_{x\in\mathbb{R}^n_+}|u(x)|$. \\Constants $C$ are  generic  and may change  from line to line.

\begin{definition} $(u, \theta)$ is called a weak solution of (\ref{1.1}) if for all $T>0$,
\begin{equation} \label{}\notag
u\in L^\infty(0, T; L^2_\sigma(\mathbb{R}^n_+))\bigcap L^2(0, T; H^1_0(\mathbb{R}^n_+))\quad  and \quad \theta\in L^\infty(0, T; L^2(\mathbb{R}^n_+))\bigcap L^2(0, T; H^1_0(\mathbb{R}^n_+))
\end{equation}
satisfy
\begin{equation} \label{}\notag
\int_0^T\int_{\mathbb{R}^n_+} (\,-u \,\partial_t\,\phi +\nabla\,
u\,\cdot\,\nabla\,\phi +u\,\cdot\,\nabla \,u\,\cdot
\phi-\theta \,e_n\,\cdot\phi)dx \,dt=\int_{\mathbb{R}^n_+} a\,\phi(0)\,dx
\end{equation}
for all $\phi\in C_0^\infty([0, T);C_{0,
\sigma}^\infty(\mathbb{R}^n_+))$; and
\begin{equation} \label{}\notag
\int_0^T\int_{\mathbb{R}^n_+} (\,-\theta\, \partial_t\,\psi+\nabla
\theta\,\cdot\nabla\,\psi+\,u\cdot\,\nabla \,\theta\,\cdot
\psi) \,dx \,dt=\int_{\mathbb{R}^n_+} \,b\,\psi(0)\,dx
\end{equation}
for all $\psi\in C_0^\infty([0, T);C_0^\infty(\mathbb{R}^n_+))$,
where $(a, b)\in L^2_\sigma(\mathbb{R}^n_+)\times L^2(\mathbb{R}^n_+)$, and $a$ satisfies a $compatibility$ $condition$: the normal component of $a$ equals to zero on $\partial\mathbb{R}^n_+$. Further, a weak solution $(u, \theta)$ of (1.1) is called a strong solution if for $t>0$, $u\in L^\infty(0, T; H^1_{0, \sigma}(\mathbb{R}^n_+))\bigcap L^2(0, T; H^2(\mathbb{R}^n_+))$  and $\theta\in L^\infty(0, T; H^1_0(\mathbb{R}^n_+))\bigcap L^2(0, T; H^2(\mathbb{R}^n_+))$.
\end{definition}
Using the Navier-Stokes regularity criteria and the standard parabolic regularity theory, we find that a strong solution of (\ref{1.1}) in fact is classical, see \cite{CaD} for example.

\begin{remark} Actually, it is not difficult to verify that one weak solution of (\ref{1.1}) is equivalent to a mild solution $(u, \theta)$, which is defined as: $u\in L^\infty(0, T; L^2_\sigma(\mathbb{R}^n_+))\bigcap L^2(0, T; H^1_0(\mathbb{R}^n_+))$ and $\theta\in L^\infty(0, T; L^2(\mathbb{R}^n_+))\bigcap L^2(0, T; H^1_0(\mathbb{R}^n_+))$ satisfy for all $T>0$ and $0\leq t<T$
\begin{equation}\notag
\left\{
\begin{array}{ll}
\displaystyle(\theta(t), \psi)_{L^2}=(E(t)b, \psi)_{L^2}-\int_0^t(E(t-s)u(s)\cdot\nabla\theta(s), \psi)_{L^2}ds,\vspace{2mm}\\
\displaystyle (u(t), \phi)_{L^2}=(e^{-tA}a, \phi)_{L^2}-\int_0^t(e^{-(t-s)A}P(u(s)\cdot\nabla u(s)-\theta(s) e_n), \phi)_{L^2}ds,
\end{array}
\right.
\end{equation}
\end{remark}
where $(\phi, \psi)$ is from the Definition 1.7, the operators $A, E(t)$ are given in section 2 below.
\section{Decay rates for the weak and strong solutions }\label{2}
\setcounter{equation}{0}
Let $A=-P\Delta$:$\,D(A)\rightarrow L^2_\sigma(\mathbb{R}^n_+)$ be the Stokes
operator for $\mathbb{R}^n_+$, where $P: L^2(\mathbb{R}^n_+)\rightarrow
L^2_\sigma(\mathbb{R}^n_+)$ is the Helmholtz projection operator. Then $A$
is positive self-adjoint operator with dense domain $D(A)\subseteq
L^2_\sigma(\mathbb{R}^n_+)$. There exists a uniquely determined resolution of  the identity
$\{E_\lambda|\,\,\lambda\geq0\}$ in $L^2_\sigma(\mathbb{R}^n_+)$ such that the positive self-adjoint operator $A^\alpha$ $(0<\alpha\leq1)$ is defined by (see \cite{So}):
\begin{equation} \label{}\notag
A^\alpha=\int_0^\infty \lambda^\alpha
dE_\lambda\,\,\,\,with\,\,\,\,domain\,\,\,\, D(A^\alpha)=\big\{v\in
L^2_\sigma(\mathbb{R}^n_+)|\,\, \int_0^\infty \lambda^{2\alpha}d\| E_\lambda
v\|^2_{L^2(\mathbb{R}^n_+)}<\infty\big\}.
\end{equation}
\begin{lemma}\label{l:2.1}\cite{DHP, FM}. For any $f\in L^q_\sigma(\mathbb{R}^n_+)$,
\begin{equation} \label{}\notag
\|\nabla^ke^{-tA}f\|_{L^r(\mathbb{R}^n_+)}\leq
Ct^{-\frac{k}{2}-\frac{n}{2}(\frac{1}{q}-\frac{1}{r})}\|f\|_{L^q(\mathbb{R}^n_+)}
\end{equation}
with $k=0, 1, \cdots,$ provided that $1\leq q<r\leq\infty$ or
$1<q\leq r \leq  \infty$.
\end{lemma}
{\bf Proof of Theorem \ref{th:1}}\quad Let $a\in L^2_\sigma(\mathbb{R}^n_+)$ and $b\in L^1(\mathbb{R}^n_+)\bigcap L^2(\mathbb{R}^n_+)$, $n\geq3$. We consider the successive approximation for $0\leq t<\infty$:
\begin{equation} \label{2.1}
\left\{
\begin{array}{ll}
\displaystyle \theta_0(t)=E(t)b,
\qquad u_0(t)=e^{-tA}a,\vspace{2mm}\\
\displaystyle
\theta_{j+1}(t)=\theta_0(t)-\int_0^tE(t-s)u_j(s)\cdot\nabla\theta_{j+1}(s)ds,\vspace{2mm}\\
\displaystyle
u_{j+1}(t)=u_0(t)-\int_0^te^{-(t-s)A}P\big(u_j(s)\cdot\nabla
u_{j+1}(s)-\theta_j(s)e_n\big)ds
\end{array}
\right.
\end{equation}
for $j=0, 1,2\cdots$. Here the operator $E(t)$ is defined by
\begin{equation} \label{2.2}
E(t)f(x)=\int_{\mathbb{R}^n_+}[G_t(x^\prime-y^\prime,
x_n-y_n)-G_t(x^\prime-y^\prime, x_n+y_n)]f(y)dy,
\end{equation}
and
$G_t(x)=(4\pi t)^{-\frac{n}{2}}e^{-\frac{|x|^2}{4t}}$ is the Gauss kernel.
Problem (\ref{2.1}) admits a unique strong solution $(\theta_{j+1}, u_{j+1})$ (see \cite{BM}), which means
$$
u_{j+1}\in L^\infty(0, T; H^1_{0, \sigma}(\mathbb{R}^n_+)) \bigcap  L^2(0, T; H^2(\mathbb{R}^n_+))
,\;\;\;\;\theta_{j+1}\in L^\infty(0, T; H^1_0(\mathbb{R}^n_+))\bigcap L^2(0, T; H^2(\mathbb{R}^n_+)).
$$
A simple calculation yields for $j=0, 1,2\cdots$ and $t>0$
\begin{equation} \label{2.3}
\left\{
\begin{array}{lll}
\partial_t\theta_{j+1}-\Delta\theta_{j+1}+(u_j\cdot\nabla)\theta_{j+1}=0& \mbox{in} & \mathbb{R}^n_+\times (0, \infty),\\
\partial_tu_{j+1}+Au_{j+1}+P(u_j\cdot\nabla)u_{j+1}=P\theta_j e_n& \mbox{in} & \mathbb{R}^n_+\times (0, \infty),\\
\nabla\cdot u_{j+1}=0 & \mbox{in} & \mathbb{R}^n_+\times (0, \infty),\\
\theta_{j+1}(x, t)=u_{j+1}(x, t)=0 & \mbox{on} & \partial\mathbb{R}^n_+\times (0, \infty),\\
u_{j+1}(x, 0)=a,\,\,\,\,\theta_{j+1}(x, 0)=b & \mbox{in} & \mathbb{R}^n_+,\\
\end{array}
\right.
\end{equation}
If $b\in L^q(\mathbb{R}^n_+)$ with $1\leq q\leq\infty$, then for $t>0$
\begin{equation} \label{2.4}
\sup\limits_{j\geq0}\|\theta_{j+1}(t)\|_{L^q(\mathbb{R}^n_+)}\leq
\|b\|_{L^1(\mathbb{R}^n_+)}\big(\Aa+\frac{Ct}{q}\big)^{-\frac{n}{2}(1-\frac{1}{q})},
\end{equation}
where $\Aa=\Aa(q, \|b\|_{L^1(\mathbb{R}^n_+)}, \|b\|_{L^q(\mathbb{R}^n_+)})$ and $C>0$ is an absolute constant. The proof of (\ref{2.4}) is same to that of Lemma 3.2 in \cite{BrS}, and as such  details are omitted. Multiply the second equation of (\ref{2.3})  by $u_{j+1}$  integrate in space, then  for $j=0,1,2\cdots$ and $t>0$
\begin{equation} \label{2.5}
\frac{d}{dt}\|u_{j+1}(t)\|^2_{L^2(\mathbb{R}^n_+)}+2\|\nabla u_{j+1}(t)\|^2_{L^2(\mathbb{R}^n_+)}=2\int_{\mathbb{R}^n_+}\theta_j(x, t)e_n\cdot u_{j+1}(x, t)dx
\end{equation}
Multiply the first  equation of (\ref{2.3})  by $\theta_{j+1}$  integrate in space, then  for $j=0,1,2\cdots$ and $t>0$
\begin{equation} \label{2.6}
\frac{d}{dt}\|\theta_{j+1}(t)\|^2_{L^2(\mathbb{R}^n_+)}+2\|\nabla \theta_{j+1}(t)\|^2_{L^2(\mathbb{R}^n_+)}=0.
\end{equation}
The following auxiliary estimate is needed first: for any $\rho, t>0$ and $j=0,1,2\cdots$
\begin{align} \label{2.7}
\|\nabla u_{j+1}(t)\|^2_{L^2(\mathbb{R}^n_+)}= \|A^\frac{1}{2}
u_{j+1}(t)\|^2_{L^2(\mathbb{R}^n_+)}
= \int_0^\infty \lambda d\|E_\lambda u_{j+1}(t)\|^2_{L^2(\mathbb{R}^n_+)}\vspace{2mm}\\
\geq \rho\int_\rho^\infty d\|E_\lambda u_{j+1}(t)\|^2_{L^2(\mathbb{R}^n_+)}\vspace{2mm}
= \rho(\|u_{j+1}(t)\|^2_{L^2(\mathbb{R}^n_+)}-\|E_\rho u_{j+1}(t)\|^2_{L^2(\mathbb{R}^n_+)}).\notag
\end{align}
Combining  (\ref{2.7}) with  (\ref{2.5}), yields for all $\rho, t>0$ and $j=0,1,2\cdots$
\begin{equation} \label{2.8}
\frac{d}{dt}\|u_{j+1}(t)\|^2_{L^2(\mathbb{R}^n_+)} +\rho\|u_{j+1}(t)\|^2_{L^2(\mathbb{R}^n_+)}\leq2\rho\|E_\rho
u_{j+1}(t)\|^2_{L^2(\mathbb{R}^n_+)}+2\|u_{j+1}(t)\|_{L^2(\mathbb{R}^n_+)} \|\theta_j(t)\|_{L^2(\mathbb{R}^n_+)}.
\end{equation}
Recall that  (see \cite{BM})
\begin{equation} \label{2.9}
\|E_\lambda P(w\cdot\nabla)u\|_{L^2(\mathbb{R}^n_+)}\leq C\lambda^\frac{n+2}{4}\|w\|_{L^2(\mathbb{R}^n_+)}\|u\|_{L^2(\mathbb{R}^n_+)},
\end{equation}
for all $\lambda>0$ and $w, u\in H^1_{0, \sigma}(\mathbb{R}^n_+)$, where $C$ is independent of $w, u$ and
$\lambda$. Combining (\ref{2.9}), (\ref{2.1}) and Lemma \ref{l:2.1}, gives for any $\rho>0$ and $t>0$, $1<r<2$
\begin{equation} \label{2.10}
\begin{array}{rcl}
&&\displaystyle \|E_\rho u_{j+1}(t)\|_{L^2(\mathbb{R}^n_+)}\vspace{2mm}\\
&\leq& \displaystyle\|E_\rho
e^{-tA}a\|_{L^2(\mathbb{R}^n_+)}+\int_0^t\|E_\rho e^{-(t-s)A}P\big(e_n\theta_j(s)\big)ds\|_{L^2(\mathbb{R}^n_+)}\vspace{2mm}\\
&& \displaystyle+\big\|E_\rho \int_0^t\,\,\, \big(\int_0^\rho+\int_\rho^\infty\!\!\!\big)e^{-(t-s)\lambda}d\big[
E_\lambda P\big((u_j(s)\cdot \,\,\nabla)u_{j+1}(s)\big)\big]ds \big\|_{L^2(\mathbb{R}^n_+)}\vspace{2mm}\\
&\leq& \displaystyle\|e^{-tA}a\|_{L^2(\mathbb{R}^n_+)}+\big\|\int_0^te^{-(t-s)\rho}E_\rho
P((u_j(s)\cdot \nabla)u_{j+1}(s))ds \big\|_{L^2(\mathbb{R}^n_+)}\vspace{2mm}\\
&& +\displaystyle\big\|\int_0^t (t-s)\big\{\int_0^\rho e^{-(t-s)\lambda}[E_\lambda
P((u_j(s)\cdot \nabla)u_{j+1}(s))]d\lambda \big\}ds \big\|_{L^2(\mathbb{R}^n_+)}\vspace{2mm}\\
&&\displaystyle+C\int_0^\frac{t}{2}(t-s)^{-\frac{n}{2}(\frac{1}{r}-\frac{1}{2})}\|\theta_j(s)\|_{L^r(\mathbb{R}^n_+)}ds +C\int^t_\frac{t}{2}\|\theta_j(s)\|_{L^2(\mathbb{R}^n_+)}ds\vspace{2mm}\\
&\leq& \displaystyle
\|e^{-tA}a\|_{L^2(\mathbb{R}^n_+)}+C\rho^\frac{n+2}{4}\int_0^t
\|u_j(s)\|_{L^2(\mathbb{R}^n_+)}\|u_{j+1}(s)\|_{L^2(\mathbb{R}^n_+)}ds\vspace{2mm}\\
&&\displaystyle+C\int_0^\frac{t}{2}(t-s)^{-\frac{n}{2}(\frac{1}{r}-\frac{1}{2})}\|\theta_j(s)\|_{L^r(\mathbb{R}^n_+)}ds +C\int^t_\frac{t}{2}\|\theta_j(s)\|_{L^2(\mathbb{R}^n_+)}ds.
\end{array}
\end{equation}
Combining (\ref{2.8}) and (\ref{2.10}), we obtain for any $\rho, t>0$ and $j=0, 1,\cdots$
\begin{equation} \label{2.11}
\begin{array}{rcl}
&&\displaystyle\frac{d}{dt}\|u_{j+1}(t)\|^2_{L^2(\mathbb{R}^n_+)} +\rho\|u_{j+1}(t)\|^2_{L^2(\mathbb{R}^n_+)}\vspace{2mm}\\
&\leq&\displaystyle C\rho\big(\|e^{-tA}a\|_{L^2(\mathbb{R}^n_+)}+\rho^\frac{n+2}{4}\int_0^t
\|u_j(s)\|_{L^2(\mathbb{R}^n_+)}\|u_{j+1}(s)\|_{L^2(\mathbb{R}^n_+)}ds\vspace{2mm}\\
&&\displaystyle+\int_0^\frac{t}{2}(t-s)^{-\frac{n}{2}(\frac{1}{r}-\frac{1}{2})}
\|\theta_j(s)\|^{\frac{2}{r}-1}_{L^1(\mathbb{R}^n_+)}\|\theta_j(s)\|^{2(1-\frac{1}{r})}_{L^2(\mathbb{R}^n_+)}ds
+\int^t_\frac{t}{2}\|\theta_j(s)\|_{L^2(\mathbb{R}^n_+)}ds\big)^2\vspace{2mm}\\
&&\displaystyle
+2\|u_{j+1}(t)\|_{L^2(\mathbb{R}^n_+)} \|\theta_j(t)\|_{L^2(\mathbb{R}^n_+)}\end{array}
\end{equation}
A similar  proof as for (\ref{2.11}), gives for $\rho, t>0$ and $j=0, 1,\cdots$
\begin{equation} \label{2.12}
\begin{array}{rcl}
&&\displaystyle\frac{d}{dt}\|\theta_{j+1}(t)\|^2_{L^2(\mathbb{R}^n_+)} +\rho\|\theta_{j+1}(t)\|^2_{L^2(\mathbb{R}^n_+)}\vspace{2mm}\\
&\leq&\displaystyle C\rho\big(\|E(t)b\|_{L^2(\mathbb{R}^n_+)}+\rho^\frac{n+2}{4}\int_0^t
\|u_j(s)\|_{L^2(\mathbb{R}^n_+)}\|\theta_{j+1}(s)\|_{L^2(\mathbb{R}^n_+)}ds\big)^2.
\end{array}
\end{equation}
Inequality  (\ref{2.4}) yields
$\sup\limits_{j\geq0}\|\theta_{j+1}(t)\|_{L^2(\mathbb{R}^n_+)}\leq C\|b\|_{L^1(\mathbb{R}^n_+)}(1+t)^{-\frac{n}{4}}.
$ for $t>0$.
From the explicit formulation of the operator $E(t)$, we infer for $t>0$
$$
\|\theta_0(t)\|_{L^2(\mathbb{R}^n_+)}=\|E(t)b\|_{L^2(\mathbb{R}^n_+)}\leq C\|b\|_{L^1(\mathbb{R}^n_+)}(1+t)^{-\frac{n}{4}}.
$$
Thus for all $t>0$
\begin{equation} \label{2.13}
\sup\limits_{j\geq0}\|\theta_j(t)\|_{L^2(\mathbb{R}^n_+)}\leq C\|b\|_{L^1(\mathbb{R}^n_+)}(1+t)^{-\frac{n}{4}}.
\end{equation}
It follows from (\ref{2.5}) that for $t>0$, $j=0, 1,\cdots$
$$
2\|u_{j+1}(t)\|_{L^2(\mathbb{R}^n_+)}\frac{d}{dt}\|u_{j+1}(t)\|_{L^2(\mathbb{R}^n_+)}+2\|\nabla u_{j+1}(t)\|^2_{L^2(\mathbb{R}^n_+)}\leq2\|u_{j+1}(t)\|_{L^2(\mathbb{R}^n_+)}
\|\theta_j(t)\|_{L^2(\mathbb{R}^n_+)},
$$
which implies
\begin{equation} \label{2.14}
\frac{d}{dt}\|u_{j+1}(t)\|_{L^2(\mathbb{R}^n_+)}\leq\|\theta_j(t)\|_{L^2(\mathbb{R}^n_+)},
\end{equation}
Note that Lemma 2.1 implies for $t>0$
\begin{equation} \label{2.15}
\|u_0(t)\|_{L^2(\mathbb{R}^n_+)}=\|e^{-tA}a\|_{L^2(\mathbb{R}^n_+)}\leq \|a\|_{L^2(\mathbb{R}^n_+)}.
\end{equation}
Combining  (\ref{2.13})-(\ref{2.15}), yields for $t>0$
\begin{equation} \label{2.16}
\sup\limits_{j\geq0}\|u_j(t)\|_{L^2(\mathbb{R}^n_+)}\leq
\left\{
\begin{array}{lll}
C(1+t)^\frac{1}{4}&\mbox{if}& n=3,\vspace{2mm}\\
C\log_e(1+t)&\mbox{if}& n=4,\vspace{2mm}\\
C&\mbox{if}& n\geq5,
\end{array}
\right.
\end{equation}
where $C$ depends on $n, \;\|a\|_{L^2(\mathbb{R}^n_+)},\;\|b\|_{L^1(\mathbb{R}^n_+)}$.
To deal with the second part of the theorem, assume that (\ref{1.4}) holds. Firstly we show that for $n=3$ and any $t>0$
\begin{equation} \label{2.17}
\sup\limits_{j\geq0}\|u_j(t)\|_{L^2(\mathbb{R}^n_+)}\leq\|a\|_{L^2(\mathbb{R}^n_+)}+\mathcal{M}(1+t)^\frac{1}{8},
\end{equation}
where $\mathcal{M}>0$ is a constant independent of $j$ to be determined below.

Using the explicit formula (2.2), we get for $t>0$
$$
E(t)b(x)=\int_{\mathbb{R}^n_+}\int_{-1}^1\partial_{x_n}G_t(x^\prime-y^\prime, x_n-sy_n)(-y_n)b(y)dsdy.
$$
So we have the following estimate for $t>0$
\begin{equation} \label{2.18}
\|\theta_0(t)\|_{L^2(\mathbb{R}^n_+)}=\|E(t)b\|_{L^2(\mathbb{R}^n_+)}\leq \|\nabla G_t\|_{L^2(\mathbb{R}^n)}\|x_nb\|_{L^1(\mathbb{R}^n_+)}\leq  C(1+t)^{-\frac{n+2}{4}}\|x_nb\|_{L^1(\mathbb{R}^n_+)}.
\end{equation}
\begin{remark} We first give the details for the case $n=3$, before $n\geq4$. As we will see it, the computations for $n\geq4$ are easier and simpler since the bounds for $\|u\|_{L^2(\mathbb{R}^n_+)}$ are better.  It should be pointed out that the smallness condition in (\ref{1.4}) is unnecessary for the case $n\geq4$, which can be found easily in the following proofs on case $n\geq4$.
\end{remark}
\noindent {\bf Case $n=3$.} By (\ref{2.15}), there exists a $j_0\geq0$ such that for $t>0$
\begin{equation} \label{2.19}
\sup\limits_{0\leq i\leq j_0}\|u_i(t)\|_{L^2(\mathbb{R}^n_+)}\leq\|a\|_{L^2(\mathbb{R}^n_+)}+\mathcal{M}(1+t)^\frac{1}{8}.
\end{equation}
For the rest of the proof of Theorem \ref{th:1},  set $\label{rho} \rho=k(t+1)^{-1}$,
where the integer $k$ will be  appropriately determined below. Multiplying (\ref{2.12}) by $(t+1)^k$ and combining with (\ref{2.13}), (\ref{2.18}), (\ref{2.19}), yields  all $i\in[0, j_0]$  and $t>0$
\begin{equation} \label{}{\notag}
\begin{array}{rcl}
&&\displaystyle \frac{d}{dt}\big((t+1)^k\|\theta_{i+1}(t)\|^2_{L^2(\mathbb{R}^n_+)}\big)\vspace{2mm}\\
&\leq&\displaystyle C(t+1)^{k-1}\big(\|\theta_0(t)\|_{L^2(\mathbb{R}^n_+)}
+C(t+1)^{-\frac{5}{4}}\int_0^t\|u_i(s)\|_{L^2(\mathbb{R}^n_+)}\|\theta_{i+1}(s)\|_{L^2(\mathbb{R}^n_+)}ds\big)^2\vspace{2mm}\\

\end{array}
\end{equation}

\begin{equation} \label{2.20}
\begin{array}{rcl}

&\leq&\displaystyle C(t+1)^{k-1}\big((1+t)^{-\frac{5}{4}}\|x_nb\|_{L^1(\mathbb{R}^n_+)}\vspace{2mm}\\
&&\displaystyle
+\|b\|_{L^1(\mathbb{R}^n_+)}(t+1)^{-\frac{5}{4}}\int_0^t(\|a\|_{L^2(\mathbb{R}^n_+)}+\mathcal{M}(1+s)^\frac{1}{8})(1+s)^{-\frac{3}{4}}ds\big)^2
\vspace{2mm}\\

&\leq&\displaystyle C(t+1)^{k-1-\frac{5}{2}}\big(\|x_nb\|_{L^1(\mathbb{R}^n_+)}
+\|b\|_{L^1(\mathbb{R}^n_+)}(\|a\|_{L^2(\mathbb{R}^n_+)}(1+t)^\frac{1}{4}+\mathcal{M}(1+t)^\frac{3}{8})\big)^2\vspace{2mm}\\

&\leq&\displaystyle C(t+1)^{k-1-\frac{7}{4}}\big(\|x_nb\|_{L^1(\mathbb{R}^n_+)}+\|b\|_{L^1(\mathbb{R}^n_+)}(\|a\|_{L^2(\mathbb{R}^n_+)}+\mathcal{M})\big)^2.
\end{array}
\end{equation}
By taking $k>0$ suitably large in (\ref{2.20}), we infer for  $i\in[0, j_0]$ and $t>0$
\begin{equation} \label{}\notag
\|\theta_{i+1}(t)\|_{L^2(\mathbb{R}^n_+)}\leq C\big(\|b\|_{L^2(\mathbb{R}^n_+)}+\|x_nb\|_{L^1(\mathbb{R}^n_+)}
+\|b\|_{L^1(\mathbb{R}^n_+)}(\|a\|_{L^2(\mathbb{R}^n_+)}+\mathcal{M})\big)(1+t)^{-\frac{7}{8}}.
\end{equation}
Note that for $t>0$,
$\|\theta_0(t)\|_{L^2(\mathbb{R}^n_+)}\leq C\|b\|_{L^1(\mathbb{R}^n_+)}(1+t)^{-\frac{n}{4}}\leq\|a\|_{L^2(\mathbb{R}^n_+)}$
provided $\|b\|_{L^1(\mathbb{R}^n_+)}\leq\delta$,  with $\delta>0$ is suitably small.
Let  $i\in[0, j_0]$ with given $j_0\geq0$, and $t>0$, then
\begin{equation} \label{2.21}
\|\theta_i(t)\|_{L^2(\mathbb{R}^n_+)}\leq C\big(\|b\|_{L^2(\mathbb{R}^n_+)}+\|x_nb\|_{L^1(\mathbb{R}^n_+)}
+\|b\|_{L^1(\mathbb{R}^n_+)}(\|a\|_{L^2(\mathbb{R}^n_+)}+\mathcal{M})\big)(1+t)^{-\frac{7}{8}}.
\end{equation}
Using (\ref{2.14}) and (\ref{2.21}), we get for $n=3$, any $i\in[0, j_0]$  and $t>0$
\begin{equation} \label{2.22}
\begin{array}{rcl}
&&\displaystyle
\|u_{i+1}(t)\|_{L^2(\mathbb{R}^n_+)}\vspace{2mm}\\
&\leq&\displaystyle  \|a\|_{L^2(\mathbb{R}^n_+)}+C_0(\|b\|_{L^2(\mathbb{R}^n_+)}+\|x_nb\|_{L^1(\mathbb{R}^n_+)}
+\|b\|_{L^1(\mathbb{R}^n_+)}(\|a\|_{L^2(\mathbb{R}^n_+)}+\mathcal{M}))(1+t)^\frac{1}{8}.
\end{array}
\end{equation}
In (\ref{2.22}), set
$$C_0(\|b\|_{L^2(\mathbb{R}^n_+)}+\|x_nb\|_{L^1(\mathbb{R}^n_+)}
+\|b\|_{L^1(\mathbb{R}^n_+)}(\|a\|_{L^2(\mathbb{R}^n_+)}+\mathcal{M}))\leq\mathcal{M}.
$$
That is, choose $\Mm>0$ sufficiently large and $\delta>0$  sufficiently small with $\|b\|_{L^1(\mathbb{R}^n_+)}\leq\delta$, so that
$(\frac{1}{C_0}-\|b\|_{L^1(\mathbb{R}^n_+)})\mathcal{M}\geq\|b\|_{L^1(\mathbb{R}^n_+)}\|a\|_{L^2(\mathbb{R}^n_+)}+\|x_nb\|_{L^1(\mathbb{R}^n_+)}.
$
Then (\ref{2.17}) follows by induction by  (\ref{2.22}) with the choices of $\Mm$ and $\delta$ given above.
From (\ref{2.17}), and following the steps for proof of (\ref{2.21}), we get
\begin{equation} \label{2.23}
\sup\limits_{j\geq0}\|\theta_j(t)\|_{L^2(\mathbb{R}^n_+)}\leq C(1+t)^{-\frac{7}{8}}.
\end{equation}
Combining (\ref{2.17}) and (\ref{2.23}), the same process that gives (\ref{2.21}), yields for all $j\geq0$, $t>0$
\begin{equation} \label{}\notag
\begin{array}{rcl}
&&\displaystyle \frac{d}{dt}\big((t+1)^k\|\theta_{j+1}(t)\|^2_{L^2(\mathbb{R}^n_+)}\big)\vspace{2mm}\\
&\leq&\displaystyle C(t+1)^{k-1}\big(\|\theta_0(t)\|_{L^2(\mathbb{R}^n_+)}
+C(t+1)^{-\frac{5}{4}}\int_0^t\|u_j(s)\|_{L^2(\mathbb{R}^n_+)}\|\theta_{j+1}(s)\|_{L^2(\mathbb{R}^n_+)}ds\big)^2\vspace{2mm}\\
&\leq&\displaystyle C(t+1)^{k-1}\big((1+t)^{-\frac{5}{4}}\|x_nb\|_{L^1(\mathbb{R}^n_+)}+(t+1)^{-\frac{5}{4}}\int_0^t(1+s)^{-\frac{3}{4}}ds\big)^2
\vspace{2mm}\\
&\leq&\displaystyle C(t+1)^{k-1-2}.
\end{array}
\end{equation}
Thus for $k$ suitably large, and $t>0$, we get
$\sup\limits_{j\geq0}\|\theta_j(t)\|_{L^2(\mathbb{R}^n_+)}\leq C(1+t)^{-1}.
$ This new estimate, combined  with (2.17), and the steps yielding (\ref{2.21}), yield for $t>0$
\begin{equation} \label{}\notag
\begin{array}{rcl}
&&\displaystyle \frac{d}{dt}\big((t+1)^k\|\theta_{j+1}(t)\|^2_{L^2(\mathbb{R}^n_+)}\big)\vspace{2mm}\\
&\leq&\displaystyle C(t+1)^{k-1}\big(\|\theta_0(t)\|_{L^2(\mathbb{R}^n_+)}
+C(t+1)^{-\frac{5}{4}}\int_0^t\|u_j(s)\|_{L^2(\mathbb{R}^n_+)}\|\theta_{j+1}(s)\|_{L^2(\mathbb{R}^n_+)}ds\big)^2\vspace{2mm}\\
&\leq&\displaystyle C(t+1)^{k-1}\big((1+t)^{-\frac{5}{4}}\|x_nb\|_{L^1(\mathbb{R}^n_+)}+(t+1)^{-\frac{5}{4}}\int_0^t(1+s)^{\frac{1}{8}-1}ds\big)^2
\vspace{2mm}\\
&\leq&\displaystyle C(t+1)^{k-1-\frac{9}{4}},
\end{array}
\end{equation}
which implies
$\sup\limits_{j\geq0}\|\theta_j(t)\|_{L^2(\mathbb{R}^n_+)}\leq C(1+t)^{-\frac{9}{8}}$.
This inequality combined with (\ref{2.14}) gives
\begin{equation} \label{2.24}
\sup\limits_{j\geq0}\|u_j(t)\|_{L^2(\mathbb{R}^n_+)}\leq C.
\end{equation}
The new bounds on $\sup\limits_{j\geq0}\|\theta_j(t)\|_{L^2(\mathbb{R}^n_+)}\leq C(1+t)^{-\frac{9}{8}} $ and $\sup\limits_{j\geq0}\|u_j(t)\|_{L^2(\mathbb{R}^n_+)}\leq C$, and the steps that yield (2.21) give
\begin{equation} \label{}\notag
\begin{array}{rcl}
&&\displaystyle \frac{d}{dt}\big((t+1)^k\|\theta_{j+1}(t)\|^2_{L^2(\mathbb{R}^n_+)}\big)\vspace{2mm}\\
&\leq&\displaystyle C(t+1)^{k-1}\big(\|\theta_0(t)\|_{L^2(\mathbb{R}^n_+)}
+C(t+1)^{-\frac{5}{4}}\int_0^t\|u_j(s)\|_{L^2(\mathbb{R}^n_+)}\|\theta_{j+1}(s)\|_{L^2(\mathbb{R}^n_+)}ds\big)^2\vspace{2mm}\\
&\leq&\displaystyle C(t+1)^{k-1}\big((1+t)^{-\frac{5}{4}}\|x_nb\|_{L^1(\mathbb{R}^n_+)}+(t+1)^{-\frac{5}{4}}\int_0^t(1+s)^{-\frac{9}{8}}ds\big)^2
\vspace{2mm}\\
&\leq&\displaystyle C(t+1)^{k-1-\frac{5}{2}},
\end{array}
\end{equation}
from which,
\begin{equation} \label{2.25}
\sup\limits_{j\geq0}\|\theta_j(t)\|_{L^2(\mathbb{R}^n_+)}\leq C(1+t)^{-\frac{5}{4}}.
\end{equation}
This concludes the proof for the case $n=3$.\\
\noindent {\bf Case $n=4$.} Setting $\rho=k(t+1)^{-1}$ with some large positive integer $k$,
multiplying both sides of (\ref{2.12}) by $(t+1)^k$, together with (\ref{2.16}), (\ref{2.18}), we conclude for any $t>0$ and $j=0, 1, 2,\cdots$
$$
\begin{array}{rcl}
&&\displaystyle \frac{d}{dt}\big((t+1)^k\|\theta_{j+1}(t)\|^2_{L^2(\mathbb{R}^n_+)}\big)\vspace{2mm}\\
&\leq&\displaystyle C(t+1)^{k-1}\big((1+t)^{-\frac{3}{2}}
+(t+1)^{-\frac{3}{2}}\int_0^t(1+s)^{-1}\log_e(1+s)ds\big)^2\vspace{2mm}\\
&\leq&\displaystyle C(t+1)^{k-1}\big((1+t)^{-\frac{3}{2}}+(t+1)^{-\frac{3}{2}}\log_e^2(1+t)\big)^2,
\end{array}
$$
which implies for any $t>0$
$$
\sup\limits_{j\geq0}\|\theta_{j+1}(t)\|_{L^2(\mathbb{R}^n_+)}\leq C(t+1)^{-\frac{3}{2}}\log_e^2(1+t).
$$
With this new estimate for $\sup\limits_{j\geq0}\|\theta_{j+1}(t)\|_{L^2(\mathbb{R}^n_+)}$, repeating the above process, we infer for any $t>0$ and $j=0, 1, 2,\cdots$
$$
\begin{array}{rcl}
&&\displaystyle \frac{d}{dt}\big((t+1)^k\|\theta_{j+1}(t)\|^2_{L^2(\mathbb{R}^n_+)}\big)\vspace{2mm}\\
&\leq&\displaystyle C(t+1)^{k-1}\big((1+t)^{-\frac{3}{2}}
+(t+1)^{-\frac{3}{2}}\int_0^t(1+s)^{-\frac{3}{2}}\log_e^3(1+s)ds\big)^2\vspace{2mm}\\
&\leq&\displaystyle C(t+1)^{k-4},
\end{array}
$$
from which, we derive for any $t>0$ and $j=0, 1,\cdots$
$$
\sup\limits_{j\geq0}\|\theta_{j+1}(t)\|_{L^2(\mathbb{R}^n_+)}\leq C(1+t)^{-\frac{3}{2}},
$$
which shows (\ref{2.25}) for $n=4$. Further using (\ref{2.14}), we find (\ref{2.24}) is true for $n=4$.\\
\noindent {\bf Case $n\geq5$.} In this case, using (\ref{2.16}), and repeating the proof process for the case $n=4$, we readily find that (\ref{2.24}), (\ref{2.25}) are true for $n\geq5$. From the above arguments, we have proved that (\ref{2.24}), (\ref{2.25}) are valid for $n\geq3$.

\medskip

{\bf $L^2$  Decay for  velocity.}\\
Recall that $\|e^{-tA}a\|_{L^2(\mathbb{R}^n_+)}\longrightarrow0$ as $t\longrightarrow\infty$; and $\sup\limits_{j\geq0}\|\theta_j(t)\|_{L^1(\mathbb{R}^n_+)}\leq C$ by (2.4). Inserting (\ref{2.24}) and (\ref{2.25}) into (\ref{2.11}) with $\frac{n+2}{n}<r<2$, taking $\rho=k(t+1)^{-1}$ with $k>1$, and integrating from $0$ to $t$, then  for $t>0$
\begin{equation} \label{}\notag
\begin{array}{rcl}
&&\displaystyle\frac{d}{dt}((1+t)^k\|u_{j+1}(t)\|^2_{L^2(\mathbb{R}^n_+)})\vspace{2mm}\\
&=&\displaystyle (1+t)^k\frac{d}{dt}\|u_{j+1}(t)\|^2_{L^2(\mathbb{R}^n_+)}+k(1+t)^{k-1}\|u_{j+1}(t)\|^2_{L^2(\mathbb{R}^n_+)}\vspace{2mm}\\
\end{array}
\end{equation}
\begin{equation} \label{}\notag
\begin{array}{rcl}
&\leq&\displaystyle Ck(1+t)^{k-1}\big(\|e^{-tA}a\|_{L^2(\mathbb{R}^n_+)}+(1+t)^{-\frac{n+2}{4}}\int_0^t
\|u_j(s)\|_{L^2(\mathbb{R}^n_+)}\|u_{j+1}(s)\|_{L^2(\mathbb{R}^n_+)}ds\vspace{2mm}\\
&&\displaystyle+\int_0^\frac{t}{2}(t-s)^{-\frac{n}{2}(\frac{1}{r}-\frac{1}{2})}
\|\theta_j(s)\|^{\frac{2}{r}-1}_{L^1(\mathbb{R}^n_+)}\|\theta_j(s)\|^{2(1-\frac{1}{r})}_{L^2(\mathbb{R}^n_+)}ds
+\int^t_\frac{t}{2}\|\theta_j(s)\|_{L^2(\mathbb{R}^n_+)}ds\big)^2\vspace{2mm}\\
&&\displaystyle
+2(1+t)^k\|u_{j+1}(t)\|_{L^2(\mathbb{R}^n_+)} \|\theta_j(t)\|_{L^2(\mathbb{R}^n_+)}
\vspace{2mm}\\
&\leq&\displaystyle C(1+t)^{k-1}\big(\|e^{-tA}a\|_{L^2(\mathbb{R}^n_+)}+(1+t)^{1-\frac{n+2}{4}}\vspace{2mm}\\
&&\displaystyle+Ct^{-\frac{n}{2}(\frac{1}{r}-\frac{1}{2})}\int_0^\frac{t}{2}(1+s)^{-\frac{n+2}{2}(1-\frac{1}{r})}ds
+\int^t_\frac{t}{2}(1+s)^{-\frac{n+2}{4}}ds\big)^2
+C(1+t)^{k-\frac{n+2}{4}}\vspace{2mm}\\
&\leq&\displaystyle C(1+t)^{k-1}\big(\|e^{-tA}a\|^2_{L^2(\mathbb{R}^n_+)}+(1+t)^{-\frac{n-2}{2}}+Ct^{-n(\frac{1}{r}-\frac{1}{2})}\big)+C(1+t)^{k-\frac{n+2}{4}},
\end{array}
\end{equation}
where we use the estimate: $\int_0^\frac{t}{2}(1+s)^{-\frac{n+2}{2}(1-\frac{1}{r})}ds<\infty$ by the choice of $\frac{n+2}{n}<r<2$. From which, if $k$ is taken sufficiently large, then
\begin{equation} \label{2.26}
\begin{array}{rcl}
&&\displaystyle
\sup\limits_{j\geq0}\|u_{j+1}(t)\|^2_{L^2(\mathbb{R}^n_+)}\vspace{2mm}\\
&\leq&\displaystyle C(1+t)^{-k}\|a\|^2_{L^2(\mathbb{R}^n_+)}+C(1+t)^{-k}\int_0^t(1+s)^{k-\frac{n+2}{4}}ds
\vspace{2mm}\\
&&\displaystyle
+C(1+t)^{-k}\int_0^t(1+s)^{k-1}\big(\|e^{-sA}a\|^2_{L^2(\mathbb{R}^n_+)}+(1+s)^{-n(\frac{1}{r}-\frac{1}{2})}+(1+s)^{-\frac{n-2}{2}}\big)ds
\vspace{2mm}\\

&\leq&\displaystyle C(1+t)^{-k}\|a\|^2_{L^2(\mathbb{R}^n_+)}+C((1+t)^{-\frac{n-2}{4}}+(1+t)^{-n(\frac{1}{r}-\frac{1}{2})}+(1+t)^{-\frac{n-2}{2}})
\vspace{2mm}\\
&&\displaystyle
+C(1+t)^{-k}\int_0^t(1+s)^{k-1}\|e^{-sA}a\|^2_{L^2(\mathbb{R}^n_+)}ds
\vspace{2mm}\\
&&\displaystyle \longrightarrow0\quad as \quad t\longrightarrow\infty.
\end{array}
\end{equation}
For the last part of the theorem, assume that $a\in L^\frac{n}{n-1}(\mathbb{R}^n_+)$. Two auxiliary estimates, needed in the sequel, are derived first.  The first one gives decay for $L^1$-norm of $\theta$. The second gives an intermediate decay for the $L^2$-norm of $u$.

\medskip

{\bf $L^1$-decay of $\theta$.}
By (\ref{2.1}), (\ref{2.24}) and (\ref{2.25}), $j=0, 1,\cdots$, $t>0$, it follows that
\begin{equation} \label{2.27}
\begin{array}{rcl}
\displaystyle
\|\theta_{j+1}(t)\|_{L^1(\mathbb{R}^n_+)}
&\leq&\displaystyle \|E(t)b\|_{L^1(\mathbb{R}^n_+)}+C\int_0^t\|\nabla G_{t-s}\|_{L^1(\mathbb{R}^n)}\|u_j(s)\theta_{j+1}(s)\|_{L^1(\mathbb{R}^n_+)}ds\vspace{2mm}\\
&\leq&\displaystyle C t^{-\frac{1}{2}}\|x_nb\|_{L^1(\mathbb{R}^n_+)}+C\int_0^t(t-s)^{-\frac{1}{2}}\|u_j(s)\|_{L^2(\mathbb{R}^n_+)}\|\theta_{j+1}(s)\|_{L^2(\mathbb{R}^n_+)}ds
\vspace{2mm}\\
&\leq&\displaystyle Ct^{-\frac{1}{2}}+C\big(\int_0^\frac{t}{2}+\int_\frac{t}{2}^t\big)
(t-s)^{-\frac{1}{2}}(1+s)^{-\frac{n+2}{4}}ds\vspace{2mm}\\
&\leq&\displaystyle Ct^{-\frac{1}{2}}+C(1+t)^{-\frac{n}{4}}
\vspace{2mm}\\
&\leq&\displaystyle C(1+t)^{-\frac{1}{2}}.
\end{array}
\end{equation}
Let $1<r<\frac{n}{n-1}$, by (\ref{2.25}) and (\ref{2.27}), $j=0, 1,\cdots$, $t>0$, we have
\begin{equation} \label{2.28}
\begin{array}{rcl}
&&\displaystyle
\int_0^\frac{t}{2}(t-s)^{-\frac{n}{2}(\frac{1}{r}-\frac{1}{2})}\|\theta_j(s)\|^{\frac{2}{r}-1}_{L^1(\mathbb{R}^n_+)}
\|\theta_j(s)\|^{2(1-\frac{1}{r})}_{L^2(\mathbb{R}^n_+)}ds
+\int^t_\frac{t}{2}\|\theta_j(s)\|_{L^2(\mathbb{R}^n_+)}ds\vspace{2mm}\\
&\leq& \displaystyle Ct^{-\frac{n}{2}(\frac{1}{r}-\frac{1}{2})}\int_0^\frac{t}{2}(1+s)^{-\frac{1}{2}(\frac{2}{r}-1)-\frac{n+2}{2}(1-\frac{1}{r})}ds
+C\int_\frac{t}{2}^t(1+s)^{-\frac{n+2}{4}}ds\vspace{2mm}\\
&\leq& \displaystyle  C(1+t)^{-\frac{1}{2}(\frac{n}{2}-1)}.
\end{array}
\end{equation}
{\bf Intermediate  $L^2$--decay of $u$.}
Combining Lemma 2.1 for $e^{-tA}$, (\ref{2.11}), (\ref{2.25}) and (\ref{2.28}) yield for $t>0$ and $j=0, 1,\cdots$
\begin{equation} \label{2.29}
\begin{array}{rcl}
&&\displaystyle\frac{d}{dt}\|u_{j+1}(t)\|^2_{L^2(\mathbb{R}^n_+)} +\rho\|u_{j+1}(t)\|^2_{L^2(\mathbb{R}^n_+)}\vspace{2mm}\\
&\leq&\displaystyle C\rho\left((1+t)^{-\frac{1}{2}(\frac{n}{2}-1)}(1+\|a\|_{L^\frac{n}{n-1}(\mathbb{R}^n_+)})+\rho^\frac{n+2}{4}\int_0^t
\|u_j(s)\|_{L^2(\mathbb{R}^n_+)}\|u_{j+1}(s)\|_{L^2(\mathbb{R}^n_+)}ds\right)^2\vspace{2mm}\\
&&\displaystyle
+C\|u_{j+1}(t)\|_{L^2(\mathbb{R}^n_+)}(1+t)^{-\frac{n+2}{4}}.
\end{array}
\end{equation}
Setting as before $\rho=k(t+1)^{-1}$ with large positive integer $k$, multiplying both sides of (\ref{2.29}) by $(t+1)^k$, using the bound (\ref{2.24}), we obtain for $t>0$ and $j=0, 1,\cdots$
\begin{equation} \label{2.30}
\frac{d}{dt}((1+t)^k\|u_{j+1}(t)\|^2_{L^2(\mathbb{R}^n_+)})\leq C(1+t)^{k-1-(\frac{n}{2}-1)}+C(1+t)^{k-\frac{n+2}{4}}.
\end{equation}
Integrating (\ref{2.30}), we have for $t>0$
\begin{equation} \label{2.31}
\sup\limits_{j\geq0}\|u_{j+1}(t)\|_{L^2(\mathbb{R}^n_+)}\leq C((1+t)^{-\frac{n}{4}+\frac{1}{2}}+(1+t)^{-\frac{n-2}{8}})\leq C(1+t)^{-\frac{n-2}{8}}.
\end{equation}
We now use the last estimate to obtain the conclusion of the theorem.
Combine (\ref{2.31}) and (\ref{2.29}), set as specified  $\rho=k(t+1)^{-1}$  with large positive integer $k$, and
multiply (\ref{2.29}) by $(t+1)^k$. Then for $t>0$ and $j=0, 1,\cdots$
\begin{equation} \label{2.32}
\begin{array}{rcl}
&&\displaystyle \frac{d}{dt}((1+t)^k\|u_{j+1}(t)\|^2_{L^2(\mathbb{R}^n_+)})\vspace{2mm}\\
&\leq&\displaystyle
C(1+t)^{k-1}\left((1+t)^{-\frac{1}{2}(\frac{n}{2}-1)}
+(1+t)^{-\frac{n+2}{4}}\int_0^t(1+s)^{-\frac{n-2}{4}}ds\right)^2\vspace{2mm}\\
&&\displaystyle+C(1+t)^{k-\frac{n-2}{8}-\frac{n+2}{4}}.
\end{array}
\end{equation}
Integrating (\ref{2.32}) yields for $k$ sufficiently large
\begin{equation} \label{2.33}
\begin{array}{rcl}
&&\displaystyle\sup\limits_{j\geq0}\|u_{j+1}(t)\|_{L^2(\mathbb{R}^n_+)}\vspace{2mm}\\
&\leq&\displaystyle C((1+t)^{-\frac{1}{2}(\frac{n}{2}-1)}+(1+t)^{-\frac{3}{8}(\frac{n}{2}-1)}+(1+t)^{-\frac{n+2}{4}}\int_0^t(1+s)^{-\frac{n-2}{4}}ds)
\vspace{2mm}\\
&\leq&\displaystyle C((1+t)^{-\frac{1}{2}(\frac{n}{2}-1)}+(1+t)^{-\frac{3}{8}(\frac{n}{2}-1)})\vspace{2mm}\\
&&\displaystyle +C(1+t)^{-\frac{n+2}{4}}\left\{
\begin{array}{lll}
1&\mbox{if}& n\geq 7,\vspace{2mm}\\
\log_e(1+t) &\mbox{if}& n=6,\vspace{2mm}\\
(1+t)^{\frac{3}{2}-\frac{n}{4}}&\mbox{if}& 3\leq n\leq 5,
\end{array}
\right.\vspace{2mm}\\
&\leq&\displaystyle C(1+t)^{-\frac{3}{8}(\frac{n}{2}-1)}.
\end{array}
\end{equation}
An inductive bootstrap  process of using the steps yielding (\ref{2.33}) yields
\begin{equation}\label{2.34}
\sup\limits_{j\geq0}\|u_{j+1}(t)\|_{L^2(\mathbb{R}^n_+)}
\leq \displaystyle C_m((1+t)^{-\frac{1}{2}(\frac{n}{2}-1)}+(1+t)^{-\frac{m-1}{m}(\frac{1}{2}(\frac{n}{2}-1)})
\end{equation}
Hence it follows that for any small $\epsilon>0$, there exists a large number $m_0=m_0(\epsilon)>0$ such that
\begin{equation} \label{2.35}
\sup\limits_{j\geq0}\|u_{j+1}(t)\|_{L^2(\mathbb{R}^n_+)}\leq C_{m_0}(1+t)^{-\frac{1}{2}(\frac{n}{2}-1)+\epsilon}.
\end{equation}
From (\ref{2.5}), (\ref{2.6}), (\ref{2.13}), (\ref{2.16}), (\ref{2.25}), (\ref{2.26}) and (\ref{2.35}), we find there exists $C>0$ independent of $j$ such that for $j=0, 1, \cdots$
$$
\|u_j\|_{L^\infty(0, \infty; L^2(\mathbb{R}^n_+))}+\|\theta_j\|_{L^\infty(0, \infty; L^2(\mathbb{R}^n_+))}+\|\nabla u_j\|_{L^2(0, \infty; L^2(\mathbb{R}^n_+))}+\|\nabla \theta_j\|_{L^2(0, \infty; L^2(\mathbb{R}^n_+))}\leq C.
$$
Using weak convergence properties, we conclude there exists a pair of function $(u, \theta)$ and select a subsequence of $(u_j, \theta_j)$ if necessary such that as $j\longrightarrow\infty$
$$
u_j\rightharpoonup u \;\;\;\mbox{weakly in }\;\; L^\infty(0, \infty; L^2(\mathbb{R}^n_+)),\;\;\;\; \nabla u_j\rightharpoonup \nabla u \;\;\;\mbox{weakly in }\;\; L^2(0, \infty; L^2(\mathbb{R}^n_+));
$$
$$
\theta_j\rightharpoonup \theta \;\;\;\mbox{weakly in }\;\; L^\infty(0, \infty; L^2(\mathbb{R}^n_+)),\;\;\;\; \nabla \theta_j\rightharpoonup \nabla\theta \;\;\;\mbox{weakly in } \;\;L^2(0, \infty; L^2(\mathbb{R}^n_+));
$$
which yields for $\;T>0$, $u\in L^\infty(0, T; L^2_\sigma(\mathbb{R}^n_+))\bigcap L^2(0, T; H^1_0(\mathbb{R}^n_+))$, $\theta\in L^\infty(0, T;
L^2(\mathbb{R}^n_+)) \bigcap L^2(0, T; H^1_0(\mathbb{R}^n_+))$, such that $(u, \theta)$ is a weak solution of (\ref{1.1}) satisfying the estimates (\ref{1.2}), (\ref{1.3}), (\ref{1.5}), (\ref{1.6}). This completes the proof of Theorem \ref{th:1}.$\,\,\,\,\,\Box$

{\bf Proof of Theorem \ref{th:4}}\quad Let  $j=0,1,2\cdots$. Consider the successive approximation of problem (1.1)
\begin{equation} \label{2.36}
\left\{
\begin{array}{ll}
\displaystyle \theta_0(t)=E(t)b,
\qquad u_0(t)=e^{-tA}a,\vspace{2mm}\\
\displaystyle
\theta_{j+1}(t)=\theta_0(t)-\int_0^tE(t-s)u_j(s)\cdot\nabla\theta_j(s)ds,\vspace{2mm}\\
\displaystyle
u_{j+1}(t)=u_0(t)-\int_0^te^{-(t-s)A}P\big(u_j(s)\cdot\nabla
u_j(s)-\theta_j(s)e_n\big)ds
\end{array}
\right.
\end{equation}
The method of  proof used is based on Kato's ideas in \cite{Ka}, for  solutions to the Navier-Stokes equations. We introduce appropriate modifications to handle the extra temperature terms (that is, the terms $L_j$ below). Let $\alpha=1+2\epsilon_1$, $\delta=1-\epsilon_1$, where the number $\epsilon_1\in(0, \frac{1}{3})$. Then $1<\alpha<2$ and $0<\delta<1$. Define for $j=0, 1,2\cdots$
\begin{equation} \label{}\notag
K_j=\sup\limits_{0<t<\infty}(t^\frac{1-\delta}{2}\|u_j(t)\|_{L^\frac{n}{\delta}(\mathbb{R}^n_+)});
\;\;
K^\prime_j=\sup\limits_{0<t<\infty}(t^\frac{1}{2}\|\nabla u_j(t)\|_{L^n(\mathbb{R}^n_+)});
\;\;
L_j=\sup\limits_{0<t<\infty}(t^\frac{3-\alpha}{2}\|\theta_j(t)\|_{L^\frac{n}{\alpha}(\mathbb{R}^n_+)}).
\end{equation}
By Lemma \ref{l:2.1}, for $t>0$
\begin{equation} \label{}\notag
\|u_0(t)\|_{L^\frac{n}{\delta}(\mathbb{R}^n_+)}=\|e^{-tA}a\|_{L^\frac{n}{\delta}(\mathbb{R}^n_+)}\leq Ct^{-\frac{1-\delta}{2}}\|a\|_{L^n(\mathbb{R}^n_+)};
\end{equation}
\begin{equation} \label{}\notag
\|\nabla u_0(t)\|_{L^n(\mathbb{R}^n_+)}=\|\nabla e^{-tA}a\|_{L^n(\mathbb{R}^n_+)}\leq Ct^{-\frac{1}{2}}\|a\|_{L^n(\mathbb{R}^n_+)};
\end{equation}
and
\begin{equation} \label{}\notag
\|\theta_0(t)\|_{L^\frac{n}{\alpha}(\mathbb{R}^n_+)}=\|E(t)b\|_{L^\frac{n}{\alpha}(\mathbb{R}^n_+)}\leq Ct^{-\frac{3-\alpha}{2}}\|b\|_{L^\frac{n}{3}(\mathbb{R}^n_+)}.
\end{equation}
Thus
\begin{equation} \label{2.37}
K_0+K_0^\prime\leq C\|a\|_{L^n(\mathbb{R}^n_+)}\quad and \quad L_0\leq C\|b\|_{L^\frac{n}{3}(\mathbb{R}^n_+)}.
\end{equation}
From (\ref{2.37}), we have $K_0+K_0^\prime+L_0\leq C\|a\|_{L^n(\mathbb{R}^n_+)}+C\|b\|_{L^\frac{n}{3}(\mathbb{R}^n_+)}<\infty$. Assume $K_j+K_j^\prime+L_j<\infty$ for some fixed $j\geq0$. By Lemma \ref{l:2.1} and (\ref{2.36}) it follows that for  $t>0$
\begin{equation} \label{}\notag
\begin{array}{rcl}
\displaystyle\|u_{j+1}(t)\|_{L^\frac{n}{\delta}(\mathbb{R}^n_+)}
&\leq&\displaystyle\|e^{-tA}a\|_{L^\frac{n}{\delta}(\mathbb{R}^n_+)}+\int_0^t\big\|e^{-(t-s)A}P(u_j(s)\cdot\nabla
u_j(s)-\theta_j(s)e_n)\big\|_{L^\frac{n}{\delta}(\mathbb{R}^n_+)}ds\vspace{2mm}\\
&\leq&\displaystyle K_0t^{-\frac{1-\delta}{2}}+C\int_0^t(t-s)^{-\frac{1}{2}}\|u_j(s)\|_{L^\frac{n}{\delta}(\mathbb{R}^n_+)}
\|\nabla u_j(s)\|_{L^n(\mathbb{R}^n_+)}ds\vspace{2mm}\\
&&\displaystyle
+C\int_0^t(t-s)^{-\frac{n}{2}(\frac{\alpha}{n}-\frac{\delta}{n})}
\|\theta_j(s)\|_{L^\frac{n}{\alpha}(\mathbb{R}^n_+)}ds
\vspace{2mm}\\
&\leq&\displaystyle K_0t^{-\frac{1-\delta}{2}}+CK_jK^\prime_j\int_0^t(t-s)^{-\frac{1}{2}}s^{-1+\frac{\delta}{2}}ds
+CL_j\int_0^t(t-s)^{-\frac{\alpha-\delta}{2}}s^{-\frac{3-\alpha}{2}}ds
\vspace{2mm}\\
&\leq&\displaystyle K_0t^{-\frac{1-\delta}{2}}+CK_jK^\prime_jt^{-\frac{1-\delta}{2}}+CL_jt^{-\frac{1-\delta}{2}}.
\end{array}
\end{equation}
Hence
\begin{equation} \label{2.38}
K_{j+1}\leq K_0+CK_jK^\prime_j+CL_j<\infty.
\end{equation}
Similarly for  $t>0$
\begin{equation} \label{}\notag
\begin{array}{rcl}
\displaystyle\|\nabla u_{j+1}(t)\|_{L^n(\mathbb{R}^n_+)}
&\leq&\displaystyle\|\nabla e^{-tA}a\|_{L^n(\mathbb{R}^n_+)}+C\int_0^t(t-s)^{-\frac{1}{2}-\frac{n}{2}(\frac{\alpha}{n}-\frac{1}{n})}
\|\theta_j(s)\|_{L^\frac{n}{\alpha}(\mathbb{R}^n_+)}ds\vspace{2mm}\\
&&\displaystyle+C\int_0^t(t-s)^{-\frac{1+\delta}{2}}\|u_j(s)\|_{L^\frac{n}{\delta}(\mathbb{R}^n_+)}
\|\nabla u_j(s)\|_{L^n(\mathbb{R}^n_+)}ds
\vspace{2mm}\\
&\leq&\displaystyle K^\prime_0t^{-\frac{1}{2}}+CL_j\int_0^t(t-s)^{-\frac{\alpha}{2}}s^{-\frac{3-\alpha}{2}}ds
\vspace{2mm}\\
&&\displaystyle +CK_jK^\prime_j\int_0^t(t-s)^{-\frac{1+\delta}{2}}s^{-1+\frac{\delta}{2}}ds
\vspace{2mm}\\
&\leq&\displaystyle K^\prime_0t^{-\frac{1}{2}}+CK_jK^\prime_jt^{-\frac{1}{2}}+CL_jt^{-\frac{1}{2}},
\end{array}
\end{equation}
which implies
\begin{equation} \label{2.39}
K^\prime_{j+1}\leq K^\prime_0+CK_jK^\prime_j+CL_j<\infty.
\end{equation}
Similar estimate as for (\ref{2.38}), yields for $t>0$
\begin{equation} \label{}\notag
\begin{array}{rcl}
\displaystyle\|\theta_{j+1}(t)\|_{L^\frac{n}{\alpha}(\mathbb{R}^n_+)}
&\leq&\displaystyle\|E(t)b\|_{L^\frac{n}{\alpha}(\mathbb{R}^n_+)}+2\int_0^t\|\nabla G_{t-s}\|_{L^{(1-\frac{\delta}{n})^{-1}}(\mathbb{R}^n)}\|(u_j(s)
\theta_j(s))\|_{L^\frac{n}{\alpha+\delta}(\mathbb{R}^n_+)}ds\vspace{2mm}\\
&\leq&\displaystyle L_0t^{-\frac{3-\alpha}{2}}+C\int_0^t(t-s)^{-\frac{1}{2}-\frac{\delta}{2}}
\|u_j(s)\|_{L^\frac{n}{\delta}(\mathbb{R}^n_+)}\|\theta_j(s)\|_{L^\frac{n}{\alpha}(\mathbb{R}^n_+)}ds\vspace{2mm}\\
&\leq&\displaystyle L_0t^{-\frac{3-\alpha}{2}}+CK_jL_j\int_0^t(t-s)^{-\frac{1+\delta}{2}}s^{-2+\frac{\alpha+\delta}{2}}ds
\vspace{2mm}\\
&\leq&\displaystyle L_0t^{-\frac{3-\alpha}{2}}+CK_jL_jt^{-\frac{3-\alpha}{2}},
\end{array}
\end{equation}
from which
\begin{equation} \label{2.40}
L_{j+1}\leq L_0+CK_jL_j<\infty.
\end{equation}
Set $M_j=K_j+K^\prime_j+L_j$ for $j=0, 1,2,\cdots$. For the above fixed $j\geq0$,  inequalities  (\ref{2.38})--(\ref{2.40}) yield
\begin{equation} \label{2.41}
M_{j+1}\leq K_0+K_0^\prime+L_0+CM_j^2+CL_j\leq
\left\{\begin{array}{lll}
\displaystyle K_0+K_0^\prime+(1+C)L_0+CM_0^2&\mbox{if}&j=0,\vspace{2mm}\\
\displaystyle K_0+K_0^\prime+C_1L_0+C_2\big(\frac{M_j+M_{j-1}}{2}\big)^2&\mbox{if}&j\geq1,
\end{array}
\right.
\end{equation}
where $C_1>1+2C$, $C_2>4C$.\\
Set $B_0=K_0+K_0^\prime+C_1L_0$. Note that   $B_0\leq C(\|a\|_{L^n(\mathbb{R}^n_+)}+\|b\|_{L^\frac{n}{3}(\mathbb{R}^n_+)})\leq C\epsilon_0$  by (\ref{2.41}) and our hypothesis.  Choose $\epsilon_0$ so small that  $B_0<\frac{1}{4C_2}$. Set $\chi_0=(2C_2)^{-1}\big(1-\sqrt{1-4B_0C_2}\big)$. Then $\chi_0=B_0+C_2\chi_0^2$ and $\chi_0\leq2B_0$. We shall show that
\begin{equation} \label{2.42}
M_j\leq\chi_0,\;\;\;\quad j=0,1,2\cdots.
\end{equation}
The proof is  by induction on $j$. For $j=0$, it easily follows that $M_0<B_0<\chi_0$.
Suppose that $M_i\leq\chi_0$ for any $0\leq i\leq j$ with the above fixed $j\geq0$.
Then from (\ref{2.41}), it follows that
\[M_1\leq  B_0+CM_0^2\leq B_0+C_2\chi_0^2\leq \chi_0,\;\;\mbox{and}\,\, M_{j+1}\leq B_0+C_2(\frac{M_j+M_{j-1}}{2})^2\leq \chi_0\;\; for\;\; j\geq1,\]
which establishes (\ref{2.42}). Now combining Lemma \ref{l:2.1}, (\ref{2.36}) and (\ref{2.42}),
then for  $j=0, 1,2,\cdots$, $t>0$ and $ n\leq q\leq\infty$
\begin{equation} \label{2.43}
\begin{array}{rcl}
\displaystyle\|u_{j+1}(t)\|_{L^q(\mathbb{R}^n_+)}
&\leq&\displaystyle\|e^{-tA}a\|_{L^q(\mathbb{R}^n_+)}+\int_0^t\|e^{-(t-s)A}P(u_j(s)\cdot\nabla
u_j(s)-\theta_j(s)e_n)\|_{L^q(\mathbb{R}^n_+)}ds\vspace{2mm}\\
&\leq&\displaystyle Ct^{-\frac{n}{2}(\frac{1}{n}-\frac{1}{q})}\|a\|_{L^n(\mathbb{R}^n_+)}+C\int_0^t(t-s)^{-\frac{n}{2}(\frac{\alpha}{n}-\frac{1}{q})}
\|\theta_j(s)\|_{L^\frac{n}{\alpha}(\mathbb{R}^n_+)}ds\vspace{2mm}\\
&&\displaystyle
+C\int_0^t(t-s)^{-\frac{n}{2}(\frac{1+\delta}{n}-\frac{1}{q})}\|u_j(s)\|_{L^\frac{n}{\delta}(\mathbb{R}^n_+)}
\|\nabla u_j(s)\|_{L^n(\mathbb{R}^n_+)}ds\vspace{2mm}\\
&\leq&\displaystyle Ct^{-\frac{n}{2}(\frac{1}{n}-\frac{1}{q})}\|a\|_{L^n(\mathbb{R}^n_+)}+CK_jK^\prime_j\int_0^t(t-s)^{-\frac{1+\delta}{2}+\frac{n}{2q}}
s^{-1+\frac{\delta}{2}}ds
\vspace{2mm}\\
&&\displaystyle +CL_j\int_0^t(t-s)^{-\frac{\alpha}{2}+\frac{n}{2q}}s^{-\frac{3-\alpha}{2}}ds
\vspace{2mm}\\
&\leq&\displaystyle C(\|a\|_{L^n(\mathbb{R}^n_+)}+M_j+M_j^2)t^{-\frac{1}{2}+\frac{n}{2q}}\vspace{2mm}\\
&\leq&\displaystyle C\epsilon_0t^{-\frac{n}{2}(\frac{1}{n}-\frac{1}{q})}.
\end{array}
\end{equation}
This establishes the $L^q, n\leq q\leq\infty$ decay for the velocity part of the $j$-solution approximation.
We now will handle the $L^q$ decay of the temperatures.
By (\ref{2.36}) and the definition of the operator $E(t)$ in (\ref{2.2}), we have for $j=0, 1,2,\cdots$ and $t>0$
\begin{equation} \label{2.44}
\theta_{j+1}(t)=E(\frac{t}{2})\theta_{j+1}(\frac{t}{2})-\int^t_\frac{t}{2}E(t-s)Pu_j(s)\cdot\nabla\theta_j(s)ds.
\end{equation}
Let $n\leq q<\infty$ and assume $n\leq q<\frac{n}{\alpha-1}$, which is possible by the choice of $\alpha=1+2\epsilon_1$, and if $\epsilon_1>0$ is sufficiently small. From Lemma \ref{l:2.1}, (\ref{2.42})--(\ref{2.44}), we get for  $j=0, 1,2,\cdots$ and $t>0$
\begin{equation} \label{2.45}
\begin{array}{rcl}
\displaystyle\|\theta_{j+1}(t)\|_{L^q(\mathbb{R}^n_+)}
&\leq&\displaystyle\|E(\frac{t}{2})\theta_{j+1}(\frac{t}{2})\|_{L^q(\mathbb{R}^n_+)}\vspace{2mm}\\
&&\displaystyle +2\int^t_\frac{t}{2}\|\nabla G_{t-s}\|_{L^{(1+\frac{1}{q}-\frac{\alpha}{n})^{-1}}(\mathbb{R}^n)}\|(u_j(s)
\theta_j(s))\|_{L^\frac{n}{\alpha}(\mathbb{R}^n_+)}ds\vspace{2mm}\\
&\leq&\displaystyle Ct^{-\frac{n}{2}(\frac{\alpha}{n}-\frac{1}{q})}\|\theta_{j+1}(t)\|_{L^\frac{n}{\alpha}(\mathbb{R}^n_+)}
\vspace{2mm}\\
&&\displaystyle
+C\int^t_\frac{t}{2}(t-s)^{-\frac{1}{2}-\frac{n}{2}(\frac{\alpha}{n}-\frac{1}{q})}
\|u_j(s)\|_{L^\infty(\mathbb{R}^n_+)}\|\theta_j(s)\|_{L^\frac{n}{\alpha}(\mathbb{R}^n_+)}ds
\vspace{2mm}\\
&\leq&\displaystyle C\epsilon_0t^{-\frac{n}{2}(\frac{\alpha}{n}-\frac{1}{q})-\frac{3-\alpha}{2}}
+C\epsilon_0^2\int^t_\frac{t}{2}(t-s)^{-\frac{1}{2}-\frac{n}{2}(\frac{\alpha}{n}-\frac{1}{q})}s^{-\frac{1}{2}-\frac{3-\alpha}{2}}ds
\vspace{2mm}\\
&\leq&\displaystyle C\epsilon_0t^{-\frac{3}{2}+\frac{n}{2q}}
+C\epsilon_0^2t^{-\frac{1}{2}-\frac{3-\alpha}{2}+1-\frac{1}{2}-\frac{n}{2}(\frac{\alpha}{n}-\frac{1}{q})}\vspace{2mm}\\
&\leq&\displaystyle C\epsilon_0t^{-\frac{3}{2}+\frac{n}{2q}}.
\end{array}
\end{equation}
Let $n\leq q\leq\infty$. From Lemma \ref{l:2.1}, (\ref{2.43}) and (\ref{2.44}), we conclude that for  $j=1,2,\cdots$ and $t>0$
\begin{equation} \label{}\notag
\begin{array}{rcl}
\displaystyle\|\nabla\theta_{j+1}(t)\|_{L^q(\mathbb{R}^n_+)}
&\leq&\displaystyle\|\nabla E(\frac{t}{2})\theta_{j+1}(\frac{t}{2})\|_{L^q(\mathbb{R}^n_+)} +\int^t_\frac{t}{2}\big\|\nabla\, E(t-s)(u_j(s)\cdot
\nabla\theta_j(s))\big\|_{L^q(\mathbb{R}^n_+)}ds\vspace{2mm}\\
&\leq&\displaystyle Ct^{-\frac{1}{2}}\|\theta_{j+1}(t)\|_{L^q(\mathbb{R}^n_+)}
+C\int^t_\frac{t}{2}(t-s)^{-\frac{1}{2}}
\|u_j(s)\|_{L^\infty(\mathbb{R}^n_+)}\|\nabla\theta_j(s)\|_{L^q(\mathbb{R}^n_+)}ds
\vspace{2mm}\\
&\leq&\displaystyle C\epsilon_0t^{-2+\frac{n}{2q}}
+C\epsilon_0N_j(t)\int^t_\frac{t}{2}(t-s)^{-\frac{1}{2}}s^{-\frac{1}{2}-2+\frac{n}{2q}}ds
\vspace{2mm}\\
&\leq&\displaystyle C\epsilon_0t^{-2+\frac{n}{2q}}+C\epsilon_0N_j(t)t^{-2+\frac{n}{2q}}.
\end{array}
\end{equation}
Here $N_j(t)=\sup\limits_{0<s\leq t}\{s^{2-\frac{n}{2q}}\|\nabla\theta_j(s)\|_{L^q(\mathbb{R}^n_+)}\}$ with $j=1,2,\cdots$.
Then  for  $j=1,2,\cdots$ and $t>0$
\begin{equation} \label{2.46}
N_{j+1}(t)\leq C\epsilon_0+C\epsilon_0N_j(t).
\end{equation}
Note that by Lemma \ref{l:2.1}, for $n\leq q\leq\infty$ and $t>0$
\begin{equation} \label{}\notag
\|\nabla\theta_0(s)\|_{L^q(\mathbb{R}^n_+)}=\|\nabla e^{-tA}b\|_{L^q(\mathbb{R}^n_+)}\leq Ct^{-2+\frac{n}{2q}}\|b\|_{L^\frac{n}{3}(\mathbb{R}^n_+)}.
\end{equation}
Let  $C\epsilon_0\leq\frac{1}{2}$ in (\ref{2.46}), then  for  $j=0, 1,2,\cdots$ and $t>0$
\begin{equation} \label{2.47}
\|\nabla\theta_{j+1}(t)\|_{L^q(\mathbb{R}^n_+)}\leq C\epsilon_0t^{-2+\frac{n}{2q}},\quad\,n\leq q\leq\infty.
\end{equation}
In addition, from (\ref{2.45}), (\ref{2.47}), and the Gagliardo-Nirenberg inequality on the half space (see (4.1) in \cite{BM} for example), one has for  $j=1,2,\cdots$, $t>0$
\begin{equation} \label{2.48}
\|\theta_{j+1}(t)\|_{L^\infty(\mathbb{R}^n_+)}\leq C\|\theta_{j+1}(t)\|^\frac{1}{2}_{L^{2n}(\mathbb{R}^n_+)}\|\nabla \theta_{j+1}(t)\|^\frac{1}{2}_{L^{2n}(\mathbb{R}^n_+)}\leq C\epsilon_0t^{-\frac{3}{2}}.
\end{equation}
Note that from (\ref{2.36}),  we have  $t>0$
\begin{equation} \label{2.49}
u_{j+1}(t)=e^{-\frac{t}{2}A}u_{j+1}(\frac{t}{2})-\int_\frac{t}{2}^te^{-(t-s)A}P(u_j(s)\cdot\nabla
u_j(s)-\theta_j(s)e_n)ds.
\end{equation}
Let $n\leq q<\infty$. Using (\ref{2.43}), (\ref{2.45}), (\ref{2.49}) and Lemma \ref{l:2.1}, we have for $j\geq1$ and $t>0$
\begin{equation} \label{2.50}
\begin{array}{rcl}
&&\displaystyle\|\nabla u_{j+1}(t)\|_{L^q(\mathbb{R}^n_+)}\vspace{2mm}\\
&\leq&\displaystyle\|\nabla e^{-\frac{t}{2}A}u_{j+1}(\frac{t}{2})\|_{L^q(\mathbb{R}^n_+)}+\int_\frac{t}{2}^t\|\nabla e^{-(t-s)A}P(u_j(s)\cdot\nabla
u_j(s)-\theta_j(s)e_n)\|_{L^q(\mathbb{R}^n_+)}ds\vspace{2mm}\\
&\leq&\displaystyle Ct^{-\frac{1}{2}}\|u_{j+1}(\frac{t}{2})\|_{L^q(\mathbb{R}^n_+)}
+C\int_\frac{t}{2}^t(t-s)^{-\frac{1}{2}}\|u_j(s)\|_{L^\infty(\mathbb{R}^n_+)}
\|\nabla u_j(s)\|_{L^q(\mathbb{R}^n_+)}ds\vspace{2mm}\\
&&\displaystyle+C\int_\frac{t}{2}^t(t-s)^{-\frac{1}{2}}\|\theta_j(s)\|_{L^q(\mathbb{R}^n_+)}ds
\vspace{2mm}\\
&\leq&\displaystyle C\epsilon_0t^{-1+\frac{n}{2q}}+C\epsilon_0\Upsilon_j(t)\int_\frac{t}{2}^t(t-s)^{-\frac{1}{2}}
s^{-\frac{3}{2}+\frac{n}{2q}}ds+C\epsilon_0\int_\frac{t}{2}^t(t-s)^{-\frac{1}{2}}s^{-\frac{3}{2}+\frac{n}{2q}}ds
\vspace{2mm}\\
&\leq&\displaystyle C\epsilon_0t^{-1+\frac{n}{2q}}+C\epsilon_0\Upsilon_j(t)t^{-1+\frac{n}{2q}},
\end{array}
\end{equation}
where $\Upsilon_j(t)=\sup\limits_{0<s\leq t}\{s^{1-\frac{n}{2q}}\|\nabla u_j(s)\|_{L^q(\mathbb{R}^n_+)}\}$ with $j=1,2,\cdots$.
Hence for $t>0$
\begin{equation} \label{2.51}
\Upsilon_{j+1}(t)\leq C\epsilon_0+C\epsilon_0\Upsilon_j(t).
\end{equation}
Note that by Lemma \ref{l:2.1}, for $n\leq q\leq\infty,\,t>0$
\begin{equation} \label{}\notag
\|\nabla u_0(s)\|_{L^q(\mathbb{R}^n_+)}=\|\nabla e^{-tA}a\|_{L^q(\mathbb{R}^n_+)}\leq Ct^{-1+\frac{n}{2q}}\|a\|_{L^n(\mathbb{R}^n_+)}.
\end{equation}
Let  $C\epsilon_0\leq\frac{1}{2}$ in (\ref{2.51}), then  for  $j=0, 1,2,\cdots$ and $t>0$
\begin{equation} \label{2.52}
\|\nabla u_{j+1}(t)\|_{L^q(\mathbb{R}^n_+)}\leq C\epsilon_0t^{-1+\frac{n}{2q}},\quad \forall\,q\in[n, \infty).
\end{equation}
In addition, by (\ref{2.43}), (\ref{2.45}), (\ref{2.49}) and Lemma \ref{l:2.1}, for any $j=1,2,\cdots$ and $t>0$
\begin{equation} \label{2.53}
\begin{array}{rcl}
&&\displaystyle\|\nabla u_{j+1}(t)\|_{L^\infty(\mathbb{R}^n_+)}\vspace{2mm}\\
&\leq&\displaystyle\|\nabla e^{-\frac{t}{2}A}u_{j+1}(\frac{t}{2})\|_{L^\infty(\mathbb{R}^n_+)}+\int_\frac{t}{2}^t\|\nabla e^{-(t-s)A}P(u_j(s)\cdot\nabla
u_j(s)-\theta_j(s)e_n)\|_{L^\infty(\mathbb{R}^n_+)}ds\vspace{2mm}\\
&\leq&\displaystyle Ct^{-\frac{1}{2}}\|u_{j+1}(\frac{t}{2})\|_{L^\infty(\mathbb{R}^n_+)}
+C\int_\frac{t}{2}^t(t-s)^{-\frac{1}{2}-\frac{1}{4}}\|u_j(s)\|_{L^{2n}(\mathbb{R}^n_+)}
\|\nabla u_j(s)\|_{L^\infty(\mathbb{R}^n_+)}ds\vspace{2mm}\\
&&\displaystyle+C\int_\frac{t}{2}^t(t-s)^{-\frac{1}{2}-\frac{1}{4}}\|\theta_j(s)\|_{L^{2n}(\mathbb{R}^n_+)}ds
\vspace{2mm}\\
&\leq&\displaystyle C\epsilon_0t^{-1}+C\epsilon_0\Pi_j(t)\int_\frac{t}{2}^t(t-s)^{-\frac{3}{4}}
s^{-\frac{5}{4}}ds+C\epsilon_0\int_\frac{t}{2}^t(t-s)^{-\frac{3}{4}}s^{-\frac{5}{4}}ds
\vspace{2mm}\\
&\leq&\displaystyle C\epsilon_0t^{-1}+C\epsilon_0\Pi_j(t)t^{-1},
\end{array}
\end{equation}
where $\Pi_j(t)=\sup\limits_{0<s\leq t}\{s\|\nabla u_j(s)\|_{L^\infty(\mathbb{R}^n_+)}\}$ with $j=0, 1,2,\cdots$.
Thus $\Pi_{j+1}(t)\leq C\epsilon_0+C\epsilon_0\Pi_j(t).$
Let again  $C\epsilon_0\leq\frac{1}{2}$. Hence for $j=0, 1,2,\cdots$ and $t>0$
\begin{equation} \label{2.54}
\|\nabla u_{j+1}(t)\|_{L^\infty(\mathbb{R}^n_+)}\leq C\epsilon_0t^{-1}.
\end{equation}
From (\ref{2.43}), (\ref{2.45}), (\ref{2.47}), (\ref{2.48}), (\ref{2.52}) and (\ref{2.54}),
using weak convergence properties, we conclude there exists a pair of function $(u, \theta)$ and select a subsequence of $(u_j, \theta_j)$ if necessary such that for each $0<T<\infty$, as $j\longrightarrow\infty$
$$
u_j\rightharpoonup u \;\;\;\mbox{weakly in }\;\; L^\infty(0, \infty; W^{1, q}(\mathbb{R}^n_+)),\;\;\;\; \theta_j\rightharpoonup \theta \;\;\;\mbox{weakly in }\;\; L^\infty(0, \infty; W^{1, q}(\mathbb{R}^n_+)).
$$
Moreover for $\,n\leq q\leq\infty$ and $t>0$
$$
\|u(t)\|_{L^q(\mathbb{R}^n_+)}\leq\liminf\limits_{j\longrightarrow\infty}\|u_{j+1}(t)\|_{L^q(\mathbb{R}^n_+)}\leq C\epsilon_0t^{-\frac{n}{2}(\frac{1}{n}-\frac{1}{q})};
$$
$$
\|\nabla u(t)\|_{L^q(\mathbb{R}^n_+)}\leq\liminf\limits_{j\longrightarrow\infty}\|\nabla u_{j+1}(t)\|_{L^q(\mathbb{R}^n_+)}\leq C\epsilon_0t^{-1+\frac{n}{2q}};
$$
$$
\|\theta(t)\|_{L^q(\mathbb{R}^n_+)}\leq\liminf\limits_{j\longrightarrow\infty}\|\theta_{j+1}(t)\|_{L^q(\mathbb{R}^n_+)}
\leq C\epsilon_0t^{-\frac{3}{2}+\frac{n}{2q}};
$$
$$
\|\nabla\theta(t)\|_{L^q(\mathbb{R}^n_+)}\leq\liminf\limits_{j\longrightarrow\infty}\|\nabla\theta_{j+1}(t)\|_{L^q(\mathbb{R}^n_+)}\leq C\epsilon_0t^{-2+\frac{n}{2q}}.
$$
Using the above estimates, parabolic regularity theory and Serrin criteria for Navier-Stokes equations, we readily find that $(u, \theta)$ is a  strong solution of (\ref{1.1}) which satisfies the estimates in Theorem \ref{th:4}. $\,\,\,\,\,\,\Box$

Before establishing Theorem \ref{th:5}, we prove an auxiliary Proposition, where part of the decay estimates for the temperature $\theta$ will be done.
\begin{proposition} \label{1:2.3} Let $a\in L^n_\sigma(\mathbb{R}^n_+)\bigcap L^\frac{n}{n-1}(\mathbb{R}^n_+)$, and $b\in L^1(\mathbb{R}^n_+)\bigcap L^2(\mathbb{R}^n_+)\bigcap L^\frac{n}{3}(\mathbb{R}^n_+)$, $n\geq3$. Assume  (\ref{1.4}) holds.
Let $(u, \theta)$ be the strong solution of (\ref{1.1}) obtained in Theorem \ref{1.4}. Then for any $\frac{n}{n-1}\leq r<2$, $r\leq q\leq \infty$ and $t>0$, $\|\nabla^k u(t)\|_{L^q(\mathbb{R}^n_+)}\leq Ct^{-\frac{k}{2}-\frac{n}{2}(\frac{1}{r}-\frac{1}{q})}\,\quad for\,\,\,\,\, k=0, 1.$
\end{proposition}
{\bf Proof}\quad For the proof we start with  estimates for $\theta \in L^q$ which  give the decay required in the case $k=0$ in Theorem \ref{th:5}, then show that the $L^r$ norms of the velocity are bounded and finally  prove the estimates of the proposition. Recall first  that $b\in L^1(\mathbb{R}^n_+)$ and by (\ref{2.27}),
\begin{equation} \label{2.55}
\|\theta(t)\|_{L^1(\mathbb{R}^n_+)}\leq\ C(1+t)^{-\frac{1}{2}},\;\;\;\;\forall\;t>0.
\end{equation}
{\bf Auxiliary estimates: }
 Let $\frac{n}{n-1}\leq r<2$, $r\leq q\leq\infty$, and  $(u, \theta)$ be the strong solution of (\ref{1.1}).
Set
\begin{equation} \label{2.56}
L(q, r,t)=\sup\limits_{0<s\leq t}(s^{\frac{n}{2}(\frac{1}{r}-\frac{1}{q})}\|u(s)\|_{L^q(\mathbb{R}^n_+)}).
\end{equation}
Let $1< r_1\leq r_2<\infty$. Then for any matrix function $F$ and $t>0$
\begin{equation} \label{2.57}
\|e^{-tA}Pdiv\,F\|_{L^{r_2}(\mathbb{R}^n_+)}\leq Ct^{-\frac{1}{2}-\frac{n}{2}(\frac{1}{r_1}-\frac{1}{r_2})}\|F\|_{L^{r_1}(\mathbb{R}^n_+)}.
\end{equation}
By Lemma \ref{l:2.1}, for $t>0$ and $\varphi\in C^\infty_{0,\sigma}(\mathbb{R}^n_+)$, we have
\begin{equation} \label{}\notag
\begin{array}{rcl}
\displaystyle
|\langle e^{-tA} Pdiv\,F, \varphi\rangle|
&=&\displaystyle|\langle F, \nabla e^{-tA}\varphi\rangle|\vspace{2mm}\\
&\leq&\displaystyle\|\nabla e^{-tA}\varphi\|_{L^\frac{r_1}{r_1-1}(\mathbb{R}^n_+)}\|F\|_{L^{r_1}(\mathbb{R}^n_+)}\vspace{2mm}\\
&\leq&\displaystyle Ct^{-\frac{1}{2}-\frac{n}{2}(\frac{1}{r_1}-\frac{1}{r_2})}\|F\|_{L^{r_1}(\mathbb{R}^n_+)}\|\varphi\|_{L^{\frac{r_2}{r_2-1}}(\mathbb{R}^n_+)},
\vspace{2mm}\\
\end{array}
\end{equation}
which yields (\ref{2.57}) is true.
By Theorem 1.1, $\|u(t)\|_{L^2(\mathbb{R}^n_+)}\leq C_\epsilon(1+t)^{\frac{n-2}{4}-\epsilon}<\infty$ for $t>0$.
Using (\ref{2.57}), and Theorem \ref{th:4}, we get for $t>0$
\begin{equation} \label{}\notag
\begin{array}{rcl}
\displaystyle
\|\int_0^t e^{-(t-s)A} Pu(s)\cdot\nabla u(s)ds\big\|_{L^r(\mathbb{R}^n_+)}
&\leq&\displaystyle\int_0^t \|e^{-(t-s)A} Pdiv\,(u\otimes u)(s)\|_{L^r(\mathbb{R}^n_+)}ds\vspace{2mm}\\
&\leq&\displaystyle C\int_0^t(t-s)^{-\frac{1}{2}}\|u(s)\|^2_{L^{2r}(\mathbb{R}^n_+)}ds\vspace{2mm}\\

&\leq&\displaystyle C\int_0^t(t-s)^{-\frac{1}{2}}\|u(s)\|^{2(1-\frac{1}{r})}_{L^\infty(\mathbb{R}^n_+)}\|u(s)\|^{\frac{2}{r}}_{L^2(\mathbb{R}^n_+)}ds
\vspace{2mm}\\
\end{array}
\end{equation}
\begin{equation} \label{}\notag
\begin{array}{rcl}
&\leq&\displaystyle C\epsilon_0^{2(1-\frac{1}{r})}(\int_0^\frac{t}{2}+\int_\frac{t}{2}^t)(t-s)^{-\frac{1}{2}}s^{-1+\frac{1}{r}}ds
\vspace{2mm}\\
&\leq&\displaystyle C\epsilon_0^{2(1-\frac{1}{r})}t^{\frac{1}{r}-\frac{1}{2}},\;\;\;\;\frac{n}{n-1}\leq r<2.
\end{array}
\end{equation}
\noindent {\bf  Uniform $L^r$ estimates for $u$.} For $t>0$, we  have
\begin{equation} \label{}\notag
u(t)=e^{-tA}a-\int_0^te^{-(t-s)A}P(u(s)\cdot\nabla
u(s)-\theta e_n)ds.
\end{equation}
Using Lemma \ref{l:2.1}, Theorem \ref{th:1} and (\ref{2.55}), we get for $t>0$
\begin{equation} \label{}\notag
\begin{array}{rcl}
\displaystyle\|u(t)\|_{L^r(\mathbb{R}^n_+)}
&\leq&\displaystyle\|e^{-tA}a\|_{L^r(\mathbb{R}^n_+)}+\|\int_0^t e^{-(t-s)A} Pu(s)\cdot\nabla u(s)ds\big\|_{L^r(\mathbb{R}^n_+)}\vspace{2mm}\\
&&\displaystyle+C\int_0^t\|\theta(s)\|^{\frac{2}{r}-1}_{L^1(\mathbb{R}^n_+)}\|\theta(s)\|^{2(1-\frac{1}{r})}_{L^2(\mathbb{R}^n_+)}ds\vspace{2mm}\\
&\leq&\displaystyle C\|a\|_{L^r(\mathbb{R}^n_+)}+C\epsilon_0^{2(1-\frac{1}{r})}t^{\frac{1}{r}-\frac{1}{2}}
+C\int_0^t (1+s)^{-\frac{1}{2}(\frac{2}{r}-1)-\frac{n+2}{2}(1-\frac{1}{r})}ds\vspace{2mm}\\
&\leq&\displaystyle C(\|a\|_{L^r(\mathbb{R}^n_+)}+\epsilon_0^{2(1-\frac{1}{r})}t^{\frac{1}{r}-\frac{1}{2}}+\log(1+t)),
\end{array}
\end{equation}
which implies
\begin{equation} \label{2.58}
L(r, r, t)=\sup\limits_{0<s\leq t}\|u(s)\|_{L^r(\mathbb{R}^n_+)}\leq C(t)<\infty.
\end{equation}
Here $\lim_{t\to \infty} C(t) = \infty$. We now use this time dependent bound to obtain a uniform bound. For the rest of the proof, let  $\frac{n}{n-1}\leq r<2$. The assumption $a\in L^n_\sigma(\mathbb{R}^n_+)\bigcap L^\frac{n}{n-1}(\mathbb{R}^n_+)$ yields $a\in L^r(\mathbb{R}^n_+)$.
Note that similar to the proof of (2.27), we also have $\|\theta(t)\|_{L^1(\mathbb{R}^n_+)}\leq C(1+t)^{-\frac{1}{2}}$, $t>0$. Whence by Lemma 2.1 and Theorem 1.1, we get for $t>0$
\begin{equation} \label{}\notag
\begin{array}{rcl}
&&\displaystyle
\int_0^t\|e^{-(t-s)A}P\theta(s) e_n\|_{L^\frac{n}{n-1}(\mathbb{R}^n_+)}ds\vspace{2mm}\\
&\leq&\displaystyle C\int_0^\frac{t}{2}(t-s)^{-\frac{n}{2}(1-\frac{n-1}{n})}\|\theta(s)\|_{L^1(\mathbb{R}^n_+)}ds
+C\int_\frac{t}{2}^t\|\theta(s)\|^{1-\frac{2}{n}}_{L^1(\mathbb{R}^n_+)}\|\theta(s)\|^\frac{2}{n}_{L^2(\mathbb{R}^n_+)}ds
\vspace{2mm}\\
&\leq&\displaystyle Ct^{-\frac{1}{2}}\int_0^\frac{t}{2}(1+s)^{-\frac{1}{2}}ds
+C\int_\frac{t}{2}^t(1+s)^{-\frac{1}{2}(1-\frac{2}{n})-\frac{n+2}{2n}}ds\vspace{2mm}\\
&\leq&\displaystyle C;
\end{array}
\end{equation}
and
\begin{equation} \label{}\notag
\begin{array}{rcl}
\displaystyle
\int_0^t\|e^{-(t-s)A}P\theta(s)e_n\|_{L^r(\mathbb{R}^n_+)}ds&\leq&\displaystyle
C\int_0^\infty\|\theta(s)\|^{\frac{2}{r}-1}_{L^1(\mathbb{R}^n_+)}\|\theta(s)\|^{2(1-\frac{1}{r})}_{L^2(\mathbb{R}^n_+)}ds\vspace{2mm}\\
&\leq&\displaystyle \int_0^\infty(1+s)^{-\frac{1}{2}(\frac{2}{r}-1)-\frac{n+2}{2}(1-\frac{1}{r})}ds\vspace{2mm}\\
&\leq&\displaystyle C.
\end{array}
\end{equation}
It follows from Lemma \ref{l:2.1}, Theorems \ref{th:1}, \ref{th:4} and (\ref{2.55}), (\ref{2.56}), (\ref{2.58}) that for $t>0$
\begin{equation} \label{}\notag
\begin{array}{rcl}
\displaystyle\|u(t)\|_{L^r(\mathbb{R}^n_+)}
&\leq&\displaystyle\|e^{-tA}a\|_{L^r(\mathbb{R}^n_+)}+\int_0^t\|e^{-(t-s)A}P\theta(s)e_n\|_{L^r(\mathbb{R}^n_+)}ds\vspace{2mm}\\
&&\displaystyle+C\int_0^t(t-s)^{-\frac{1}{2}}\|u(s)\|_{L^r(\mathbb{R}^n_+)}\|\nabla u(s)\|_{L^n(\mathbb{R}^n_+)}ds\vspace{2mm}\\
&\leq&\displaystyle C(1+\|a\|_{L^r(\mathbb{R}^n_+)})+C\epsilon_0L(r, r, t)\int_0^t(t-s)^{-\frac{1}{2}}s^{-\frac{1}{2}}ds\vspace{2mm}\\
&\leq&\displaystyle C(1+\|a\|_{L^r(\mathbb{R}^n_+)})+C_0\epsilon_0L(r, r, t).
\end{array}
\end{equation}
Thus,
\begin{equation} \label{2.59}
L(r, r, t)\leq C(1+\|a\|_{L^r(\mathbb{R}^n_+)})+C_0\epsilon_0L(r, r, t).
\end{equation}
Take $\epsilon_0>0$ such that $C_0\epsilon_0\leq\frac{1}{2}$ in (\ref{2.59}). Then from (\ref{2.58}) and (\ref{2.59}), one has for $\frac{n}{n-1}\leq r<2$ and $t>0$,
\begin{equation} \label{2.60}
\|u(t)\|_{L^r(\mathbb{R}^n_+)}\leq L(r, r)\leq 2C(1+\|a\|_{L^r(\mathbb{R}^n_+)}).
\end{equation}

\noindent {\bf $L^{\infty}$ estimate for $\theta$}. Let $m>n$. By Theorem \ref{th:4}, we conclude for $t>0$
\begin{equation} \label{}\notag
\begin{array}{rcl}
\displaystyle
\|\theta(t)\|_{L^\infty(\mathbb{R}^n_+)}&\leq&\displaystyle \|E(\frac{t}{2})\theta(\frac{t}{2})\|_{L^\infty(\mathbb{R}^n_+)}+C\int_\frac{t}{2}^t\|\nabla G_{t-s}\|_{L^\frac{m}{m-1}(\mathbb{R}^n)}\|u(s)\theta(s)\|_{L^m(\mathbb{R}^n_+)}ds
\vspace{2mm}\\
&\leq&\displaystyle C t^{-\frac{n}{2}}\|\theta(\frac{t}{2})\|_{L^1(\mathbb{R}^n_+)}
+C\int_\frac{t}{2}^t(t-s)^{-\frac{1}{2}-\frac{n}{2m}}\|u(s)\|_{L^m(\mathbb{R}^n_+)}
\|\theta(s)\|_{L^\infty(\mathbb{R}^n_+)}ds
\vspace{2mm}\\
&\leq&\displaystyle Ct^{-\frac{n+1}{2}}+C\epsilon_0\Phi(t)\int_\frac{t}{2}^t(t-s)^{-\frac{1}{2}-\frac{n}{2m}}
s^{-\frac{n+1}{2}-\frac{n}{2}(\frac{1}{n}-\frac{1}{m})}ds\vspace{2mm}\\
&\leq&\displaystyle Ct^{-\frac{n+1}{2}}+C\epsilon_0\Phi(t) t^{-\frac{n+1}{2}},
\end{array}
\end{equation}
from which,
\begin{equation} \label{2.61}
\Phi(t)\leq C+C\epsilon_0\Phi(t), \quad where\quad \Phi(t)=\sup\limits_{0<s\leq t}(s^{\frac{n+1}{2}}\|\theta(s)\|_{L^\infty(\mathbb{R}^n_+)}).
\end{equation}
Take $\epsilon_0>0$ suitably small such that $C\epsilon_0\leq\frac{1}{2}$ in (\ref{2.61}). Then we get for $t>0$
\begin{equation} \label{2.62}
\|\theta(t)\|_{L^\infty(\mathbb{R}^n_+)}\leq Ct^{-\frac{n+1}{2}}.
\end{equation}
\noindent{\bf Decay for $\theta \in L^q, 1<q< \infty $}. By  (\ref{2.55}) and  (\ref{2.62}), we have the
straightforward estimate for $t>0$:
\begin{equation} \label{2.63}
\|\theta(t)\|_{L^q(\mathbb{R}^n_+)} \leq \|\theta(t)\|_{L^1(\mathbb{R}^n_+)}^{\frac{1}{q}}\|\theta(t)\|_{L^\infty(\mathbb{R}^n_+)}^\frac{q-1}{q}
\leq Ct^{-\frac{(n+1)(q-1)}{2q}}t^{-\frac{1}{2q}}=Ct^{-\frac{n+1}{2}+\frac{n}{2q}}=Ct^{-\frac{1}{2}-\frac{n}{2}(1-\frac{1}{q})}.
\end{equation}
\noindent{\bf Decay for $u \in L^q,$ $\frac{n}{n-1}\leq r<2$, $r\leq q\leq\infty$.}
By (\ref{2.63}), we have the auxiliary estimates for $k=0,1$ and any $t>0$
\begin{equation} \label{2.64}
\begin{array}{rcl}
\displaystyle\int_\frac{t}{2}^t (t-s)^{-\frac{k}{2}}\|\theta(s)\|_{L^q(\mathbb{R}^n_+)}ds
&\leq&\displaystyle C\int_\frac{t}{2}^t (t-s)^{-\frac{k}{2}}s^{-\frac{1}{2}-\frac{n}{2}(1-\frac{1}{q})}ds\\
&=&\displaystyle Ct^{-\frac{k}{2}-\frac{n}{2}(\frac{n-1}{n}-\frac{1}{r})-\frac{n}{2}(\frac{1}{r}-\frac{1}{q})}\vspace{2mm}\\
&\leq&\displaystyle Ct^{-\frac{k}{2}-\frac{n}{2}(\frac{1}{r}-\frac{1}{q})}.
\end{array}
\end{equation}
\noindent These estimates combined with  Lemma \ref{l:2.1}, Theorem \ref{th:4} and (\ref{2.60}) yield the decays for $u$ and $\nabla u$ in $L^q$. Let $\frac{n}{n-1}\leq r<2, \;r\leq q\leq \infty$, then for $k=0, 1$ and any $t>0$
\begin{equation} \label{2.65}
\begin{array}{rcl}
&&\displaystyle\|\nabla^k u(t)\|_{L^q(\mathbb{R}^n_+)}\vspace{2mm}\\
&\leq&\displaystyle\|\nabla^k e^{-\frac{t}{2}A}u(\frac{t}{2})\|_{L^q(\mathbb{R}^n_+)}+\int_\frac{t}{2}^t\|\nabla^k e^{-(t-s)A}P(u(s)\cdot\nabla
u(s)-\theta e_n)\|_{L^q(\mathbb{R}^n_+)}ds\vspace{2mm}\\
&\leq&\displaystyle Ct^{-\frac{k}{2}-\frac{n}{2}(\frac{1}{r}-\frac{1}{q})}\|u(\frac{t}{2})\|_{L^r(\mathbb{R}^n_+)}
+C\int_\frac{t}{2}^t(t-s)^{-\frac{k}{2}}\|\theta(s)\|_{L^q(\mathbb{R}^n_+)}ds\vspace{2mm}\\
&&\displaystyle
+C\int_\frac{t}{2}^t(t-s)^{-\frac{k}{2}-\frac{n}{2}(\frac{\beta}{n}+\frac{1}{q}-\frac{1}{q})}\|u(s)\|_{L^q(\mathbb{R}^n_+)}
\|\nabla u(s)\|_{L^\frac{n}{\beta}(\mathbb{R}^n_+)}ds,\;\;\;\;0<\beta\leq1.
\vspace{2mm}\\
\end{array}
\end{equation}
Let $k=0$ and $\beta=1$ in (\ref{2.65}). Combining the decay estimate for $\|\nabla u\|_{L^n(\mathbb{R}^n_+)}$ in Theorem \ref{th:4} yields
\begin{equation} \label{2.66}
\begin{array}{rcl}
\displaystyle\| u(t)\|_{L^q(\mathbb{R}^n_+)}
&\leq&\displaystyle Ct^{-\frac{n}{2}(\frac{1}{r}-\frac{1}{q})}+C\epsilon_0K(t)\int_\frac{t}{2}^t(t-s)^{-\frac{1}{2}}
s^{-\frac{1}{2}-\frac{n}{2}(\frac{1}{r}-\frac{1}{q})}ds
\vspace{2mm}\\
&\leq&\displaystyle Ct^{-\frac{n}{2}(\frac{1}{r}-\frac{1}{q})}+C\epsilon_0K(t)t^{-\frac{n}{2}(\frac{1}{r}-\frac{1}{q})},
\end{array}
\end{equation}
which implies $K(t)\leq C+C\epsilon_0K(t)$, with $\,K(t)=\sup\limits_{0<s\leq t}\{s^{\frac{n}{2}(\frac{1}{r}-\frac{1}{q})}\|u(s)\|_{L^q(\mathbb{R}^n_+)}\}.$
Taking $\epsilon_0>0$ such that $C\epsilon_0\leq\frac{1}{2}$, then for $t>0$
\begin{equation} \label{2.67}
\|u(t)\|_{L^q(\mathbb{R}^n_+)}\leq Ct^{-\frac{n}{2}(\frac{1}{r}-\frac{1}{q})}.
\end{equation}
\noindent Let $k=1$ and $\beta=\frac{1}{2}$ in (\ref{2.65}). (\ref{2.67}), (\ref{2.65}) and Theorem \ref{th:4} yield
\begin{equation} \label{2.68}
\|\nabla u(t)\|_{L^q(\mathbb{R}^n_+)}\leq  Ct^{-\frac{1}{2}-\frac{n}{2}(\frac{1}{r}-\frac{1}{q})}+C\int_\frac{t}{2}^t(t-s)^{-\frac{1}{2}-\frac{1}{4}}
s^{-\frac{1}{2}-\frac{1}{4}-\frac{n}{2}(\frac{1}{r}-\frac{1}{q})}ds\leq  Ct^{-\frac{1}{2}-\frac{n}{2}(\frac{1}{r}-\frac{1}{q})}.
\end{equation}
This completes the proof of the Proposition.\\
\noindent{\bf Proof of Theorem \ref{th:5}}. By (\ref{2.55}), (\ref{2.62}) and (\ref{2.63}), we have the decay estimate for $\theta \in L^q(\mathbb{R}^n_+)$, $1\leq q\leq\infty$. The decay for $\nabla\theta \in L^q(\mathbb{R}^n_+)$ is obtained as follows: Combining Theorem \ref{th:4} and (\ref{2.63}) yields for $1\leq q\leq\infty$ and any $t>0$
\begin{equation} \label{}\notag
\begin{array}{rcl}
\displaystyle\|\nabla\theta(t)\|_{L^q(\mathbb{R}^n_+)}
&\leq&\displaystyle\|\nabla E(\frac{t}{2})\theta(\frac{t}{2})\|_{L^q(\mathbb{R}^n_+)}+\int^t_\frac{t}{2}\big\|\nabla E(t-s)(u(s)\cdot\nabla\theta(s))\|_{L^q(\mathbb{R}^n_+)}ds\vspace{2mm}\\
&\leq&\displaystyle Ct^{-\frac{1}{2}}
\|\theta(\frac{t}{2})\|_{L^q(\mathbb{R}^n_+)}
+C\int^t_\frac{t}{2}(t-s)^{-\frac{3}{4}}
\|u(s)\|_{L^{2n}(\mathbb{R}^n_+)}\|\nabla\theta(s)\|_{L^q(\mathbb{R}^n_+)}ds\\
&\leq&\displaystyle Ct^{-1-\frac{n}{2}(1-\frac{1}{q})}
+C\epsilon_0\Theta_1(t)\int^t_\frac{t}{2}(t-s)^{-\frac{3}{4}}
s^{-\frac{5}{4}-\frac{n}{2}(1-\frac{1}{q})}ds\\
&\leq&\displaystyle Ct^{-1-\frac{n}{2}(1-\frac{1}{q})}
+C\epsilon_0\Theta_1(t) t^{-1-\frac{n}{2}(1-\frac{1}{q})},
\end{array}
\end{equation}
which yields
\begin{equation} \label{2.69}
\Theta_1(t)\leq C+C\epsilon_0\Theta_1(t),\quad where\quad \Theta_1(t)=\sup\limits_{0<s\leq t}\{s^{1+\frac{n}{2}(1-\frac{1}{q})}\|\nabla\theta(s)\|_{L^q(\mathbb{R}^n_+)}\}.
\end{equation}
Take $\epsilon_0>0$ in (\ref{2.69}) such that $C\epsilon_0\leq\frac{1}{2}$. Then for $1\leq q\leq\infty$ and any $t>0$
\begin{equation} \label{2.70}
\|\nabla\theta(t)\|_{L^q(\mathbb{R}^n_+)}\leq Ct^{-1-\frac{n}{2}(1-\frac{1}{q})}.
\end{equation}
Let $r=\frac{n}{n-1}$ in (\ref{2.67}) and (\ref{2.68}) respectively. We find for $\frac{n}{n-1}\leq q\leq\infty$ and any $t>0$
\begin{equation} \label{2.71}
\|\nabla^ku(t)\|_{L^q(\mathbb{R}^n_+)}\leq Ct^{-\frac{k-1}{2}-\frac{n}{2}(1-\frac{1}{q})},\;\;\;\;k=0, 1.
\end{equation}
From (\ref{2.55}), (\ref{2.62}), (\ref{2.63}), (\ref{2.70}) and (\ref{2.71}), the proof of the theorem is complete.  $\,\,\,\,\,\,\Box$

\section{Decay estimates for the second spatial order derivatives }\label{3}

\setcounter{equation}{0}

We first establish an auxiliary lemma that gives the proof of the first part of the Theorem \ref{th:6}.
\begin{lemma} \label{l:3.1}Let $a\in L^n_\sigma(\mathbb{R}^n_+)\bigcap L^\frac{n}{n-1}(\mathbb{R}^n_+)$ and $b\in L^1(\mathbb{R}^n_+)\bigcap L^2(\mathbb{R}^n_+)$ $\bigcap L^\frac{n}{3}(\mathbb{R}^n_+)$, $n\geq3$. Let $(u, p, \theta)$ be the strong solution of (\ref{1.1}) obtained in Theorem \ref{th:4}. Then for $t>0$
\begin{equation} \label{}\notag
\begin{array}{rcl}
&&\displaystyle
\|\nabla^2u(t)\|_{L^r (\mathbb{R}^n_+)}+\|Au(t)\|_{L^r(\mathbb{R}^n_+)}+\|\partial_tu(t)\|_{L^r(\mathbb{R}^n_+)}+\|\nabla p(t)\|_{L^r(\mathbb{R}^n_+)}\vspace{2mm}\\
&\leq&\displaystyle Ct^{-\frac{1}{2}-\frac{n}{2}(1-\frac{1}{r})}(1+t^{-n+\frac{3}{2}})\quad for\quad \frac{n}{n-1}\leq r<\infty;
\end{array}
\end{equation}
and
\begin{equation} \label{}\notag
\|\nabla^2\theta(t)\|_{L^r (\mathbb{R}^n_+)}\leq Ct^{-\frac{3}{2}-\frac{n}{2}(1-\frac{1}{r})}(1+t^{-\frac{n-2}{2}})\quad for \quad 1\leq r\leq\infty.
\end{equation}
\end{lemma}
{\bf Proof.}\quad Let $0<\alpha<1$ and $0<\delta<1-\alpha$. We first show that for  $1<r<\infty$ and  $t>0, \,\,h>0$
\begin{equation} \label{3.1}
\|A^\alpha u(t+h)-A^\alpha u(t)\|_{L^r(\mathbb{R}^n_+)}\leq C(h^\delta t^{-\alpha-\delta-\frac{n}{2}+\frac{1}{2}+\frac{n}{2r}}
+h^{1-\alpha}t^{-\frac{1}{2}-\frac{n}{2}+\frac{n}{2r}})(1+t^{-\frac{n-2}{2}}).
\end{equation}
Here it should be pointed out that (3.1) is crucial in estimating $J_4$ below. \\
Observe that for $t>0$, $h>0$
\begin{equation} \label{3.2}
\begin{array}{rcl}
\displaystyle
u(t+h)&=&\displaystyle e^{-(t+h-\frac{t}{2})A}u(\frac{t}{2})-\int_\frac{t}{2}^{t+h} e^{-(t+h-s)A}(Pu(s)\cdot\nabla u(s)-Pe_n\theta(s))ds\vspace{2mm}\\
&=&\displaystyle e^{-hA}e^{-\frac{t}{2}A}u(\frac{t}{2})-\big(\int_t^{t+h}+\int_\frac{t}{2}^t\big)e^{-(t-s)A}e^{-hA}(Pu(s)\cdot\nabla u(s)-Pe_n\theta(s))ds.
\end{array}
\end{equation}
Note that for any $\varphi\in D(A^\delta)$, $h\geq0$ and $1<q<\infty$
\begin{equation} \label{}\notag
\begin{array}{rcl}
\displaystyle\|(e^{-hA}-I)\varphi\|_{L^q(\mathbb{R}^n_+)}
=\displaystyle h\|\int_0^1A^{1-\delta}e^{-shA}A^\delta\varphi ds\|_{L^q(\mathbb{R}^n_+)}
&\leq&\displaystyle Ch\int_0^1(sh)^{\delta-1}ds\|A^\delta\varphi\|_{L^q(\mathbb{R}^n_+)}\vspace{2mm}\\
&\leq&\displaystyle\frac{C}{\delta}h^\delta\|A^\delta\varphi\|_{L^q(\mathbb{R}^n_+)}.
\end{array}
\end{equation}
Therefore, from (\ref{3.2}) and Lemma \ref{l:2.1}, Theorem \ref{th:5}, we conclude for all $1<r<\infty$ and $t>0, \,\,h>0$
\begin{equation} \label{}\notag
\begin{array}{rcl}
&&\displaystyle
\Hh_{\alpha}(h,t):=\|A^\alpha u(t+h)-A^\alpha u(t)\|_{L^r(\mathbb{R}^n_+)}\vspace{2mm}\\
&\leq&\displaystyle\|A^\alpha(e^{-hA}-I)e^{-\frac{t}{2}A}u(\frac{t}{2})\|_{L^r(\mathbb{R}^n_+)}
+\int_t^{t+h}\|A^\alpha e^{-(t+h-s)A}(Pu(s)\cdot\nabla u(s)-Pe_n\theta(s))\|_{L^r(\mathbb{R}^n_+)}ds\vspace{2mm}\\
&&\displaystyle+\int^t_\frac{t}{2}\|(e^{-hA}-I)A^\alpha e^{-(t-s)A}(Pu(s)\cdot\nabla u(s)-Pe_n\theta(s))\|_{L^r(\mathbb{R}^n_+)}ds.
\vspace{2mm}\\
\end{array}
\end{equation}
Hence,
\begin{equation} \label{}\notag
\begin{array}{rcl}
&&\displaystyle \Hh_{\alpha}(h,t)\vspace{2mm}\\
&\leq&\displaystyle Ch^\delta\|A^{\alpha+\delta}e^{-\frac{t}{2}A}u(\frac{t}{2})\|_{L^r(\mathbb{R}^n_+)}
\vspace{2mm}\\
&&\displaystyle+C\int_t^{t+h}(t+h-s)^{-\alpha}(\|u(s)\|_{L^{2r}(\mathbb{R}^n_+)}\|\nabla u(s)\|_{L^{2r}(\mathbb{R}^n_+)}+\|\theta(s)\|^\frac{1}{r}_{L^1(\mathbb{R}^n_+)}\|\theta(s)\|^{1-\frac{1}{r}}_{L^\infty(\mathbb{R}^n_+)})ds\vspace{2mm}\\
&&\displaystyle
+Ch^\delta\int^t_\frac{t}{2}(t-s)^{-\alpha-\delta}(\|u(s)\|_{L^{2r}(\mathbb{R}^n_+)}\|\nabla u(s)\|_{L^{2r}(\mathbb{R}^n_+)}+\|\theta(s)\|^\frac{1}{r}_{L^1(\mathbb{R}^n_+)}\|\theta(s)\|^{1-\frac{1}{r}}_{L^\infty(\mathbb{R}^n_+)})ds
\vspace{2mm}\\
\end{array}
\end{equation}
\begin{equation} \label{}\notag
\begin{array}{rcl}
&\leq&\displaystyle Ch^\delta t^{-\alpha-\delta}\|u(\frac{t}{2})\|_{L^r(\mathbb{R}^n_+)}
+C\int_t^{t+h}(t+h-s)^{-\alpha}
\big(s^{\frac{1}{2}-{n}(1-\frac{1}{2r})}+s^{-\frac{1}{2r}-\frac{n+1}{2}(1-\frac{1}{r})}\big)ds
\vspace{2mm}\\
&&\displaystyle
+Ch^\delta\int^t_\frac{t}{2}(t-s)^{-\alpha-\delta}
\big(s^{\frac{1}{2}-{n}(1-\frac{1}{2r})}+s^{-\frac{1}{2r}-\frac{n+1}{2}(1-\frac{1}{r})}\big)ds
\vspace{2mm}\\
&\leq&\displaystyle C(h^\delta t^{-\alpha-\delta-\frac{n}{2}+\frac{1}{2}+\frac{n}{2r}}
+h^{1-\alpha}(t^{-n+\frac{1}{2}+\frac{n}{2r}}+t^{-\frac{1}{2}-\frac{n}{2}+\frac{n}{2r}})+h^\delta (t^{-\alpha-\delta-n+\frac{3}{2}+\frac{n}{2r}}+t^{-\alpha-\delta-\frac{n}{2}+\frac{1}{2}+\frac{n}{2r}}))\vspace{2mm}\\
&\leq&\displaystyle C(h^\delta t^{-\alpha-\delta-\frac{n}{2}+\frac{1}{2}+\frac{n}{2r}}
+h^{1-\alpha}t^{-\frac{1}{2}-\frac{n}{2}+\frac{n}{2r}})(1+t^{-\frac{n-2}{2}}),
\end{array}
\end{equation}
this proves (\ref{3.1}). In particular for the sequel we will need the case: $\alpha=\frac{1}{2}$,
\begin{equation} \label{3.3}
\Hh_{\frac{1}{2}}(h,t)\leq C(h^\delta t^{-\delta-\frac{n}{2}+\frac{n}{2r}}
+h^\frac{1}{2}t^{-\frac{1}{2}-\frac{n}{2}+\frac{n}{2r}})(1+t^{-\frac{n-2}{2}}),\;\;\;0<\delta<\frac{1}{2},\;\;\;1<r<\infty.
\end{equation}
We decompose $Au(t) $ as follows.

It is not difficult to verify that it holds true for any $t>0$
\begin{equation} \label{3.4}
\begin{array}{rcl}
\displaystyle Au(t)&=&\displaystyle Ae^{-\frac{3t}{4}A}u(\frac{t}{4})-(I-e^{-\frac{t}{2}A})(P(u\cdot\nabla u)(t)-Pe_n\theta(t))\vspace{2mm}\\
&&\displaystyle-\int_\frac{t}{4}^\frac{t}{2}Ae^{-(t-s)A}(P(u\cdot\nabla u)(s)-Pe_n\theta(s))ds\vspace{2mm}\\
&&\displaystyle-\int_\frac{t}{2}^tAe^{-(t-s)A}(P(u\cdot\nabla u)(s)-P(u\cdot\nabla u)(t)+Pe_n\theta(t)-Pe_n\theta(s))ds\vspace{2mm}\\
&=&\displaystyle J_1(t)+J_2(t)+J_3(t)+J_4(t).
\end{array}
\end{equation}
We now obtain $L^r$ bounds for $J_i, i=1,2,3,4$. By Lemma \ref{l:2.1} and Theorem \ref{th:5} that for $\frac{n}{n-1}\leq r<\infty$ and each $t>0$
\begin{equation} \label{3.5}
\|J_1(t)\|_{L^r(\mathbb{R}^n_+)}\leq Ct^{-1}\|u(\frac{t}{4})\|_{L^r(\mathbb{R}^n_+)}\leq Ct^{-\frac{1}{2}-\frac{n}{2}(1-\frac{1}{r})};
\end{equation}
\begin{equation} \label{3.6}
\begin{array}{rcl}
\displaystyle
\|J_2(t)\|_{L^r(\mathbb{R}^n_+)}&\leq& C\|u(t)\|_{L^{2r}(\mathbb{R}^n_+)}\|\nabla u(t)\|_{L^{2r}(\mathbb{R}^n_+)}+C\|\theta(t)\|^\frac{1}{r}_{L^1(\mathbb{R}^n_+)}\|\theta(t)\|^{1-\frac{1}{r}}_{L^\infty(\mathbb{R}^n_+)}\vspace{2mm}\\
&\leq&\displaystyle Ct^{\frac{1}{2}-\frac{n}{2}(1-\frac{1}{2r})-\frac{n}{2}(1-\frac{1}{2r})}+Ct^{-\frac{1}{2r}-\frac{n+1}{2}(1-\frac{1}{r})}\vspace{2mm}\\
&\leq&\displaystyle Ct^{-n+\frac{1}{2}+\frac{n}{2r}}+Ct^{-\frac{1}{2}-\frac{n}{2}(1-\frac{1}{r})}\vspace{2mm}\\
&\leq&\displaystyle Ct^{-\frac{1}{2}-\frac{n}{2}(1-\frac{1}{r})}(1+t^{-\frac{n-2}{2}});
\end{array}
\end{equation}
\begin{equation} \label{3.7}
\begin{array}{rcl}
\displaystyle
\|J_3(t)\|_{L^r(\mathbb{R}^n_+)}&\leq&\displaystyle C\int_\frac{t}{4}^\frac{t}{2}(t-s)^{-1}(\|u(s)\|_{L^{2r}(\mathbb{R}^n_+)}\|\nabla u(s)\|_{L^{2r}(\mathbb{R}^n_+)}+\|\theta(s)\|^\frac{1}{r}_{L^1(\mathbb{R}^n_+)}\|\theta(s)\|^{1-\frac{1}{r}}_{L^\infty(\mathbb{R}^n_+)})ds\vspace{2mm}\\
&\leq&\displaystyle C\int_\frac{t}{4}^\frac{t}{2}(t-s)^{-1}(s^{-\frac{1}{2}-n+\frac{n}{2r}}+s^{-\frac{1}{2}-\frac{n}{2}+\frac{n}{2r}})ds\vspace{2mm}\\
&\leq&\displaystyle Ct^{-\frac{1}{2}-\frac{n}{2}+\frac{n}{2r}}(1+t^{-\frac{n-2}{2}}).
\end{array}
\end{equation}
Since we don't have a priori local bound for $\|\nabla^2\theta(t)\|_{L^q(\mathbb{R}^n_+)}$, we use (\ref{2.36}) to establish the decay of $\|\nabla^2\theta_j(t)\|_{L^q(\mathbb{R}^n_+)}$ by induction. Let $1\leq q\leq\infty$. Using the structure of the operator of $E(t)$, it is not difficult to verify that for any $t>0$
\begin{equation} \label{}\notag
\begin{array}{rcl}
\displaystyle
\|\nabla^2\theta_0(t)\|_{L^q(\mathbb{R}^n_+)}&=&\displaystyle\|\nabla^2E(t)b\|_{L^q(\mathbb{R}^n_+)}\vspace{2mm}\\
&=&\displaystyle\|\nabla^2_x\partial_{x_n}\int_{\mathbb{R}^n_+}\int_{-1}^1G_t(x^\prime-y^\prime, x_n-sy_n)y_nb(y)dsdy\|_{L^q(\mathbb{R}^n_+)}\vspace{2mm}\\
&\leq&\displaystyle C\|\nabla^3G_t\|_{L^q(\mathbb{R}^n)}\|y_nb\|_{L^1(\mathbb{R}^n_+)}\vspace{2mm}\\
&\leq&\displaystyle C_0t^{-\frac{3}{2}-\frac{n}{2}(1-\frac{1}{q})}.
\end{array}
\end{equation}
Assume that there exists $j\geq0$, such that
$$
\|\nabla^2\theta_j(t)\|_{L^q(\mathbb{R}^n_+)}\leq C_\ast t^{-\frac{3}{2}-\frac{n}{2}(1-\frac{1}{q})}.
$$
where $C_\ast\geq C_0$ is independent of $j$, which will be determined later.

From (\ref{2.36}), we know for any $j\geq0$ and $t>0$
\begin{equation} \label{}\notag
\theta_{j+1}(t)=E(\frac{t}{2})\theta_{j+1}(\frac{t}{2})-\int_\frac{t}{2}^tE(t-s)u_j(s)\cdot\nabla\theta_j(s)ds.
\end{equation}

Let $1\leq q\leq\infty$. It follows from Theorems \ref{th:4}, \ref{th:5} (where the decay estimates are also valid for $(u_j, \theta_j)$) that for any $t>0$
\begin{equation} \label{3.8}
\begin{array}{rcl}
&&\displaystyle
\|\nabla^2\theta_{j+1}(t)\|_{L^q(\mathbb{R}^n_+)}\vspace{2mm}\\
&\leq&\displaystyle\|\nabla^2E(\frac{t}{2})\theta_{j+1}(\frac{t}{2})\|_{L^q(\mathbb{R}^n_+)}+2\int_\frac{t}{2}^t\|\nabla G_{t-s}\|_{L^1(\mathbb{R}^n)}\|\nabla\big(u_j(s)\cdot\nabla \theta_j(s)\big)\|_{L^q(\mathbb{R}^n_+)}ds\vspace{2mm}\\
&\leq&\displaystyle Ct^{-1}\|\theta_{j+1}(\frac{t}{2})\|_{L^q(\mathbb{R}^n_+)}+C\int_\frac{t}{2}^t(t-s)^{-\frac{1}{2}}\|\nabla u_j(s)\|_{L^\infty(\mathbb{R}^n_+)}\|\nabla \theta_j(s)\|_{L^q(\mathbb{R}^n_+)}ds\vspace{2mm}\\
&&\displaystyle+C\int_\frac{t}{2}^t(t-s)^{-\frac{1}{2}}\|u_j(s)\|_{L^\infty(\mathbb{R}^n_+)}\|\nabla^2\theta_j(s)\|_{L^q(\mathbb{R}^n_+)}ds
\vspace{2mm}\\
&\leq&\displaystyle Ct^{-\frac{3}{2}-\frac{n}{2}(1-\frac{1}{q})}+C\int_\frac{t}{2}^t(t-s)^{-\frac{1}{2}}
s^{-2-\frac{n}{2}(1-\frac{1}{q})}ds\vspace{2mm}\\
&&\displaystyle+C\epsilon_0C_\ast\int_\frac{t}{2}^t(t-s)^{-\frac{1}{2}}s^{-2-\frac{n}{2}(1-\frac{1}{q})}ds
\vspace{2mm}\\
&\leq&\displaystyle C_1t^{-\frac{3}{2}-\frac{n}{2}(1-\frac{1}{q})}+C\epsilon_0C_\ast t^{-\frac{3}{2}-\frac{n}{2}(1-\frac{1}{q})}.
\end{array}
\end{equation}
Taking $\epsilon_0>0$ suitably small such that $C\epsilon_0\leq\frac{1}{2}$ in (\ref{3.8}), and take $C_\ast=C_0+2C_1$. Then for $1\leq q\leq\infty$ and $t>0$
\begin{equation} \label{}\notag
\|\nabla^2\theta_{j+1}(t)\|_{L^q(\mathbb{R}^n_+)}\leq C_\ast t^{-\frac{3}{2}-\frac{n}{2}(1-\frac{1}{q})}.
\end{equation}
From the above induction argument, we conclude that for any $t>0$
\begin{equation} \label{}\notag
\sup\limits_{j\geq0}\|\nabla^2\theta_j(t)\|_{L^q(\mathbb{R}^n_+)}\leq C_\ast t^{-\frac{3}{2}-\frac{n}{2}(1-\frac{1}{q})}.
\end{equation}
Since we have established the estimates of $\|\nabla^k\theta_j(t)\|_{L^q(\mathbb{R}^n_+)}$ ($k=0, 1$) in (2.45), (2.47), (2.48), and like using the weak convergence arguments in the proof of Theorem 1.5, we infer for each $t>0$ and function $\theta(t)$, as $j\longrightarrow\infty$
$$
\nabla^2\theta_j(t)\rightharpoonup \nabla^2\theta(t) \;\;\mbox{weakly in}\;\; L^q(\mathbb{R}^n_+), \;\;1<q\leq\infty.
$$
Moreover,
\begin{equation} \label{3.9}
\|\nabla^2\theta(t)\|_{L^q(\mathbb{R}^n_+)}\leq\liminf\limits_{j\longrightarrow\infty}\|\nabla^2\theta_j(t)\|_{L^q(\mathbb{R}^n_+)}\leq C_\ast t^{-\frac{3}{2}-\frac{n}{2}(1-\frac{1}{q})}, \;\;1<q\leq\infty.
\end{equation}
Now we deal with $\|\nabla^2\theta(t)\|_{L^1(\mathbb{R}^n_+)}$, $t>0$. Observe that
\begin{equation} \label{}\notag
\theta(t)=E(\frac{t}{2})\theta(\frac{t}{2})-\int_\frac{t}{2}^tE(t-s)u(s)\cdot\nabla\theta(s)ds.
\end{equation}
Like the proof of (3.8), using (3.9) with $q=2$, and Theorem 1.6, we have for $t>0$
\begin{equation} \label{}\notag
\begin{array}{rcl}
&&\displaystyle
\|\nabla^2\theta(t)\|_{L^1(\mathbb{R}^n_+)}\vspace{2mm}\\
&\leq&\displaystyle\|\nabla^2E(\frac{t}{2})\theta(\frac{t}{2})\|_{L^1(\mathbb{R}^n_+)}+2\int_\frac{t}{2}^t\|\nabla G_{t-s}\|_{L^1(\mathbb{R}^n)}\|\nabla\big(u(s)\cdot\nabla \theta(s)\big)\|_{L^1(\mathbb{R}^n_+)}ds\vspace{2mm}\\
&\leq&\displaystyle Ct^{-1}\|\theta(\frac{t}{2})\|_{L^1(\mathbb{R}^n_+)}+C\int_\frac{t}{2}^t(t-s)^{-\frac{1}{2}}\|\nabla u(s)\|_{L^\infty(\mathbb{R}^n_+)}\|\nabla \theta(s)\|_{L^1(\mathbb{R}^n_+)}ds\vspace{2mm}\\
&&\displaystyle+C\int_\frac{t}{2}^t(t-s)^{-\frac{1}{2}}\|u(s)\|_{L^n(\mathbb{R}^n_+)}\|\nabla^2\theta(s)\|_{L^\frac{n}{n-1}(\mathbb{R}^n_+)}ds
\vspace{2mm}\\
&\leq&\displaystyle Ct^{-\frac{3}{2}}+C\int_\frac{t}{2}^t(t-s)^{-\frac{1}{2}}
s^{-2}ds+C\int_\frac{t}{2}^t(t-s)^{-\frac{1}{2}}(1+s)^{\frac{1}{2}-\frac{n}{2}(1-\frac{1}{n})}s^{-\frac{3}{2}-\frac{n}{2}(1-\frac{n-1}{n})}ds
\vspace{2mm}\\
&\leq&\displaystyle C_1t^{-\frac{3}{2}}+C_2t^{-\frac{3}{2}}(1+t)^{-\frac{n-2}{2}}
\vspace{2mm}\\
&\leq&\displaystyle C_3t^{-\frac{3}{2}},
\end{array}
\end{equation}
from which, we know (3.9) is valid for $q=1$.\\
From (\ref{3.9}), Theorem \ref{th:5} and the equation on $\theta$ in (\ref{1.1}), we get for all $1\leq q\leq\infty$ and $t>0$
\begin{equation} \label{3.10}
\begin{array}{rcl}
\displaystyle
\|\partial_t\theta(t)\|_{L^q(\mathbb{R}^n_+)}&\leq&\displaystyle\|\Delta\theta(t)\|_{L^q(\mathbb{R}^n_+)}+\|u(t)\cdot\nabla \theta(t)\|_{L^q(\mathbb{R}^n_+)}
\vspace{2mm}\\
&\leq&\displaystyle Ct^{-\frac{3}{2}-\frac{n}{2}(1-\frac{1}{q})}+\|u(t)\|_{L^\infty(\mathbb{R}^n_+)}\|\nabla \theta(t)\|_{L^q(\mathbb{R}^n_+)}\vspace{2mm}\\
&\leq&\displaystyle Ct^{-\frac{3}{2}-\frac{n}{2}(1-\frac{1}{q})}(1+t^{-\frac{n-2}{2}}).
\end{array}
\end{equation}
It follows from (\ref{3.10}) that for  $1\leq q\leq\infty$ and  $0<s\leq t$
\begin{equation} \label{3.11}
\begin{array}{rcl}
\displaystyle\|\theta(t)-\theta(s)\|_{L^q(\mathbb{R}^n_+)}
\leq\displaystyle \int_s^t\|\partial_\tau\theta(\tau)\|_{L^q(\mathbb{R}^n_+)}d\tau
&\leq&\displaystyle  C\int_s^t\tau^{-\frac{3}{2}-\frac{n}{2}(1-\frac{1}{q})}(1+\tau^{-\frac{n-2}{2}})d\tau\vspace{2mm}\\
&\leq&\displaystyle  C(t-s)s^{-\frac{3}{2}-\frac{n}{2}(1-\frac{1}{q})}(1+s^{-\frac{n-2}{2}}).
\end{array}
\end{equation}
Note that for $1<q<n$
\begin{equation} \label{3.12}
\|\varphi\|_{L^\frac{nq}{n-q}(\mathbb{R}^n_+)}\leq C\|A^\frac{1}{2}\varphi\|_{L^q(\mathbb{R}^n_+)},\quad \forall \varphi\in D(A^\frac{1}{2});
\end{equation}
We also recall that  (see \cite{BM})
\begin{equation} \label{}\notag\|\nabla u(t)\|_{L^q(\mathbb{R}^n_+)}\approx\|A^\frac{1}{2}u(t)\|_{L^q(\mathbb{R}^n_+)},\,\,\,1<q<\infty.
\end{equation}
Using the above estimate, and combining (\ref{3.11}), (\ref{3.12}), Lemma \ref{l:2.1}, and Theorem \ref{th:5} yield for any $\frac{n}{n-1}\leq r<\infty$,  $t>0$
\begin{equation} \label{}\notag
\begin{array}{rcl}
\displaystyle\|J_4(t)\|_{L^r(\mathbb{R}^n_+)}
&\leq& \displaystyle C\int_\frac{t}{2}^t(t-s)^{-1}(\|Pu(s)\cdot\nabla (u(s)-u(t))\|_{L^r(\mathbb{R}^n_+)}+\|Pe_n(\theta(t)-\theta(s))\|_{L^r(\mathbb{R}^n_+)})ds\vspace{2mm}\\
&&\displaystyle+C\int_\frac{t}{2}^t(t-s)^{-1}\|e^{-\frac{(t-s)A}{2}}P(u(s)-u(t))\cdot\nabla u(t)\|_{L^r(\mathbb{R}^n_+)}ds\vspace{2mm}\\
&\leq&\displaystyle C\int_\frac{t}{2}^t(t-s)^{-1}\|u(s)\|_{L^{2r}(\mathbb{R}^n_+)}\|\nabla u(s)-\nabla u(t)\|_{L^{2r}(\mathbb{R}^n_+)}ds\vspace{2mm}\\
&&\displaystyle+C\int_\frac{t}{2}^t(t-s)^{-1}\|\theta(t)-\theta(s)\|_{L^r(\mathbb{R}^n_+)}ds
\vspace{2mm}\\
&&\displaystyle+C\int_\frac{t}{2}^t(t-s)^{-1-\frac{n}{2}(\frac{1}{q}-\frac{1}{r})}\|u(s)-u(t)\|_{L^\frac{nq}{n-q}(\mathbb{R}^n_+)}\|\nabla u(t)\|_{L^n(\mathbb{R}^n_+)}ds\vspace{2mm}\\
&=&\displaystyle I_1+I_2+I_3.
\end{array}
\end{equation}
Using (\ref{3.3}), (\ref{3.11}) and (\ref{3.12}), we estimate $I_i,\;\;i=1,2,3$ separately.
\begin{equation} \label{}\notag
\begin{array}{rcl}
\displaystyle
I_1 &\leq& \displaystyle C\int_\frac{t}{2}^t(t-s)^{-1} [(t-s)^\delta s^{-\delta-\frac{n}{2}+\frac{n}{4r}}+(t-s)^\frac{1}{2}s^{-\frac{1}{2}-\frac{n}{2}+\frac{n}{4r}}]
s^{\frac{1}{2}-\frac{n}{2}(1-\frac{1}{2r})}
(1+s^{-\frac{n-2}{2}})ds \vspace{2mm}\\
&\leq&\displaystyle  Ct^{-\frac{1}{2}-\frac{n}{2}(1-\frac{1}{r})}(1+t^{-n+2});
\end{array}
\end{equation}
\begin{equation} \label{}\notag
I_2 \leq   C\int_\frac{t}{2}^t s^{-\frac{3}{2}-\frac{n}{2}(1-\frac{1}{r})}(1+s^{-\frac{n-2}{2}})ds \leq t^{-\frac{1}{2}-\frac{n}{2}(1-\frac{1}{r})}(1+t^{-\frac{n-2}{2}});
\end{equation}
\begin{equation} \label{}\notag
\begin{array}{rcl}
\displaystyle I_3 &\leq& \displaystyle Ct^{-\frac{n}{2}(1-\frac{1}{n})}\int_\frac{t}{2}^t((t-s)^{\delta-1-\frac{n}{2}(\frac{1}{q}-\frac{1}{r})} s^{-\delta-\frac{n}{2}+\frac{n}{2q}}
+(t-s)^{-\frac{1}{2}-\frac{n}{2}(\frac{1}{q}-\frac{1}{r})}s^{-\frac{1}{2}-\frac{n}{2}+\frac{n}{2q}})(1+s^{-\frac{n-2}{2}}) ds\vspace{2mm}\\
&\leq&\displaystyle Ct^{-\frac{n}{2}(1-\frac{1}{n})}(1+t^{-\frac{n-2}{2}})\vspace{2mm}\\
&&\displaystyle\times \Big(t^{-\delta-\frac{n}{2}+\frac{n}{2q}}\int_\frac{t}{2}^t(t-s)^{\delta-1-\frac{n}{2}(\frac{1}{q}-\frac{1}{r})}ds
+t^{-\frac{1}{2}-\frac{n}{2}+\frac{n}{2q}}\int_\frac{t}{2}^t(t-s)^{-\frac{1}{2}-\frac{n}{2}(\frac{1}{q}-\frac{1}{r})}ds\Big).
\end{array}
\end{equation}
To proceed, we need to guarantee the following two integrals
$$
\int_\frac{t}{2}^t(t-s)^{\delta-1-\frac{n}{2}(\frac{1}{q}-\frac{1}{r})}ds<\infty\;\;\;\mbox{and}\;\;\;
\;\;\int_\frac{t}{2}^t(t-s)^{-\frac{1}{2}-\frac{n}{2}(\frac{1}{q}-\frac{1}{r})}ds<\infty.
$$
To do this, $\delta-\frac{n}{2}(\frac{1}{q}-\frac{1}{r})>0$ and $\frac{1}{2}-\frac{n}{2}(\frac{1}{q}-\frac{1}{r})>0$ are required. Whence
we take $q=r$ if $r\in[\frac{n}{n-1}, n)$. If $n\leq r<\infty$, one can find $\delta\in(0, \frac{1}{2})$ such that $n\leq r<\frac{n}{1-2\delta}$, which is equivalent to $\frac{1}{r}+\frac{2\delta}{n}>\frac{1}{n}$. Then we take a number $q\in(\frac{n}{n-1}, n)$ such that $\frac{1}{r}+\frac{2\delta}{n}>\frac{1}{q}>\frac{1}{n}$, which implies $\delta-\frac{n}{2}(\frac{1}{q}-\frac{1}{r})>0$, and then $\frac{1}{2}-\frac{n}{2}(\frac{1}{q}-\frac{1}{r})>0$ due to the choice of $\delta\in(0,\frac{1}{2})$. From these arguments, we get for $t>0$
$$
\begin{array}{rcl}
\displaystyle I_3 &\leq& \displaystyle  Ct^{-\frac{n}{2}(1-\frac{1}{n})}(1+t^{-\frac{n-2}{2}})\big(t^{-\delta-\frac{n}{2}+\frac{n}{2q}+\delta-\frac{n}{2}(\frac{1}{q}-\frac{1}{r})}
+t^{-\frac{1}{2}-\frac{n}{2}+\frac{n}{2q}+\frac{1}{2}-\frac{n}{2}(\frac{1}{q}-\frac{1}{r})}\big)\vspace{2mm}\\
&\leq&\displaystyle Ct^{-\frac{n-1}{2}-\frac{n}{2}(1-\frac{1}{r})}(1+t^{-\frac{n-2}{2}}).
\end{array}
$$
Combining the estimates for $I_i, i=1,2,3$ yields for $t>0$
\begin{equation} \label{3.13}
\|J_4(t)\|_{L^r(\mathbb{R}^n_+)} \leq Ct^{-\frac{1}{2}-\frac{n}{2}(1-\frac{1}{r})}(1+t^{-n+2}).
\end{equation}
Combining the estimates (\ref{3.5})--(\ref{3.7}), (\ref{3.13}) for $\|J_i(t)\|_{L^r(\mathbb{R}^n_+)}$, $i=1,2,3,4$ with (\ref{3.4}) gives
for $t>0$
\begin{equation} \label{3.14}
\|Au(t)\|_{L^r(\mathbb{R}^n_+)}\leq Ct^{-\frac{1}{2}-\frac{n}{2}(1-\frac{1}{r})}(1+t^{-n+2})\quad for\quad \frac{n}{n-1}\leq r<\infty.
\end{equation}
Note that (see \cite{BM}): $\|\nabla^2 u(t)\|_{L^q(\mathbb{R}^n_+)}\leq C\|Au(t)\|_{L^q(\mathbb{R}^n_+)}$ for each $q\in (1, \infty)$.
Hence from (\ref{3.14}), one has for $t>0$
\begin{equation} \label{3.15}
\|\nabla^2 u(t)\|_{L^r(\mathbb{R}^n_+)}\leq Ct^{-\frac{1}{2}-\frac{n}{2}(1-\frac{1}{r})}(1+t^{-n+2})\quad for\quad \frac{n}{n-1}\leq r<\infty.
\end{equation}
Since $\partial_tu(t)=-Au(t)-P\big((u\cdot\nabla u)(t)-e_n\theta(t)\big)$, from (\ref{3.14}), (\ref{3.15}) and Theorem \ref{th:5}, we infer for $t>0$
\begin{equation} \label{3.16}
\begin{array}{rcl}
\displaystyle
\|\partial_tu(t)\|_{L^r(\mathbb{R}^n_+)}&\leq&\displaystyle\|Au(t)\|_{L^r(\mathbb{R}^n_+)}+\|u(t)\|_{L^{2r}(\mathbb{R}^n_+)}\|\nabla u(t)\|_{L^{2r}(\mathbb{R}^n_+)}+\|\theta(t)\|_{L^r(\mathbb{R}^n_+)}\vspace{2mm}\\
&\leq&\displaystyle Ct^{-\frac{1}{2}-\frac{n}{2}(1-\frac{1}{r})}(1+t^{-n+2})\quad for\quad  \frac{n}{n-1}\leq  r<\infty.
\end{array}
\end{equation}
Then it follows from Theorem \ref{th:5} and (\ref{3.14})--(\ref{3.16}) that for $t>0$
\begin{equation} \label{3.17}
\begin{array}{rcl}
\displaystyle
\|\nabla p(t)\|_{L^r(\mathbb{R}^n_+)}&\leq&\displaystyle \|\partial_tu(t)\|_{L^r(\mathbb{R}^n_+)}+\|\nabla^2 u(t)\|_{L^r(\mathbb{R}^n_+)}\vspace{2mm}\\
&&\displaystyle +\|u(t)\|_{L^{2r}(\mathbb{R}^n_+)}\|\nabla u(t)\|_{L^{2r}(\mathbb{R}^n_+)}+\|\theta(t)\|_{L^r(\mathbb{R}^n_+)}\vspace{2mm}\\
&\leq&\displaystyle Ct^{-\frac{1}{2}-\frac{n}{2}(1-\frac{1}{r})}(1+t^{-n+2})\quad for\quad  \frac{n}{n-1}\leq r<\infty.
\end{array}
\end{equation}
From (\ref{3.9}), (\ref{3.14})--(\ref{3.17}), we complete the proof of Lemma \ref{l:3.1}. $\,\,\,\,\Box$\\
We now prove the rest of the estimates for the conclusion of Theorem \ref{th:6}. We start with some auxiliary estimates.
Let  \begin{equation} \label{3.18}
\mathcal{N}=\int_0^\infty F(\tau)d\tau,
\end{equation}
where the operator $F$ is defined: $F(t)f(x)=\int_{\mathbb{R}^n_+}[G_t(x^\prime-y^\prime,
x_n-y_n)+G_t(x^\prime-y^\prime, x_n+y_n)]f(y)dy$.
Then  $g=\mathcal{N}f$ is the solution to the Neumann problem
\begin{equation} \label{}\notag
-\Delta g=f \;\;in \;\; \mathbb{R}^n_+,\;\;\;\partial_\nu g\,|_{\partial\mathbb{R}^n_+}=0.
\end{equation}
(see \cite{Han}). An easy calculation gives for  $u\in C^\infty_{0, \sigma}(\mathbb{R}^n_+)$ and a scalar function $\theta\in C^\infty_0(\mathbb{R}^n_+)$
\begin{equation} \label{3.19}
P(u\cdot\nabla u)=u\cdot\nabla u+\sum\limits_{i,
j=1}^n\nabla\mathcal{N}\partial_i\partial_j(u_iu_j);
\end{equation}
and
\begin{equation} \label{3.20}
P(e_n\theta)=e_n\theta+\nabla\mathcal{N}div\,(e_n\theta).
\end{equation}
\begin{lemma}\label{1:3.2}  Let $0\leq\eta<1$ and $1\leq k\leq n$. Then for any $1\leq q\leq\infty$
\begin{equation} \label{3.21}
\|\sum\limits_{i,j=1}^nx_n^{-\eta}\partial_k\mathcal{N}\partial_i\partial_j(u_iu_j)\|_{L^q(\mathbb{R}^n_+)}\leq C(\|u\|^2_{L^{2q}(\mathbb{R}^n_+)}+\|\nabla u\|^2_{L^{2q}(\mathbb{R}^n_+)}),\quad\forall\,u\in C^\infty_{0,\sigma}(\mathbb{R}^n_+);
\end{equation}
and
\begin{equation} \label{3.22}
\|x_n^{-\eta}\partial_k\mathcal{N}div\,(e_n\theta)\|_{L^q(\mathbb{R}^n_+)}\leq C(\|y_n\theta\|_{L^q(\mathbb{R}^n_+)}+\|\partial_n \theta\|_{L^q(\mathbb{R}^n_+)}),\quad\forall\,\theta\in C^\infty_0(\mathbb{R}^n_+).
\end{equation}
\end{lemma}
{\bf Proof.}\quad The proof of (\ref{3.21}) is given in \cite{Han} with $\eta=0$, and in \cite{Han2} with $0<\eta<1$ respectively. It remains to prove (\ref{3.22}). Denote the odd and even extensions of a function $f$ from $\mathbb{R}^n_+$ to $\mathbb{R}^n$, respectively
by
\begin{equation} \label{}\notag
f^\ast(x^\prime, x_n)=\left\{
\begin{array}{lll}
\displaystyle f(x^\prime, x_n) & if &  x_n\geq0 ,\vspace{2mm}\\
\displaystyle -f(x^\prime, -x_n) & if & x_n<0,
\end{array}
\right.
\quad
and
\quad \label{}\notag
f_\ast(x^\prime, x_n)=\left\{
\begin{array}{lll}
\displaystyle f(x^\prime, x_n) & if &  x_n\geq0 ,\vspace{2mm}\\
\displaystyle f(x^\prime, -x_n) & if & x_n<0.
\end{array}
\right.
\end{equation}
Note that for any $x\in\mathbb{R}^n$ and $\tau>0$,
\begin{equation} \label{3.23}
\begin{array}{rcl}
\displaystyle G_\tau\ast[div\,(e_n\theta)]_\ast(x)&=&\displaystyle\int_{\mathbb{R}^n}G_\tau(x^\prime-y^\prime, x_n-y_n)(\partial_n\theta)_\ast(y) dy\vspace{2mm}\\
&=&\displaystyle\int_{\mathbb{R}^n}G_\tau(x^\prime-y^\prime, x_n-y_n)\partial_{y_n}\theta^\ast(y) dy\vspace{2mm}\\
&=&\displaystyle\partial_n\int_{\mathbb{R}^n}G_\tau(x^\prime-y^\prime, x_n-y_n)\theta^\ast(y) dy\vspace{2mm}\\
&=&\displaystyle\partial_n\int_{\mathbb{R}^n}(G_\tau(x^\prime-y^\prime, x_n-y_n)-G_\tau(x^\prime-y^\prime, x_n))\theta^\ast(y) dy\vspace{2mm}\\
&=&\displaystyle\partial_n\int_{\mathbb{R}^n}\int_0^1\frac{d}{ds}G_\tau(x^\prime-y^\prime, x_n-sy_n)ds\theta^\ast(y) dy\vspace{2mm}\\
&=&\displaystyle\partial_n\int_{\mathbb{R}^n}\int_0^1\partial_{x_n}G_\tau(x^\prime-y^\prime, x_n-sy_n)(-y_n)\theta^\ast(y) dsdy\vspace{2mm}\\
&=&\displaystyle\partial_{x_n}\partial_{x_n}\int_{\mathbb{R}^n}\int_0^1G_\tau(x^\prime-y^\prime, x_n-sy_n)(-y_n)\theta^\ast(y) dsdy\vspace{2mm}\\
\end{array}
\end{equation}
Let $1\leq q\leq\infty$. It follows from (\ref{3.18}) and (\ref{3.23}) that for any $1\leq k\leq n$.
\begin{equation} \label{}\notag
\begin{array}{rcl}
&&\displaystyle\|x_n^{-\eta}\partial_k\mathcal{N}div\,(e_n\theta)\|_{L^q(\mathbb{R}^n_+)}\vspace{2mm}\\
&\leq&\displaystyle \|x_n^{-\eta}\partial_k\int_0^1 G_\tau\ast[div\,(e_n\theta)]_\ast d\tau\|_{L^q(\mathbb{R}^{n-1}\times(0,1))}\vspace{2mm}\\
&&\displaystyle+\|x_n^{-\eta}\partial_k\int_1^\infty G_\tau\ast[div\,(e_n\theta)]_\ast d\tau\|_{L^q(\mathbb{R}^{n-1}\times(0,1))}+\|\partial_k\mathcal{N}div\,(e_n\theta)\|_{L^q(\mathbb{R}^n_+)}\vspace{2mm}\\
&\leq&\displaystyle C\sup\limits_{y\in\mathbb{R}^n}\|\int_0^1x_n^{-\eta}\partial_{x_k}G_\tau(x-y)d\tau\|_{L^1(\mathbb{R}^{n-1}\times(0,1))}
\|\partial_n\theta\|_{L^q(\mathbb{R}^n_+)}+\|\partial_k\mathcal{N}div\,(e_n\theta)\|_{L^q(\mathbb{R}^n_+)}\vspace{2mm}\\
\end{array}
\end{equation}
\begin{equation} \label{3.24}
\begin{array}{rcl}
&&\displaystyle+C\sup\limits_{y\in\mathbb{R}^n}\|\int_0^1\int_1^\infty x_n^{-\eta}\partial_{x_k}\partial_{x_n}\partial_{x_n}
G_\tau(x^\prime-y^\prime, x_n-sy_n)d\tau ds\|_{L^1(\mathbb{R}^{n-1}\times(0,1))}\|y_n\theta\|_{L^q(\mathbb{R}^n_+)}\vspace{2mm}\\
&=&\displaystyle CK_1+K_2+CK_3.
\end{array}
\end{equation}
First we estimate $K_1$ and $K_3$. Let $1\leq k\leq n$. $0\leq\eta<1$, $q_1, q_2\in(1, \infty)$ such that $\frac{1}{q_1}+\frac{1}{q_2}=1$ and $\eta q_1<1$. Then for any $y=(y^\prime, y_n)\in\mathbb{R}^n$
\begin{equation} \label{3.25}
\begin{array}{rcl}
&&\displaystyle\|\int_0^1
x_n^{-\eta}\partial_{x_k}G_\tau(x-y)d\tau\|_{L^1(\mathbb{R}^{n-1}\times(0,1))}\vspace{2mm}\\
&\leq&\displaystyle\int_0^1\int_{\mathbb{R}^{n-1}}\int_0^1\tau^{-\frac{1}{2}}x_n^{-\eta}\frac{|x_k-y_k|}{2\sqrt{\tau}}G_\tau(x-y)dxd\tau
\vspace{2mm}\\
&\leq&\displaystyle C\int_0^1\int_0^1\tau^{-1}x_n^{-\eta}e^{-\frac{(x_n-y_n)^2}{8\tau}}dx_nd\tau
\vspace{2mm}\\
&\leq&\displaystyle C\int_0^1\tau^{-1}\big(\int_0^1x_n^{-\eta q_1}dx_n\big)^\frac{1}{q_1}\big(\int_0^1e^{-\frac{(x_n-y_n)^2q_2}{8\tau}}dx_n\big)^\frac{1}{q_2}d\tau\vspace{2mm}\\
&\leq&\displaystyle C\int_0^1\tau^{-1+\frac{1}{2q_2}}d\tau\leq\displaystyle C;
\end{array}
\end{equation}
and
\begin{equation} \label{3.26}
\begin{array}{rcl}
&&\displaystyle\|\int_0^1\int_1^\infty
x_n^{-\eta}\partial_{x_k}\partial_{x_n}\partial_{x_n}G_\tau(x^\prime-y^\prime, x_n-sy_n)d\tau ds\|_{L^1(\mathbb{R}^{n-1}\times(0,1))}\vspace{2mm}\\
&\leq&\displaystyle C\int_0^1\int_1^\infty\tau^{-2}\int_{\mathbb{R}^{n-1}}\int_0^1x_n^{-\eta}\big(\frac{|x_k-\delta_k(s)y_k|}{\sqrt{\tau}}

+\frac{|x_n-sy_n|}{\sqrt{\tau}}\frac{|x_n-sy_n|}{\sqrt{\tau}}\frac{|x_k-\delta_k(s)y_k|}{\sqrt{\tau}}\big)\vspace{2mm}\\
&&\displaystyle\times(4\pi \tau)^{-\frac{n-1}{2}}e^{-\frac{|x^\prime-y^\prime|^2}{4\tau}}e^{-\frac{|x_n-sy_n|^2}{4\tau}}dx^\prime dx_nd\tau ds
\vspace{2mm}\\
&\leq&\displaystyle C\int_0^1\int_1^\infty\int_0^1\tau^{-2}x_n^{-\eta}e^{-\frac{(x_n-sy_n)^2}{8\tau}}dx_nd\tau ds
\vspace{2mm}\\
&\leq&\displaystyle C\int_0^1\int_1^\infty\tau^{-2}\big(\int_0^1x_n^{-\eta q_1}dx_n\big)^\frac{1}{q_1}\big(\int_0^1e^{-\frac{(x_n-sy_n)^2q_2}{8\tau}}dx_n\big)^\frac{1}{q_2}d\tau ds\vspace{2mm}\\
&\leq&\displaystyle C\int_1^\infty\tau^{-2+\frac{1}{2q_2}}d\tau\leq\displaystyle C,
\end{array}
\end{equation}
where
\[\delta_k(s) = \left\{
\begin{array}{l l}
s& \quad \text{if $k=n$}\\
1 & \quad \text{if $1\leq k\leq n-1$}
\end{array} \right.\]
Let us now estimate $K_2$. Let $1\leq q\leq\infty$. Using (\ref{3.18}) and (\ref{3.23}), we have for any $1\leq k\leq n$
\begin{equation} \label{3.27}
\begin{array}{rcl}
&&\displaystyle\|\partial_k\mathcal{N}div\,(e_n\theta)\|_{L^q(\mathbb{R}^n_+)}=\|\partial_k\int_0^\infty F(\tau)  div\,(e_n\theta)d\tau\|_{L^q(\mathbb{R}^n_+)}\vspace{2mm}\\
&\leq&\displaystyle \|\partial_k\int_0^1 G_\tau\ast[div\,(e_n\theta)]_\ast d\tau\|_{L^q(\mathbb{R}^n)}+\|\partial_k\int_1^\infty G_\tau\ast[div\,(e_n\theta)]_\ast d\tau\|_{L^q(\mathbb{R}^n)}\vspace{2mm}\\
&\leq&\displaystyle C\|\int_0^1
\partial_kG_\tau d\tau\|_{L^1(\mathbb{R}^n)}
\|\partial_n\theta\|_{L^q(\mathbb{R}^n_+)}\vspace{2mm}\\
&&\displaystyle+C\sup\limits_{(y^\prime, y_n)\in\mathbb{R}^n}\|\int_0^1\int_1^\infty\partial_{x_k}\partial_{x_n}\partial_{x_n}G_\tau(x^\prime-y^\prime, x_n-sy_n)d\tau ds\|_{L^1(\mathbb{R}^n)}\|y_n\theta\|_{L^q(\mathbb{R}^n_+)}\vspace{2mm}\\

&\leq&\displaystyle C\int_0^1\tau^{-\frac{1}{2}}d\tau
\|\partial_n\theta\|_{L^q(\mathbb{R}^n_+)}+C\int_1^\infty\tau^{-\frac{3}{2}} d\tau
\|y_n\theta\|_{L^q(\mathbb{R}^n_+)}\vspace{2mm}\\
&\leq&\displaystyle C(\|\partial_n\theta\|_{L^q(\mathbb{R}^n_+)}+\|y_n\theta\|_{L^q(\mathbb{R}^n_+)}).
\end{array}
\end{equation}
From (\ref{3.24})--(\ref{3.27}), we conclude that (\ref{3.22}) holds for $0\leq\eta<1$ and $1\leq k\leq n$.
$\,\,\,\,\,\,\Box$
To proceed, we need the following known results (see \cite{Han2} for (\ref{3.28}); \cite{Han1} for (\ref{3.29})):\\
Let $1<q<\infty$. Assume that $f=(f_1, f_2, \cdots, f_n)\in W^{1,q}(\mathbb{R}^n_+)\,\,(n\geq2)$ satisfies $\nabla\cdot f=0$ in $\mathbb{R}^n_+$ and $f_n|_{\partial \mathbb{R}^n_+}=0$. Then for any $\,0<\delta<1$ and $t>0$
\begin{equation} \label{3.28}
\|\nabla^2e^{-tA}f\|_{L^\infty(\mathbb{R}^n_+)}\leq
C(t^{-\frac{1}{2}-\frac{n}{2q}}\|\nabla f\|_{L^q(\mathbb{R}^n_+)}+t^{-1+\frac{\delta}{2}-\frac{n}{2q}}\|y_n^{-\delta}f\|_{L^q(\mathbb{R}^n_+)});
\end{equation}
and
\begin{equation} \label{3.29}
\|\sum\limits_{i,j=1}^n\nabla^2\mathcal{N}\partial_i\partial_j(u_iu_j)\|_{L^q(\mathbb{R}^n_+)}\leq
C(\|u\|^2_{L^{2q}(\mathbb{R}^n_+)}+\|\nabla
u\|^2_{L^{2q}(\mathbb{R}^n_+)}+\|\nabla^2u\|^2_{L^{2q}(\mathbb{R}^n_+)})
\end{equation}
for all $u\in C^\infty_{0,\sigma}(\mathbb{R}^n_+)$; and
\begin{equation} \label{3.30}
\|\nabla^2\mathcal{N}div\,(e_n\theta)\|_{L^q(\mathbb{R}^n_+)}\leq
C(\|\theta\|_{L^q(\mathbb{R}^n_+)}+\|\nabla^2\theta\|_{L^q(\mathbb{R}^n_+)})
\end{equation}
for any scalar function $\theta\in C^\infty_0(\mathbb{R}^n_+)$.\\
(\ref{3.28}) is the result of Lemma 3.4 in \cite{Han2}, and (\ref{3.29}) is from Lemma 2.3 in \cite{Han1}. The proof of (\ref{3.30}) is similar to that of (\ref{3.29}), here we give a sketch of its proof for readers' convenience.

Using (\ref{3.18}), we have for all $1\leq k,m\leq n$ and $t>0$,
$$
\begin{array}{rcl}
&&\displaystyle
\|\partial_k\partial_m\mathcal{N}div\;(e_n\theta)\|_{L^q(\mathbb{R}^n_+)}
=\|\partial_k\partial_m\int_0^\infty F(\tau)\partial_n\theta d\tau\|_{L^q(\mathbb{R}^n_+)}\vspace{2mm}\\
&=&\displaystyle\|\partial_k\partial_m\big(\int_0^1+\int_1^\infty\big)G_\tau\ast(\partial_n\theta)_\ast d\tau\|_{L^q(\mathbb{R}^n)}
\vspace{2mm}\\
&\leq&\displaystyle\|\partial_k\int_0^1G_\tau\ast[\partial_m(\partial_n\theta)_\ast]d\tau\|_{L^q(\mathbb{R}^n)}
+\|\partial_n\partial_k\partial_m\int_1^\infty
G_\tau\ast \theta^\ast d\tau\|_{L^q(\mathbb{R}^n)}\vspace{2mm}\\
&\leq&\displaystyle
C\int_0^1\|\partial_kG_\tau\|_{L^1(\mathbb{R}^n)}d\tau\|\nabla^2
\theta\|_{L^q(\mathbb{R}^n_+)}
+C\int_1^\infty\|\partial_n\partial_k\partial_mG_\tau\|_{L^1(\mathbb{R}^n)}d\tau\|\theta\|_{L^{q}(\mathbb{R}^n_+)}\vspace{2mm}\\
&\leq&\displaystyle
C\int_0^1\tau^{-\frac{1}{2}}d\tau\|\nabla^2\theta\|_{L^{q}(\mathbb{R}^n_+)}
+C\int_1^\infty\tau^{-\frac{3}{2}}d\tau\|\theta\|_{L^{q}(\mathbb{R}^n_+)}\vspace{2mm}\\
&\leq&\displaystyle C(\|\theta\|_{L^{q}(\mathbb{R}^n_+)}+\|\nabla^2
\theta\|_{L^{q}(\mathbb{R}^n_+)}\big),
\end{array}
$$
which is (\ref{3.30}). Here $f^\ast$, $f_\ast$ denote the odd and even extensions from $\mathbb{R}^n_+$ to $\mathbb{R}^n$, and their definitions are given in the proof of Lemma 3.2.\\
{\bf Proof of Theorem \ref{th:6}.}\quad The first part of the proof  establishes the decay of $\|\nabla^2u(t)\|_{L^\infty(\mathbb{R}^n_+)}$. We assume $0<\eta<1$, and $(u, \theta)$ is a strong solution of (\ref{1.1}). Then for $1<r<\infty$ and $t>0$
\begin{equation} \label{3.31}
\begin{array}{rcl}
&&\displaystyle\|x_n^{-\eta}u(t)\cdot\nabla u(t)\|_{L^r(\mathbb{R}^n_+)}+\|x_n^{-\eta}\theta(t)\|_{L^r(\mathbb{R}^n_+)}\vspace{2mm}\\
&\leq&\displaystyle \|x_n^{-\eta}u(t)\|_{L^{2r}(\mathbb{R}^{n-1}\times(0,1))}\|\nabla u(t)\|_{L^{2r}(\mathbb{R}^n_+)}+\|x_n^{-\eta}\theta(t)\|_{L^r(\mathbb{R}^{n-1}\times(0,1))}\vspace{2mm}\\
&&\displaystyle +\|u(t)\|_{L^{2r}(\mathbb{R}^n_+)}\|\nabla u(t)\|_{L^{2r}(\mathbb{R}^n_+)}+\|\theta(t)\|_{L^r(\mathbb{R}^n_+)}.
\end{array}
\end{equation}
Note that for  $t>0$, since  $u(x',0,t)=0$
\begin{equation} \label{3.32}
\begin{array}{rcl}
\displaystyle
\|x_n^{-\eta}u(t)\|^{2r}_{L^{2r}(\mathbb{R}^{n-1}\times(0,1))}
&\leq&\displaystyle\int_0^1\int_{\mathbb{R}^{n-1}}x_n^{-2\eta r}|u(x^\prime, x_n, t)-u(x^\prime, 0, t)|^{2r}dx^\prime dx_n\vspace{2mm}\\
&\leq&\displaystyle\int_0^1\int_{\mathbb{R}^{n-1}}x_n^{2r-1-2\eta r}\int_0^{x_n}|\partial_nu(x^\prime, z_n, t)|^{2r}dz_ndx^\prime dx_n\vspace{2mm}\\
&\leq&\displaystyle\int_{\mathbb{R}^{n-1}}\int_0^1|\partial_nu(x^\prime, z_n, t)|^{2r}dz_ndx^\prime\vspace{2mm}\\
&\leq&\displaystyle\|\nabla u(t)\|^{2r}_{L^{2r}(\mathbb{R}^n_+)}.
\end{array}
\end{equation}
Similarly,
\begin{equation} \label{3.33}
\|x_n^{-\eta}\theta(t)\|_{L^r(\mathbb{R}^{n-1}\times(0,1))}\leq\|\nabla\theta(t)\|_{L^r(\mathbb{R}^n_+)},\quad\forall\,t>0.
\end{equation}
Hence from (\ref{3.31})--(\ref{3.33}), we obtain for $t>0$.
\begin{equation} \label{3.34}
\begin{array}{rcl}
&&\displaystyle
\|x_n^{-\eta}u(t)\cdot\nabla u(t)\|_{L^r(\mathbb{R}^n_+)}+\|x_n^{-\eta}\theta(t)\|_{L^r(\mathbb{R}^n_+)}\vspace{2mm}\\
&\leq&\displaystyle C(\|u(t)\|^2_{L^{2r}(\mathbb{R}^n_+)}+\|\theta(t)\|_{L^r(\mathbb{R}^n_+)}+\|\nabla u(t)\|^2_{L^{2r}(\mathbb{R}^n_+)}+\|\nabla\theta(t)\|_{L^r(\mathbb{R}^n_+)}).
\end{array}
\end{equation}
To establish the decay $\|\nabla^2u(t)\|_{L^\infty(\mathbb{R}^n_+)}$, we use Lemma 3.2 combined with an estimate for the term $\|x_n\theta(t)\|_{L^r(\mathbb{R}^n_+)}$ with some $r\in(1,\infty)$. Since we don't have a local bound on $\|x_n\theta(t)\|_{L^r(\mathbb{R}^n_+)}$, we have to go back and consider the approximate solutions in (\ref{2.36}).\\
Using the definition of the operator of $E(t)$, we have for $t>0$
\begin{equation} \label{}\notag
\begin{array}{rcl}
\displaystyle
x_nE(t)b(x)&=&\displaystyle x_n\int_{\mathbb{R}^n_+}[G_t(x^\prime-y^\prime, x_n-y_n)-G_t(x^\prime-y^\prime, x_n+y_n)]b(y)dy\vspace{2mm}\\
&=&\displaystyle x_n\int_{\mathbb{R}^n_+}\int_{-1}^1\partial_{x_n}G_t(x^\prime-y^\prime, x_n-sy_n)(-y_n)b(y)dsdy\vspace{2mm}\\
&=&\displaystyle \int_{-1}^1\int_{\mathbb{R}^n_+}(x_n-sy_n)\partial_{x_n}G_t(x^\prime-y^\prime, x_n-sy_n)(-y_n)b(y)dsdyds\vspace{2mm}\\
&&\displaystyle+\int_{-1}^1\int_{\mathbb{R}^n_+}\partial_{x_n}G_t(x^\prime-y^\prime, x_n-sy_n)sy_n(-y_n)b(y)dyds,
\end{array}
\end{equation}
which together with the assumptions: $y_nb, y_n^2b\in L^1(\mathbb{R}^n_+)$, yield for $1<r<\infty$ and $t>0$
\begin{equation} \label{}\notag
\begin{array}{rcl}
\displaystyle
\|x_n\theta_0(t)\|_{L^r(\mathbb{R}^n_+)}&=&\displaystyle\|x_nE(t)b\|_{L^r(\mathbb{R}^n_+)}\vspace{2mm}\\
&\leq&\displaystyle C(\|x_n\partial_nG_t\|_{L^r(\mathbb{R}^n)}\|y_nb\|_{L^1(\mathbb{R}^n_+)}+\|\partial_n G_t\|_{L^r(\mathbb{R}^n)}\|y_n^2b\|_{L^1(\mathbb{R}^n_+)})\vspace{2mm}\\
&\leq&\displaystyle \overline{C_0}t^{-\frac{n}{2}(1-\frac{1}{r})}.
\end{array}
\end{equation}
Proceeding by induction, assume there exists $j\geq0$, such that
$$
\|x_n\theta_j(t)\|_{L^r(\mathbb{R}^n_+)}\leq C_{\ast\ast}t^{-\frac{n}{2}(1-\frac{1}{r})},
$$
where $C_{\ast\ast}\geq \overline{C_0}$ is independent of $j$.\\
Let $\frac{1}{r_0}=\frac{1}{n}+\frac{1}{r}$, $1<r<\infty$. Then from the second equation in (\ref{2.36}), we have for any $t>0$
\begin{equation} \label{3.35}
\begin{array}{rcl}
&&\displaystyle\|x_n\theta_{j+1}(t)\|_{L^r(\mathbb{R}^n_+)}\vspace{2mm}\\
&\leq&\displaystyle \|x_nE(t)b\|_{L^r(\mathbb{R}^n)}+C\int_0^\frac{t}{2}\||x_n|\nabla G_{t-s}\|_{L^r(\mathbb{R}^n)}\|u_j(s)\|_{L^2(\mathbb{R}^n_+)}\|\theta_j(s)\|_{L^2(\mathbb{R}^n_+)}ds\vspace{2mm}\\
&&\displaystyle+C\int_0^\frac{t}{2}\|\nabla G_{t-s}\|_{L^{r_0}(\mathbb{R}^n)}
\|u_j(s)\|_{L^{(1-\frac{1}{r_0})^{-1}}(\mathbb{R}^n_+)}\|y_n\theta_j(s)\|_{L^r(\mathbb{R}^n_+)}ds\vspace{2mm}\\
&&\displaystyle+C\int_\frac{t}{2}^t\||x_n|\nabla G_{t-s}\|_{L^1(\mathbb{R}^n)}\|u_j(s)\|_{L^\infty(\mathbb{R}^n_+)}\|\theta_j(s)\|_{L^r(\mathbb{R}^n_+)}ds\vspace{2mm}\\
&&\displaystyle+C\int_\frac{t}{2}^t\|\nabla G_{t-s}\|_{L^1(\mathbb{R}^n)}\|u_j(s)\|_{L^\infty(\mathbb{R}^n_+)}\|y_n\theta_j(s)\|_{L^r(\mathbb{R}^n_+)}ds\vspace{2mm}\\
&=&\displaystyle I_1(t)+I_2(t)+I_3(t)+I_4(t).
\end{array}
\end{equation}
Using Theorems \ref{th:4}, \ref{th:5}, where the decay estimates are also true for $(u_j, \theta_j)$, we get for $t>0$
\begin{equation} \label{3.36}
\begin{array}{rcl}
\displaystyle
I_1(t)&\leq&\displaystyle Ct^{-\frac{n}{2}(1-\frac{1}{r})}+ C\int_0^\frac{t}{2}(t-s)^{-\frac{n}{2}(1-\frac{1}{r})}(1+s)^{\frac{1}{2}-\frac{n}{4}-\frac{n+2}{4}}ds\vspace{2mm}\\
&\leq&\displaystyle Ct^{-\frac{n}{2}(1-\frac{1}{r})}+Ct^{-\frac{n}{2}(1-\frac{1}{r})}\int_0^\infty(1+s)^{-\frac{n}{2}}ds\vspace{2mm}\\
&\leq&\displaystyle Ct^{-\frac{n}{2}(1-\frac{1}{r})};
\end{array}
\end{equation}
\begin{equation} \label{3.37}
I_3(t)\leq\int^t_\frac{t}{2}s^{-\frac{1}{2}-\frac{1}{2}-\frac{n}{2}(1-\frac{1}{r})}ds\leq Ct^{-\frac{n}{2}(1-\frac{1}{r})};
\end{equation}
and
\begin{equation} \label{3.38}
I_4(t)\leq C\epsilon_0C_{\ast\ast} \int_\frac{t}{2}^t(t-s)^{-\frac{1}{2}}s^{-\frac{1}{2}-\frac{n}{2}(1-\frac{1}{r})}ds\leq C\epsilon_0C_{\ast\ast}t^{-\frac{n}{2}(1-\frac{1}{r})}.
\end{equation}
Now we estimate $I_2(t)$. Let $(\frac{n}{n-1}<)\frac{2n}{n-2}<r<\infty$. By the choice of $\frac{1}{r_0}=\frac{1}{n}+\frac{1}{r}$, we find $2<r_0<n$ and $\frac{n}{n-1}<(1-\frac{1}{r_0})^{-1}<2$. Note that $y_nb\in L^1(\mathbb{R}^n_+)\bigcap L^\infty(\mathbb{R}^n_+)$.
Whence for any $0<\gamma<(1-\frac{1}{r_0})^{-1}-\frac{n}{n-1}$ and $s>0$
\begin{equation} \label{3.39}
\begin{array}{rcl}
&&\displaystyle\|u_j(s)\|_{L^{(1-\frac{1}{r_0})^{-1}}(\mathbb{R}^n_+)}\|y_n\theta_j(s)\|_{L^r(\mathbb{R}^n_+)}\vspace{2mm}\\
&\leq&\displaystyle  C_{\ast\ast}(1+s)^{-\frac{n}{2}(1-\frac{1}{r})}\|u_j(s)\|^{\gamma(1-\frac{1}{r_0})}_{L^\infty(\mathbb{R}^n_+)}
\|u_j(s)\|^{1-\gamma(1-\frac{1}{r_0})}_{L^{(1-\frac{1}{r_0})^{-1}-\gamma}(\mathbb{R}^n_+)}\vspace{2mm}\\
&\leq&\displaystyle C\epsilon_0^{\gamma(1-\frac{1}{r_0})}C_{\ast\ast}s^{-\frac{\gamma}{2}(1-\frac{1}{r_0})}
(1+s)^{-\frac{n}{2}(1-\frac{1}{r})+[\frac{1}{2}-\frac{n}{2}(1-\frac{1}{(1-\frac{1}{r_0})^{-1}-\gamma})](1-\gamma(1-\frac{1}{r_0}))}
\vspace{2mm}\\
&\leq&\displaystyle C\epsilon_0^{\gamma(1-\frac{1}{r_0})}C_{\ast\ast}(1+s)^{\mu(\gamma)}s^{-\frac{\gamma}{2}(1-\frac{1}{r_0})},
\end{array}
\end{equation}
where $\mu(\gamma)=-\frac{n}{2}(1-\frac{1}{r})+\frac{1}{2}-\frac{n}{2}(1-\frac{1}{(1-\frac{1}{r_0})^{-1}-\gamma})
+\frac{n\gamma}{2}(1-\frac{1}{r_0})(1-\frac{1}{(1-\frac{1}{r_0})^{-1}-\gamma})$.
Note that
\begin{equation} \label{}\notag
\lim\limits_{\gamma\longrightarrow0}\mu(\gamma)=-\frac{n}{2}(1-\frac{1}{r})+\frac{1}{2}-\frac{n}{2r_0}
=-\frac{n}{2}(1-\frac{1}{r})+\frac{1}{2}-\frac{n}{2}(\frac{1}{n}+\frac{1}{r})=-\frac{n}{2}.
\end{equation}
Thus there exists $\gamma\in(0, (1-\frac{1}{r_0})^{-1}-\frac{n}{n-1})$ such that $1+\mu(\gamma)<0$ and $1-\frac{\gamma}{2}(1-\frac{1}{r_0})>0$.\\
Let $t>2$. Then
\begin{equation} \label{}\notag
\int_0^\frac{t}{2}(1+s)^{\mu(\gamma)}s^{-\frac{\gamma}{2}(1-\frac{1}{r_0})}ds
\leq\int_0^1s^{-\frac{\gamma}{2}(1-\frac{1}{r_0})}ds+\int_1^\frac{t}{2}(1+s)^{\mu(\gamma)}ds\leq C.
\end{equation}
If $0<t\leq2$, then
\begin{equation} \label{}\notag
\int_0^\frac{t}{2}(1+s)^{\mu(\gamma)}s^{-\frac{\gamma}{2}(1-\frac{1}{r_0})}ds
\leq\int_0^1s^{-\frac{\gamma}{2}(1-\frac{1}{r_0})}ds\leq C.
\end{equation}
Whence for $t>0$
\begin{equation} \label{}\notag
\int_0^\frac{t}{2}(1+s)^{\mu(\gamma)}s^{-\frac{\gamma}{2}(1-\frac{1}{r_0})}ds\leq C.
\end{equation}
Therefore from (\ref{3.39}), we derive for any $t>0$
\begin{equation} \label{3.40}
\begin{array}{rcl}
I_2(t)&\leq&\displaystyle  C\epsilon_0^{\gamma(1-\frac{1}{r_0})}C_{\ast\ast}\int_0^\frac{t}{2}(t-s)^{-\frac{1}{2}-\frac{n}{2}(1-\frac{1}{r_0})}
(1+s)^{\mu(\gamma)}s^{-\frac{\gamma}{2}(1-\frac{1}{r_0})}ds\vspace{2mm}\\
&\leq&\displaystyle C\epsilon_0^{\gamma(1-\frac{1}{n}-\frac{1}{r})}C_{\ast\ast}t^{-\frac{n}{2}(1-\frac{1}{r})}.
\end{array}
\end{equation}
From (\ref{3.35})--(\ref{3.38}) and (\ref{3.40}), we obtain for the given $j$ and any $t>0$
\begin{equation} \label{3.41}
\|x_n\theta_{j+1}(t)\|_{L^r(\mathbb{R}^n_+)}\leq \overline{C_1}t^{-\frac{n}{2}(1-\frac{1}{r})}+\overline{C_1}(\epsilon_0+\epsilon_0^{\gamma(1-\frac{1}{n}-\frac{1}{r})})C_{\ast\ast}t^{-\frac{n}{2}(1-\frac{1}{r})}.
\end{equation}
Take $\epsilon_0>0$ suitably small so that $\overline{C_1}(\epsilon_0+\epsilon_0^{\gamma(1-\frac{1}{n}-\frac{1}{r})})\leq\frac{1}{2}$ in (\ref{3.41}), and  $C_{\ast\ast}=\overline{C_0}+2\overline{C_1}$. Then for $\frac{2n}{n-2}<r<\infty$ and any $t>0$
\begin{equation} \label{}\notag
\|x_n\theta_{j+1}(t)\|_{L^r(\mathbb{R}^n_+)}\leq C_{\ast\ast} t^{-\frac{n}{2}(1-\frac{1}{r})}.
\end{equation}
From the above induction argument, we conclude that for $\frac{2n}{n-2}<r<\infty$ and any $t>0$
\begin{equation} \label{}\notag
\sup\limits_{j\geq0}\|x_n\theta_{j+1}(t)\|_{L^r(\mathbb{R}^n_+)}\leq C_{\ast\ast} t^{-\frac{n}{2}(1-\frac{1}{r})}.
\end{equation}
By a standard weak convergence procedure, we conclude for $\frac{2n}{n-2}<r<\infty$ and any $t>0$
\begin{equation} \label{3.42}
\|x_n\theta(t)\|_{L^r(\mathbb{R}^n_+)}\leq C_{\ast\ast} t^{-\frac{n}{2}(1-\frac{1}{r})}.
\end{equation}
From (\ref{3.28})--(\ref{3.30}), (\ref{3.34}), (\ref{3.42}) and Lemmata 2.1, 3.1, 3.2, Theorem \ref{th:5}, we deduce for any $t\geq 1$
\begin{equation} \label{}\notag
\begin{array}{rcl}
&&\displaystyle\|\nabla^2 u(t)\|_{L^\infty(\mathbb{R}^n_+)}\vspace{2mm}\\
&\leq&\displaystyle\|\nabla^2e^{-\frac{t}{2}A}u(\frac{t}{2})\|_{L^\infty(\mathbb{R}^n_+)}+\int_\frac{t}{2}^t\|\nabla^2e^{-(t-s)A}(Pu(s)\cdot\nabla u(s)-Pe_n\theta(s))\|_{L^\infty(\mathbb{R}^n_+)}ds\vspace{2mm}\\
&\leq&\displaystyle Ct^{-1-\frac{n}{2r}}\|u(\frac{t}{2})\|_{L^r(\mathbb{R}^n_+)}
+C\int_\frac{t}{2}^t(t-s)^{-\frac{1}{2}-\frac{n}{2r}}\|\nabla (Pu(s)\cdot\nabla u(s)-Pe_n\theta(s))\|_{L^r(\mathbb{R}^n_+)}ds\vspace{2mm}\\
&&\displaystyle+C\int_\frac{t}{2}^t(t-s)^{-1+\frac{\eta}{2}-\frac{n}{2r}}\|x_n^{-\eta}(Pu(s)\cdot\nabla u(s)-Pe_n\theta(s))\|_{L^r(\mathbb{R}^n_+)}ds\vspace{2mm}\\
&\leq&\displaystyle Ct^{-1-\frac{n}{2r}}\|u(\frac{t}{2})\|_{L^r(\mathbb{R}^n_+)}
+C\int_\frac{t}{2}^t(t-s)^{-\frac{1}{2}-\frac{n}{2r}}(\|\nabla(u(s)\cdot\nabla u(s))\|_{L^r(\mathbb{R}^n_+)}+\|\nabla\theta(s)\|_{L^r(\mathbb{R}^n_+)}\vspace{2mm}\\
&&\displaystyle
+\|\nabla(\sum\limits_{i,j=1}^n\nabla\mathcal{N}\partial_i\partial_j(u_iu_j)(s))\|_{L^r(\mathbb{R}^n_+)}
+\|\nabla^2\mathcal{N}div\,(e_n\theta(s))\|_{L^r(\mathbb{R}^n_+)})ds
\vspace{2mm}\\
&&\displaystyle+C\int_\frac{t}{2}^t(t-s)^{-1+\frac{\eta}{2}-\frac{n}{2r}}(\|x_n^{-\eta}u(s)\cdot\nabla u(s)\|_{L^r(\mathbb{R}^n_+)}
+\|x_n^{-\eta}\theta(s)\|_{L^r(\mathbb{R}^n_+)}\vspace{2mm}\\
&&\displaystyle+\|\sum\limits_{i,j=1}^nx_n^{-\eta}\nabla\mathcal{N}\partial_i\partial_j(u_iu_j)(s)\|_{L^r(\mathbb{R}^n_+)}
+\|x_n^{-\eta}\nabla\mathcal{N}div\,(e_n\theta(s))\|_{L^r(\mathbb{R}^n_+)})ds\vspace{2mm}\\
&\leq&\displaystyle Ct^{-1-\frac{n}{2r}}\|u(\frac{t}{2})\|_{L^r(\mathbb{R}^n_+)}\vspace{2mm}\\
&&\displaystyle+C\int_\frac{t}{2}^t(t-s)^{-\frac{1}{2}-\frac{n}{2r}}(\|u(s)\|^2_{L^{2r}(\mathbb{R}^n_+)}+\|\nabla u(s)\|^2_{L^{2r}(\mathbb{R}^n_+)}+\|\nabla^2 u(s)\|^2_{L^{2r}(\mathbb{R}^n_+)}\vspace{2mm}\\
&&\displaystyle+\|\theta(s)\|_{L^{r}(\mathbb{R}^n_+)}+\|\nabla \theta(s)\|_{L^{r}(\mathbb{R}^n_+)}+\|\nabla^2 \theta(s)\|_{L^{r}(\mathbb{R}^n_+)})ds\vspace{2mm}\\
&&\displaystyle+C\int_\frac{t}{2}^t(t-s)^{-1+\frac{\eta}{2}-\frac{n}{2r}}(\|u(s)\|^2_{L^{2r}(\mathbb{R}^n_+)}+\|\nabla u(s)\|^2_{L^{2r}(\mathbb{R}^n_+)}\vspace{2mm}\\
&&\displaystyle+\|y_n\theta(s)\|_{L^{r}(\mathbb{R}^n_+)}+\|\theta(s)\|_{L^{r}(\mathbb{R}^n_+)}
+\|\nabla\theta(s)\|_{L^{r}(\mathbb{R}^n_+)})ds\vspace{2mm}\\
&\leq&\displaystyle Ct^{-\frac{n+1}{2}}+C\int_\frac{t}{2}^t(t-s)^{-\frac{1}{2}-\frac{n}{2r}}(s^{1-n(1-\frac{1}{2r})}
+s^{-\frac{1}{2}-\frac{n}{2}(1-\frac{1}{r})})ds\vspace{2mm}\\
&&\displaystyle+C\int_\frac{t}{2}^t(t-s)^{-1+\frac{\eta}{2}-\frac{n}{2r}}(s^{1-n(1-\frac{1}{2r})}
+s^{-\frac{n}{2}(1-\frac{1}{r})})ds\vspace{2mm}\\
\vspace{2mm}\\
&\leq&\displaystyle Ct^{-\frac{n+1}{2}}+Ct^{-\frac{n}{2}}(1+t^{-\frac{n-3}{2}})+Ct^{-\frac{n}{2}+\frac{\eta}{2}}(1+t^{-\frac{n-2}{2}})\vspace{2mm}\\
\vspace{2mm}\\
&\leq&\displaystyle Ct^{-\frac{n}{2}+\frac{\eta}{2}}\;\; require \;\;\max\{\frac{n}{\eta},\; \frac{2n}{n-2}\}<r<\infty,
\end{array}
\end{equation}
from which, we give the decay of $\|\nabla^2 u(t)\|_{L^\infty(\mathbb{R}^n_+)}$.\\
Now we use (\ref{2.36}) by induction to establish the decay of $\|\nabla^3\theta(t)\|_{L^r(\mathbb{R}^n_+)}$.\\
Let $1\leq q\leq\infty$. Using the structure of the operator of $E(t)$, it is not difficult to verify that for any $t>0$
\begin{equation} \label{3.43}
\begin{array}{rcl}
\displaystyle
\|\nabla^3\theta_0(t)\|_{L^q(\mathbb{R}^n_+)}&=&\displaystyle\|\nabla^3E(t)b\|_{L^q(\mathbb{R}^n_+)}\vspace{2mm}\\
&=&\displaystyle\|\nabla^3_x\partial_{x_n}\int_{\mathbb{R}^n_+}\int_{-1}^1G_t(x^\prime-y^\prime, x_n-sy_n)y_nb(y)dsdy\|_{L^q(\mathbb{R}^n_+)}\vspace{2mm}\\
&\leq&\displaystyle C\|\nabla^4G_t\|_{L^q(\mathbb{R}^n)}\|y_nb\|_{L^1(\mathbb{R}^n_+)}\vspace{2mm}\\
&\leq&\displaystyle C_2t^{-2-\frac{n}{2}(1-\frac{1}{q})}.
\end{array}
\end{equation}
Assume that there exists $j\geq0$, such that for any $t>0$
\begin{equation} \label{3.44}
\|\nabla^3\theta_j(t)\|_{L^q(\mathbb{R}^n_+)}\leq C_{\ast\ast\ast} t^{-2-\frac{n}{2}(1-\frac{1}{q})}.
\end{equation}
where $C_{\ast\ast\ast}\geq C_2$ is independent of $j$, which will be determined later.\\
To proceed, we need to find some properties of the operators of $E(t)$ and $F(t)$.
\begin{equation} \label{3.45}
\partial^\ell_m[E(t)g]=E(t)\partial^\ell_m g\quad and\quad \partial^\ell_m[F(t)g]=F(t)\partial^\ell_m g,\;\;\;\forall\;1\leq m\leq n-1,\;\;\ell=1, 2,\cdots;
\end{equation}
\begin{equation} \label{3.46}
\partial_n[E(t)f]=F(t)\partial_nf\;\;\;for\;\;\;any\;\;\;f\;|_{\partial\mathbb{R}^n_+}=0;
\end{equation}
\begin{equation} \label{3.47}
\partial_n[F(t)g]=E(t)\partial_n g.
\end{equation}
(\ref{3.45}) is obvious, because no boundary arises for $x_j\in\mathbb{R}^1$ for $1\leq j\leq n-1$. Inequalities (\ref{3.46}) and (\ref{3.47}) are obtained as follows.
$$
\begin{array}{rcl}
\displaystyle\partial_{x_n}[E(t)f](x)&=&\displaystyle \int_{\mathbb{R}^n_+}\{\partial_{x_n}[G_t(x^\prime-y^\prime,
x_n-y_n)-G_t(x^\prime-y^\prime, x_n+y_n)]\}f(y)dy\vspace{2mm}\\
&=&\displaystyle \int_{\mathbb{R}^n_+}\{-\partial_{y_n}[G_t(x^\prime-y^\prime,
x_n-y_n)+G_t(x^\prime-y^\prime, x_n+y_n)]\}f(y)dy\vspace{2mm}\\
&=&\displaystyle \int_{\mathbb{R}^n_+}[G_t(x^\prime-y^\prime,
x_n-y_n)+G_t(x^\prime-y^\prime, x_n+y_n)]\partial_{y_n}f(y)dy\vspace{2mm}\\
&=&\displaystyle F(t)(\partial_nf)(x),
\end{array}
$$
yielding (\ref{3.46}). Here we used that  $f(y)=0\;\;on\;\;\partial\mathbb{R}^n_+$. To obtain  (\ref{3.47}) we proceed as follows,
$$
\begin{array}{rcl}
\displaystyle\partial_{x_n}[F(t)g](x)&=&\displaystyle \int_{\mathbb{R}^n_+}\{\partial_{x_n}[G_t(x^\prime-y^\prime,
x_n-y_n)+G_t(x^\prime-y^\prime, x_n+y_n)]\}g(y)dy\vspace{2mm}\\
&=&\displaystyle \int_{\mathbb{R}^n_+}\{-\partial_{y_n}[G_t(x^\prime-y^\prime,
x_n-y_n)-G_t(x^\prime-y^\prime, x_n+y_n)]\}g(y)dy\vspace{2mm}\\
&=&\displaystyle \int_{\mathbb{R}^n_+}[G_t(x^\prime-y^\prime,
x_n-y_n)-G_t(x^\prime-y^\prime, x_n+y_n)]\partial_{y_n}g(y)dy\vspace{2mm}\\
&=&\displaystyle E(t)(\partial_ng)(x).
\end{array}
$$
Here we used that: $[G_t(x^\prime-y^\prime,
x_n-y_n)-G_t(x^\prime-y^\prime, x_n+y_n)]\,|_{y\in\partial\mathbb{R}^n_+}=0$.\\
Recall that for any $j\geq0$ and $t>0$
\begin{equation} \label{}\notag
\theta_{j+1}(t)=E(\frac{t}{2})\theta_{j+1}(\frac{t}{2})-\int_\frac{t}{2}^tE(t-s)u_j(s)\cdot\nabla\theta_j(s)ds.
\end{equation}
It follows from (\ref{2.43})-- (\ref{2.47}), and Theorems \ref{th:4}, \ref{th:5} (where the decay estimates are also valid for $(u_j, \theta_j)$) that for $1\leq q\leq\infty$ and $t\geq1$
\begin{equation} \label{3.48}
\begin{array}{rcl}
&&\displaystyle
\|\nabla^3\theta_{j+1}(t)\|_{L^q(\mathbb{R}^n_+)}\vspace{2mm}\\
&\leq&\displaystyle\|\nabla^3E(\frac{t}{2})\theta_{j+1}(\frac{t}{2})\|_{L^q(\mathbb{R}^n_+)}+2\int_\frac{t}{2}^t\|\nabla G_{t-s}\|_{L^1(\mathbb{R}^n)}\|\nabla^2\big(u_j(s)\cdot\nabla \theta_j(s)\big)\|_{L^q(\mathbb{R}^n_+)}ds\vspace{2mm}\\
&\leq&\displaystyle Ct^{-\frac{3}{2}}\|\theta_{j+1}(\frac{t}{2})\|_{L^q(\mathbb{R}^n_+)}+C\int_\frac{t}{2}^t(t-s)^{-\frac{1}{2}}\|\nabla^2 u_j(s)\|_{L^{2q}(\mathbb{R}^n_+)}\|\nabla \theta_j(s)\|_{L^{2q}(\mathbb{R}^n_+)}ds\vspace{2mm}\\
&&\displaystyle+C\int_\frac{t}{2}^t(t-s)^{-\frac{1}{2}}\|\nabla u_j(s)\|_{L^\infty(\mathbb{R}^n_+)}\|\nabla^2\theta_j(s)\|_{L^q(\mathbb{R}^n_+)}ds
\vspace{2mm}\\
&&\displaystyle+C\int_\frac{t}{2}^t(t-s)^{-\frac{1}{2}}\|u_j(s)\|_{L^\infty(\mathbb{R}^n_+)}\|\nabla^3\theta_j(s)\|_{L^q(\mathbb{R}^n_+)}ds
\vspace{2mm}\\
&\leq&\displaystyle Ct^{-2-\frac{n}{2}(1-\frac{1}{q})}+C\int_\frac{t}{2}^t(t-s)^{-\frac{1}{2}}
s^{-1-\frac{n}{2}(1-\frac{1}{2q})}Q(q, s)ds\vspace{2mm}\\
&&\displaystyle+C\int_\frac{t}{2}^t(t-s)^{-\frac{1}{2}}
s^{-\frac{n}{2}-\frac{3}{2}-\frac{n}{2}(1-\frac{1}{q})}ds\vspace{2mm}\\
&&\displaystyle+C\epsilon_0C_{\ast\ast\ast}\int_\frac{t}{2}^t(t-s)^{-\frac{1}{2}}s^{-\frac{1}{2}-2-\frac{n}{2}(1-\frac{1}{q})}ds
\vspace{2mm}\\
&\leq&\displaystyle Ct^{-2-\frac{n}{2}(1-\frac{1}{q})}+Ct^{-1-\frac{n}{2}(1-\frac{1}{2q})}Q(q, t)+Ct^{-\frac{n+2}{2}-\frac{n}{2}(1-\frac{1}{q})}
+C\epsilon_0C_{\ast\ast\ast} t^{-2-\frac{n}{2}(1-\frac{1}{q})},
\end{array}
\end{equation}
where $Q(q, t)=\left\{
\begin{array}{lll}
t^{-\frac{1}{2}-\frac{n}{2}(1-\frac{1}{2q})} &if& 1\leq q<\infty,\\
t^{-\frac{n}{2}+\epsilon}  &if& q=\infty.
\end{array}
\right.$ The number $\epsilon$ comes from the decay of $\|\nabla^2 u(t)\|_{L^\infty(\mathbb{R}^n_+)}$ in Theorem \ref{th:6}.\\
Let $1\leq q<\infty$. Then from (\ref{3.48}), we get for any $t\geq1$
\begin{equation} \label{3.49}
\|\nabla^3\theta_{j+1}(t)\|_{L^q(\mathbb{R}^n_+)}\leq C_3t^{-2-\frac{n}{2}(1-\frac{1}{q})}+C\epsilon_0C_{\ast\ast\ast} t^{-2-\frac{n}{2}(1-\frac{1}{q})}
\end{equation}
Taking $\epsilon_0>0$ suitably small such that $C\epsilon_0\leq\frac{1}{2}$ in (\ref{3.49}), and take $C_{\ast\ast\ast}=C_2+2C_3$. Then for $1\leq q<\infty$ and $t\geq1$
\begin{equation} \label{3.50}
\|\nabla^3\theta_{j+1}(t)\|_{L^q(\mathbb{R}^n_+)}\leq C_{\ast\ast\ast} t^{-2-\frac{n}{2}(1-\frac{1}{q})}.
\end{equation}
From (\ref{3.43}), (\ref{3.44}) and (\ref{3.50}), we conclude that for $1\leq q<\infty$ and $t\geq1$
\begin{equation} \label{}\notag
\sup\limits_{j\geq0}\|\nabla^3\theta_j(t)\|_{L^q(\mathbb{R}^n_+)}\leq C_{\ast\ast\ast} t^{-2-\frac{n}{2}(1-\frac{1}{q})}.
\end{equation}
By a standard weak convergence procedure, we get for $1\leq q<\infty$ and $t\geq1$
\begin{equation} \label{3.51}
\|\nabla^3\theta(t)\|_{L^q(\mathbb{R}^n_+)}\leq C_{\ast\ast\ast} t^{-2-\frac{n}{2}(1-\frac{1}{q})}.
\end{equation}
If $q=\infty$. From (\ref{3.48}), we have for any $t\geq1$
\begin{equation} \label{3.52}
\|\nabla^3\theta_{j+1}(t)\|_{L^\infty(\mathbb{R}^n_+)}\leq Ct^{-2-\frac{n}{2}}+Ct^{-2-\frac{n}{2}+\epsilon-\frac{n-2}{2}}+C\epsilon_0C_{\ast\ast\ast} t^{-2-\frac{n}{2}}
\end{equation}
Take $0<\epsilon<\frac{n-2}{2}$ in (\ref{3.52}). Repeating the above arguments, we readily find there exists a constant $\widetilde{C}_{\ast\ast\ast}>0$ such that for any $t\geq1$
\begin{equation} \label{3.53}
\|\nabla^3\theta(t)\|_{L^\infty(\mathbb{R}^n_+)}\leq \widetilde{C}_{\ast\ast\ast} t^{-2-\frac{n}{2}}.
\end{equation}
Combining (\ref{3.51}) and (\ref{3.53}), we establish the decay of $\|\nabla^3\theta(t)\|_{L^q(\mathbb{R}^n_+)}$ with $1\leq q\leq \infty$. From the above arguments, and together with Lemma \ref{l:3.1}, we complete the proof of the theorem. $\,\,\,\,\,\,\Box$

\vskip 2mm

{\bf Acknowledgements:} This work was completed with the support by NSFC: 11071239 for P. Han.  M. Schonbek was partially
supported by  NSF Grant DMS-0900909.

\end{document}